\documentclass[11pt]{amsart}
\usepackage[hmargin=3cm,vmargin=3cm]{geometry}
\usepackage[backend=bibtex,doi=false,isbn=false,style=alphabetic,url=false,maxnames=6]{biblatex}
\usepackage{etex}
\bibliography{../../../bibliography/bib}

\usepackage[latin1]{inputenc}
\usepackage{xspace,amssymb,amsfonts}
\usepackage{amsthm,amsmath,amscd,stmaryrd,latexsym,youngtab,mathtools,comment}
\usepackage[scaled]{helvet} % ss
\usepackage{courier} % tt

\usepackage{combelow}

\usepackage{euscript}

\usepackage{graphicx}
\usepackage[all]{xy}
\SelectTips{cm}{}

\usepackage{tikz, tikz-cd}
\usetikzlibrary{matrix, arrows, calc}
\tikzset{frontline/.style={preaction={draw=white,-,line width=6pt}},}  %%% for 3d commutative diagrams

\usepackage[dvips]{epsfig}
\usepackage{pinlabel}
\usepackage{psfrag}
\tikzset{ssmtf/.style={matrix of math nodes, row sep=3em,
               column sep=3em, text height=2.5ex, text depth=1.25ex} }

\usepackage[nohug]{diagrams}

\newcommand{\idempotent}{\EB}
\newcommand{\unital}{\UB}
\newcommand{\counital}{\CB}
\newcommand{\HHH}{\operatorname{HHH}}
\newcommand{\HH}{\operatorname{HH}}
\newcommand{\Hecke}{\HB}

\newcommand{\inv}{^{-1}}
\newcommand{\R}{\mathbb{R}}
\newcommand{\Z}{\mathbb{Z}}
\newcommand{\ring}{\Z}

\newcommand{\Q}{\mathbb{Q}}
\newcommand{\C}{\mathbb{C}}
\renewcommand{\d}{\delta}
\renewcommand{\emptyset}{\varnothing}
\newcommand{\Endg}{\underline{\operatorname{End}}}

\newcommand{\DHom}{\operatorname{DHom}}
\newcommand{\DEnd}{\operatorname{DEnd}}

\newcommand{\Homg}{\underline{\operatorname{Hom}}}

\newcommand{\FT}{\operatorname{FT}}

\newcommand{\Id}{\operatorname{Id}}
\newcommand{\Hom}{\operatorname{Hom}}

\newcommand{\ip}[1]{\langle #1 \rangle}
\renewcommand{\matrix}[1]{\begin{bmatrix}#1\end{bmatrix}}
\newcommand{\one}{\mathbbm{1}}
\newcommand{\K}{\mathcal{K}}

\newcommand{\Sym}{\operatorname{Sym}}

\newcommand{\xx}{\mathbf{x}}
\newcommand{\yy}{\mathbf{y}}

\renewcommand{\deg}{\operatorname{deg}}

\newcommand{\im}{\operatorname{im}}
\newcommand{\Tot}{\operatorname{Tot}}
\newcommand{\TL}{\operatorname{TL}}
\newcommand{\Cone}{\operatorname{Cone}}
\newcommand{\spic}[1]{\begin{minipage}{6pt}\includegraphics[scale=1]{fig/#1}\end{minipage}}
\newcommand{\sspic}[2]{\begin{minipage}{#1}\includegraphics[scale=1]{fig/#2}\end{minipage}}

\IfFileExists{comments.sty}{\usepackage{comments}}

\RequirePackage{color}
\definecolor{myred}{rgb}{0.75,0,0}
\definecolor{mygreen}{rgb}{0,0.5,0}
\definecolor{myblue}{rgb}{0,0.25,0.65}
\definecolor{references}{rgb}{0,0,1}
\RequirePackage{ifpdf}
\ifpdf
  \IfFileExists{pdfsync.sty}{\RequirePackage{pdfsync}}{}
  \RequirePackage[pdftex,
   colorlinks = true,
   urlcolor = references, % \href{...}{...} external (URL)
   citecolor = references, % \cite{...}
   linkcolor = references, % \ref{...} and \pageref{...}
   bookmarksopen=true]{hyperref}
\else
  \RequirePackage[hypertex]{hyperref}
\fi

\newtheorem{theorem}{Theorem}[section]
\newtheorem{lemma}[theorem]{Lemma}
\newtheorem{deth}[theorem]{Definition/Theorem}
\newtheorem{proposition}[theorem]{Proposition}

\newtheorem{corollary}[theorem]{Corollary}

\newtheorem{conjecture}[theorem]{Conjecture}

\theoremstyle{definition}

\newtheorem{definition}[theorem]{Definition}
\newtheorem{construction}[theorem]{Construction}
\newtheorem{notation}[theorem]{Notation}

\newtheorem{example}[theorem]{Example}

\newtheorem{hypothesis}[theorem]{Hypothesis}

\theoremstyle{remark}
\newtheorem{remark}[theorem]{Remark}

\numberwithin{equation}{section}

    \def\SM{{\mathbb{S}}}

    \def\AC{{\mathcal{A}}}
    
\def\CB{{\mathbf C}}  %\def\cb{{\mathbf c}} 
\def\CC{{\mathcal{C}}}
    \def\DC{{\mathcal{D}}}
\def\EB{{\mathbf E}}

\def\HB{{\mathbf H}}    
\def\IB{{\mathbf I}}    
    
    \def\KC{{\mathcal{K}}}

    \def\NC{{\mathcal{N}}}
    
\def\PB{{\mathbf P}}    \def\PC{{\mathcal{P}}}

\def\UB{{\mathbf U}}

\def\AS{{\EuScript A}}
\def\CS{{\EuScript C}}
\def\DS{{\EuScript D}}

\def\PS{{\EuScript P}}

\def\a{\alpha}
\def\b{\beta}

\def\d{\delta}

\def\e{\varepsilon}

\def\l{\lambda}

\let\phi=\varphi

\usepackage{bbm}
\def\C{{\mathbbm C}}

\def\R{{\mathbbm R}}
\def\Z{{\mathbbm Z}}
\def\Q{{\mathbbm Q}}
\def\1{\mathbbm{1}}

\newcommand{\rank}{\operatorname{rank}}

\renewcommand{\to}{\rightarrow}

\renewcommand{\sl}{\mathfrak{sl}}

\renewcommand{\mod}{\mathrm{-mod}}

\def\gmod{{\text{-gmod}}}

\newcommand{\refequal}[1]{\xy {\ar@{=}^{#1}
(-1,0)*{};(1,0)*{}};
\endxy}

\newcommand{\End}{\operatorname{End}}
\newcommand{\Tr}{\operatorname{Tr}}

\newcommand{\Ext}{\operatorname{Ext}}

\newcommand{\Ch}{\textrm{Ch}}

\renewcommand{\ker}{\operatorname{ker}}

\newcommand{\hocolim}{\operatorname{hocolim}}
\DeclarePairedDelimiter\floor{\lfloor}{\rfloor}

\newcommand{\SBim}{\SM\textrm{Bim}}

\newcommand{\Br}{\textrm{Br}}

\begin{document}

\begin{abstract}
We show that the triply graded Khovanov-Rozansky homology of the torus link $T_{n,k}$ stablizes as $k\to \infty$.  We explicitly compute the stable homology, as a ring, which proves a conjecture of Gorsky-Oblomkov-Rasmussen-Shende.  To accomplish this, we construct complexes $P_n$ of Soergel bimodules which categorify the Young symmetrizers corresponding to one-row partitions and show that $P_n$ is a stable limit of Rouquier complexes.  A certain derived endomorphism ring of $P_n$ computes the aforementioned stable homology of torus links. 
\end{abstract}

\title[Categorified Young symmetrizers and torus links]{Categorified Young symmetrizers and stable homology of torus links}

\author{Matthew Hogancamp} \address{University of Southern California}

\maketitle

\setcounter{tocdepth}{1}
\tableofcontents

\section{Introduction}
\label{sec-intro}
Over the past several years, there have appeared a number of fascinating conjectures \cite{ORS12,ObSh12,GOR12,GORS12,GorNeg-pp,GNR16} relating certain link invariants with certain Hilbert schemes and rational Cherednik algebras.  Conceptually, these conjectures can be thought of as concrete manifestations of a deep connection between the known, mathematically rigorous constructions of link homology and physical approaches to link homology \cite{Cherednik13,NawObl-pp}.  One such conjecture \cite{GORS12} states that Khovanov-Rozansky homology of the $(n,k)$ torus link $T_{n,k}$ can be computed from the Hilbert scheme of points on the complex surface $f(z,w)=0$, where $f(z,w):=z^n+w^k$.  From the definition of Khovanov-Rozansky homology, it is not at all clear that there should be such a connection.\footnote{It should be noted that D.~Maulik \cite{Maul15} has proven the decategorified version of this conjecture, as well as its generalization to algebraic links.  Maulik's proof uses skein theory techniques which are currently unavailable in the categorical context.}  Nonetheless, this and many other conjectures are firmly supported by the experimental evidence.

In this paper, we show that the Khovanov-Rozansky homology $H_{KR}(T_{n,k})$ stabilizes as as $k\to\infty$, and prove a limiting version of the above quoted conjecture:

\begin{theorem}\label{thm-introHHH}
After an overall shift in the trigrading, the integral triply graded Khovanov-Rozansky homology of the torus links $T_{n,k}$ approach the following limit as $k\to\infty$:
\[
H_{\text{KR}}(T_{n,\infty}) \cong \ring[u_1,\ldots,u_n,\xi_1,\ldots,\xi_n].
\]
This is an isomorphism of triply graded rings, where $u_k$ is an even indeterminate of degree $\deg(u_k)=q^{2k}t^{2-2k}$, and $\xi_k$ is an odd indeterminate of tridegree $q^{2k-4}at^{2-2k}$.   Here, $a$ and $t$ denote Hochschild degree and homological degree, respectively.
\end{theorem}
Note that the $\xi_k$, being odd variables, are assumed to anti-commute and square to zero.  The algebra structure on $H_{\text{KR}}(T_{n,\infty})$ comes from its being identified with a certain derived endomorphism ring of the categorified Young idempotent $P_n$, which is discussed below.

To prove this theorem, we adopt the approach of \cite{Kh07}\footnote{For Khovanov-Rozansky homology over the integers, see \cite{Kras10}}, which constructs triply graded link homology from Hochschild homology---equivalently, Hochschild cohomology---of Soergel bimodules.  Below, we let $\SBim_n$ denote the category of Soergel bimodules (over $\Z$) associated to the symmetric group $S_n$.  Associated to each $n$-strand braid $\b$ one has the \emph{Rouquier complex} $F(\b)\in\Ch(\SBim_n)$, which is a chain complex of Soergel bimodules, well defined up to canonical homotopy equivalence.  Hochschild cohomology of bimodules gives a functor $\HHH:\K(\SBim_n)\rightarrow \Z\mod^{\Z\times \Z\times \Z}$ from the homotopy category to the category of triply graded abelian groups.   In \cite{Kh07}, Khovanov proves that $\HHH(F(\b))$ is a well-defined invariant of the braid-closure $\hat{\b}$, up to isomorphism and an overall shift in tridegree.   The shift can be fixed by a normalization\footnote{See our Corollary \ref{cor:invariance}.}.  The resulting link invariant is isomorphic to Khovanov-Rozansky homology \cite{KR08b}.

Let $x=\sigma_{n-1}\cdots \sigma_2\sigma_1$ denote the positive braid lift of the $n$-cycle $(n,n-1,\ldots,1)$, and let $X=F(x)$ denote the corresponding Rouquier complex.  Note that the closure of $x^m$ is the $(n,m)$ torus link.  Let $\one = R$ denote the trivial Soergel bimodule.  Motivated by work of Rozansky \cite{Roz10a}, we show that powers $X^{\otimes k}$ approach a well-defined limit:
\begin{theorem}\label{thm-introHocolim}
There is a family of chain maps $\{f_k:X^{\otimes k}\rightarrow X^{\otimes k+1}\}_{k=0}^\infty$ whose homotopy colimit $P_n\in \Ch^-(\SBim_n)$ satisfies:
\begin{itemize}
\item[(P1)] $P_n$ kills Bott-Samelson bimodules.
\item[(P2)] $\Cone(\eta)$ is constructed from Bott-Samelson bimodules, where $\eta:\one=X^{\otimes 0}\rightarrow P_n$, is the structure map associated to homotopy colimits.
\end{itemize}
Furthermore, the pair $(P_n,\eta)$ is uniquely characterized by (P1) and (P2) up to canonical equivalence: if $(P_n',\eta')$ is another pair satisfying (P1) and (P2) then there is a unique chain map $\phi:P_n\rightarrow P_n'$ up to homotopy such that $\phi\circ\eta \simeq \eta'$, and this map is a chain homotopy equivalence.
\end{theorem}
This theorem is restated and proven in \S \ref{sec-projector}.  Axiom (P1) means that $P_n\otimes B_i \simeq 0 \simeq B_i \otimes P_n$ for each Bott-Samelson bimodule $B_1,\ldots,B_{n-1}$, while axiom (P2) means that $\Cone(\eta)$ is homotopy equivalent to a complex whose chain bimodules are direct sums of tensor products of the $B_i$ with grading shifts.  The axioms ensure that $P_n\otimes \Cone(\eta)\simeq 0 \simeq \Cone(\eta)\otimes P_n$, which implies that $\eta:\one\rightarrow P_n$ becomes a homotopy equivalence after tensoring on the left or right with $P_n$.  We call a morphism $\eta:\one\rightarrow P$ with this property a \emph{unital idempotent}.  We develop some theory of such idempotents in a separate note \cite{Hog17a}.  There is a dual object $P_n^\vee\in \Ch^+(\SBim_n)$ which is supported in non-negative homological degrees.  This complex is equipped with a map $\eta^\vee:P_n^\vee\rightarrow \one$ which makes $P_n^\vee$ into a \emph{counital idempotent}.  The general theory of such idempotents implies that $P_n$ is a unital algebra in the homotopy category of Soergel bimodules, while $P_n^\vee$ is a coalgebra.  We prefer algebras to coalgebras, hence we prefer $P_n$ to $P_n^\vee$.  The two complexes are related by the application of a contravariant duality functor.

The idempotent complex constructed here categorifies the Young idempotent $p_{(n)}\in \HB_n$ labeled by the 1-row partition $(n)$, which is the Hecke algebra lift of the Jones-Wenzl projector.  Thus, this theorem lifts previous categorifications of the Jones-Wenzl projector to the setting of Soergel bimodules.  See \S \ref{subsec:introIdempts} for more.

Given Theorem \ref{thm-introHocolim}, it follows that $\HHH(P_n)$ is a colimit of triply graded homologies of of the $(n,k)$-torus links up to shifts (see also Corollary \ref{cor:stabilization}).  In \S \ref{subsec-HHH} we introduce the triply graded hom space $\DHom(C,D)$ between complexes of Soergel bimodules, and we show that $\HHH(C)\cong \DHom(R,C)$.  This reformulation of $\HHH$ will be quite useful, and is our main reason for preferring Hochschild cohomology over homology.  Indeed, some general arguments \cite{Hog17a} show that if $\eta:\one\rightarrow P$ is a unital idempotent in a monoidal category $\AC$, then $\End_{\AC}(P)\cong \Hom_{\AC}(\one,P)$ is a commutative $\End_{\AC}(\one)$-algebra.  Our main theorem is:

\begin{theorem}\label{thm-introEndP}
There are isomorphisms of triply graded algebras
\[
\HHH(P_n)\cong \DEnd(P_n)\cong \ring[u_1,\ldots,u_n,\xi_1,\ldots,\xi_n]
\]
where the tridegrees are $\deg(u_k)=q^{2k}t^{2-2k}$, and $\deg(\xi_k)=q^{2k-4}t^{2-2k}a$.
\end{theorem}

%  It is remarkable that the complex $P_n$, which is quite a difficult object to describe explicitly, has such simple endomorphism ring.

%Let $\sigma_i^\pm$ denote an elementary generator of the $n$-stranded braid group, and let $\underlin{b}=\sigma_{i_1}^\pm\cdots \sigma_{i_r}^\pm$.  Associated to $\underline{b}$ one has the \emph{Rouquier complex} $F_{\underline{b}}\in \Ch(\SBim_n)$, which is a complex of Soergel bimodules.  This complex $F_{\underline{b}}$ only depends on the braid r$b$ epresented by $\underline{b}$ up to a canonical homotopy equivalence \cite{}.

%The Hochschild homology of bimodules determines a functor from $\HH:\Ch(\SBim_n)\rightarrow \Ch(\Z\vect^{\Z\times \Z)$ to the category of complexes of bigraded $\ring$-vector spaces.  One of these gradings is the Hochschild grading, and the second is an internal grading owing to the fact that objects of $\SBim_n$ are in fact graded bimodules.  Together with the homological degree of complexes, these form the three gradings mentioned in Theorem \ref{thm-introHHH}.

\begin{remark}
In particular the stable triply graded homology of torus links with $\Z$ coefficients has no torsion.  In contrast, the stable $\sl_N$ homology is expected to have $p$-torsion for all sufficiently large $p$,depending on $N$ (see \S 2.5 of \cite{GOR12}). 
\end{remark}

\subsection{The main idea}
\label{subsec-introIdea}
Theorem \ref{thm-endP} says that $P_n$ is acted on by a polynomial ring.  We will construct $P_n$ in such a way that this structure is evident: the action of the variable $u_n$ represents a certain periodicity in an expression of $P_n$ in terms of $P_{n-1}$. 

To motivate the basic idea, we comment first on the decategorified Young symmetrizers.  Fix an integer $n\geq 1$, and let $j_k\in \Br_n$ denote the braid obtained by winding a single strand around $k-1$ parallel strands.  That is to say,
\[
j_k = \sigma_{k-1}\sigma_{k-2}\cdots \sigma_1\sigma_1\cdots \sigma_{k-2}\sigma_{k-1}.
\]
We will refer to the $j_2,\ldots,j_{n}\in \Br_n$ as Jucys-Murphy braids; they generate an abelian subgroup of $\Br_n$.  Note that the full twist can be written as $j_2j_3\ldots j_{n}$.  The Young symmetrizers $p_T\in \HB_n$ (see \S \ref{subsec:introIdempts}) can be defined as projections onto simultaneous eigenspaces of the action of $j_2,\ldots,j_n$ on on $j_k\in\HB_n$.  A special case of this construction yields the idempotents $p_n$ of interest to us here:%We will regard braids also as linear operators on the Hecke algebra $\HB_n$.  The % The relationship between these elements and the idempotent $p_n\in\HB_n$ is described in terms of the group homomorphism $\rho:\Br_n\rightarrow \HB_n^\times$ which sends the standard braid generator $\sigma_i$ to the standard Hecke algebra generator $T_i$.  Then
\begin{proposition}\label{prop-decatConstruction}
There is a unique family of elements $p_n\in \HB_n$ ($n \geq 1$) satisfying $p_1=1$ together with
\begin{subequations}
\begin{equation}\label{eq-minpoly}
(j_{n}-q^{2n})(j_n-1)p_{n-1}=0
\end{equation}
\begin{equation}\label{eq-decatPnExpr}
p_{n}=\frac{j_n-q^{2n}}{1-q^{2n}}p_{n-1}
\end{equation}
\end{subequations}
for $n\geq 2$.  Here, we are regarding $p_{n-1}$ as an element of $\HB_n$ via the standard inclusion $\HB_{n-1}\subset \HB_n$.  In other words, these equations state that left multiplication by $j_{n}$ acts diagonalizably on $p_{n-1}\HB_{n}p_{n-1}$ with eigenvalues $q^{2n}$ and $1$, and $p_{n}$ is the projection onto the $1$ eigenspace.
\end{proposition}
A proof of the above can be obtained by a ``decategorification'' of arguments in this paper.   In the author's joint work \cite{ElHog17a-pp} with B.~Elias, we propose a categorification of all of the Young symmetrizers by analyzing the categorical spectral theory of the Jucys-Murphy braids acting on the Soergel category.  The construction below is a special case of that story.

Next we state our categorical analogue of Proposition \ref{prop-decatConstruction}. Regard $P_{n-1}$ as an object of $\Ch^-(\SBim_{n})$ by extending scalars, and let $J_n$ denote the Rouquier complex associated to the Jucys-Murphy braid $j_n$.  We will denote tensor product $\otimes_R$ by $\otimes$, or sometimes simply by juxtaposition.  Introduce the diagram notation:
\[
P_{n-1} \ = \ \sspic{.4in}{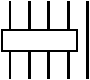} \ \ \ \ \ \ \ \ \ \ \ \ \ \ \ \ \ \ \ \ \ P_{n-1} J_n \ = \ \sspic{.4in}{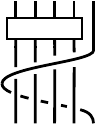},
\]
In \S \ref{sec:structure} we show that there are chain maps $\b_{(n-1,1)}', \b_{(n)}' \in \Hom(P_{n-1}, P_{n-1}J_n)$ of degree\footnote{Note that we write our degrees multiplicatively; '$q$' denotes bimodule, or quantum, degree, and '$t$' denotes homological degree.} $q^{2n}t^{2-2n}$, respectively $q^0t^0$, such that the mapping cones satisfy:
\begin{itemize}
\item $\Cone(\b_{(n-1,1)}')$ kills Bott-Samelson bimodules, up to homotopy equivalence.
\item $\Cone(\b_{(n)}')$ is constructed from Bott-Samelson bimodules, up to homotopy equivalence. 
\end{itemize}
These properties imply that
\begin{equation}
\Cone(\b_{(n-1,1)}')\otimes \Cone(\b_{(n)}')\simeq 0 \simeq \Cone(\b_{(n)}')\otimes \Cone(\b_{(n-1,1)}')
\end{equation} 
This is the categorical analogue of Equation (\ref{eq-minpoly}).  Now, how do we recover the categorified projection from its categorified minimal polynomial?  Expanding the right-hand side of expression (\ref{eq-decatPnExpr}) into a power series in positive powers of $q$ suggests we should consider the following semi-infinite diagram of chain complexes and chain maps:

\vskip 7pt
\begin{equation}\label{eq-introPn}
\begin{minipage}{3.8in}
\begin{tikzpicture}[baseline=-0.25em]
\matrix (m) [ssmtf,column sep=5em, row sep=5em]
{\sspic{.4in}{Pone}(2n)\ip{2-2n} & \sspic{.4in}{PJ}(0)\ip{0} \\
\sspic{.4in}{Pone}(4n)\ip{4-4n} &  \sspic{.4in}{PJ}(2n)\ip{2-2n} \\
\sspic{.4in}{Pone}(6n)\ip{6-6n} &  \sspic{.4in}{PJ}(4n)\ip{4-4n} \\
\cdots & \cdots\\};
\node (11) at (m-1-1)[xshift=0em,yshift=-1em] {};
\node (21) at (m-2-1)[xshift=0em,yshift=-1em] {};
\node (31) at (m-3-1)[xshift=0em,yshift=-1em] {};
\node (22) at (m-2-2)[xshift=-5em,yshift=1em] {};
\node (32) at (m-3-2)[xshift=-5em,yshift=1em] {};
\node (42) at (m-4-2)[xshift=-5em,yshift=1em] {};
\path[->,font=\scriptsize]
(m-1-1) edge node[auto] {} (m-1-2)
(11) edge node[auto] {} (22)
(m-2-1) edge node[auto] {} (m-2-2)
(21) edge node[auto] {} (32)
(m-3-1) edge node[auto] {} (m-3-2)
(31) edge node[auto] {} (42);
\end{tikzpicture}
\end{minipage}
\end{equation}
\vskip 7pt\noindent
where $(k)\ip{\ell}$ is our notation for a shift in $q$-degree and homological degree by $+k$, $+\ell$, respectively.  The horizontal arrows above are given by $\b_{(n-1,1)}'$ and the diagonal arrows are $\b_{(n)}'$, with the appropriate shift functors applied.  This diagram defines a chain map
\[
\Psi\ :\ \ring[u_n]\otimes P_{n-1}(2n)\ip{2-2n} \ \ \longrightarrow \ \ \ring[u_n]\otimes P_{n-1}J_n,
\]
where $u_{n}$ is a formal indeterminate of bidegree $q^{2n}t^{2-2n}$.  In \S \ref{sec:structure} we prove that $P_{n}\simeq \Cone(\Psi)$.  This is the categorical analogue of expression (\ref{eq-decatPnExpr}).  %Now, all of the essential properties of $P_{n}$ become simple observations about $P_{n-1}$, $\b_1$, and $\b_2$.  This is taken up in \S \ref{subsec-constructionOfPn}.

Since $P_n$ kills Bott-Samelson bimodules, it follows that $P_n$ absorbs Rouquier complexes up to homotopy.  Thus we may tensor our description above on the right with the Rouquier complex associated to $\sigma_{n-1}\inv\cdots \sigma_2\inv\sigma_1\inv$, obtaining an equivalent expression for $P_n$ in terms of the complexes
\[
P_{n-1}Y_n  \ = \  \sspic{.4in}{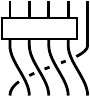} \ \ \ \ \ \ \ \ \ \ \ \ \ \ \ \ \ \ \ \ \ P_{n-1} X_n \ = \ \sspic{.4in}{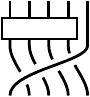}
\]
The calculation of $\HHH(P_n)$ is then proven by induction, with the inductive step provided by the Markov move.

\subsection{Colored HOMFLYPT homology}
\label{subsec:coloredHomology}
The colored HOMFLYPT polynomial is an invariant of oriented links whose components are colored by partitions or, equivalently, irreducible $\sl_N$ representations for some $N$.  Our projector $P_n\in \K^-(\SBim_n)$ slides over and under crossings (Proposition \ref{prop-Psliding}), hence can be used to define a $\lambda$-colored triply graded link homology, where $\lambda$ is a one-row partition (corresponding to the representations $\Sym^n(\C^N)$).  This is the first such construction, though it should be mentioned that a $\Lambda^n$-colored triply graded link homology was constructed by Webster-Williamson in \cite{WW09}, and Cautis \cite{Cau12} has constructed colored $\sl_N$-link homology for all colors. 

 Let $L\subset S^3$ be a link whose components are labelled by non-negative integers $\ell_1,\ldots,\ell_r$, called the colors.  Choose a presentation of $L$ as the closure of an $n$-stranded braid, and let $k_1,\ldots,k_n$ be the colors as one reads from left to right along, say, the bottom of the braid diagram.  Then replace the $i$-th strand by $k_i$ parallel copies, and insert an idempotent $P_{k_i}$.  The result is a chain complex in $\K^-(\SBim_m)$, where $m=k_1+\cdots +k_n$.  Evaluating $\HHH(-)$ on this complex defines the colored homology $H(L; \ell_1,\ldots,\ell_r)$, up to an overall shift in tridegree, depending on the braid index, the colors, and the writhe of the braid.

The value of the $n$-colored unknot is $\HHH(P_n)\cong \ring[u_1,\ldots,u_n,\xi_1,\ldots,\xi_n]$ with degrees as in Theorem \ref{thm-introHHH}.  The homology $H(L;\ell_1,\ldots,\ell_r)$ is a module over the homology of an unlink $H(U; \ell_1)\otimes_\ring \cdots \otimes_\ring H(U;\ell_r)$, by general arguments.  This action does not depend on any choices, up to nonzero scalars.  Thus, $\ring[u_1,\ldots,u_n,\xi_1,\ldots,\xi_n]$ acts as cohomology operations on our $\Sym^n$-colored link homology.  We omit the details in the interest of length, and also since they are similar to results in our earlier work \cite{H14a}.

\subsection{Young idempotents, Jones-Wenzl projectors, and \texorpdfstring{$\sl_N$}{sl(N)} homology}
\label{subsec:introIdempts}

Within the Hecke algebra $\HB_n$ one has a canonically defined, complete set of central idempotents $p_\lambda$, indexed by partitions of $n$.  A standard tableau $T$ on $n$ boxes can be thought of as a sequence of Young diagrams (equivalently, partitions) $\emptyset, \lambda_1,\ldots,\lambda_n$, in which each $\lambda_i$ differs from $\lambda_{i-1}$ in the addition of a single box.  By multiplying the corresponding idempotents $p_{\lambda_i}\in\HB_i\subset \HB_n$ together, one obtains an idempotent $p_T\in \HB_n$.   The $p_T$ are primitive, and $p_\lambda=\sum_T p_T$, where the sum is over all standard tableaux with shape $\lambda$.  The idempotents $p_T$ are $q$-analogues of the classical Young symmetrizers \cite{Gyo86}.  In this paper we are concerned with the categorification of $p_{(n)}$, where $(n)$ is the one-row partition of $n$.

For each integer $N\geq 1$, one may consider the quotient of the Hecke algebra $\HB_n$ by the ideal generated by the $p_\lambda$, where $\lambda$ has more than $N$ parts.  We call this the $\sl_N$ quotient of $\HB_n$, since it is isomorphic to $\End_{U_q(\sl_N)}(V^{\otimes n})$, where $V$ is a $q$-version of the $N$-dimensional standard representation of $\sl_N$.  In case $N=2$, the endomorphism ring is the Temperley-Lieb algebra $\TL_n$.  The image of $p_{(n)}$ in $\End_{U_q(\sl_N)}(V^{\otimes n})$ is called the Jones-Wenzl projector.  The Jones-Wenzl projector can also be defined as the projection operator of the $U_q(\sl_N)$-representation $V^{\otimes n}$ onto the $q$-symmetric power $\Sym^n(V)$.  The Jones-Wenzl projectors play an important role quantum topology, where they can be used to define the $\Sym^n(V)$-colored $\sl_N$ link polynomial in terms of the uncolored polynomial.

The idempotent complex $P_n$ of Soergel bimodules which we construct here is a lift of the categorified Jones-Wenzl projectors constructed by Cooper-Krushkal \cite{CK12a} and Rozansky \cite{Roz10a} for $\sl_2$, D.~Rose \cite{Rose14} for $\sl_3$, and later S.~Cautis \cite{Cau12} for all $\sl_N$.  More precisely, for any integer $N\geq 1$, there is a monoidal functor from Soergel bimodules to the category of $\sl_N$ matrix factorizations \cite{Ras06,KR08}, or $\sl_N$ foams \cite{QuRo14}.  The image of our $P_n$ under these functors is homotopy equivalent to the corresponding categorified Jones-Wenzl projectors.  The simplifications in this paper go through in the $\sl_N$ specializations, and we obtain:

\begin{theorem}\label{thm-introSLNhomology}
There is some $\ring[u_1,\ldots,u_n]$-equivariant differential $d_N$ on $\ring[u_1,\ldots,u_n,\xi_1,\ldots,\xi_n]$ which computes the limiting $\sl_N$ homology of the $(n,k)$ torus links as $k\to \infty$.  The degrees of the generators are obtained from those in Theorem \ref{thm-introHHH} by specializing $a\mapsto t\inv q^{2N+2}$.% $(\deg_q,\deg_h)$ are given by $\deg_{\sl_N}(u_k) = (2k,2-2k)$ and $\deg_{\sl_N}(\xi_k)=(2k+2N-2,1-2k)$.
\end{theorem}
%See \S \ref{} for more of a discussion on the relationship with $\sl_N$ link homology.
 The following is conjectured in \cite{GOR12}:
\begin{conjecture}
The differential $d_N$ is determined by $d_N(\xi_k)=\sum u_{i_1}\cdots u_{i_k}$, where the sum is over sequences of integers $1\leq i_{1},\ldots,i_k\leq n$ which sum to $k+N-1$.  This differential is extended to all of $\ring[u_1,\ldots,u_n,\xi_1,\ldots,\xi_n]$ by the graded Leibniz rule with respect to the standard multiplication.
\end{conjecture}
%There are two main difficulties in proving this conjecture.  Firstly, we do not know if the isomorphism $\HHH(P_n)\cong \Z[u_1,\ldots,u_n,\xi_1,\ldots,\xi_n]$ is an isomorphism of algebras.  This means that we do not know if $d_N$ respects the Leibniz rule with respect to the obvious multiplication of polynomials.  Secondly the precise formula for $d_N(\xi_k)$ requires more explicit knowledge of the projector $P_n$ than we currently have.
We do not pursue the connections to $\sl_N$ homology in this paper, so we omit the proof of Theorem \ref{thm-introSLNhomology}.

\subsection{Relation to other work}
This paper extends our previous work \cite{H14a} from the $\sl_2$ setting to the HOMFLYPT setting.  Decategorified, this means working with the Hecke algebra $\HB_n$ instead of its quotient $\TL_n$.  Categorified, this means working with Soergel bimodules rather than the Bar-Natan/Khovanov category of tangles and cobordisms.  In the setting of Soergel bimodules, we are able to give a clean formula for the homology of $P_n$.

\textbf{The one column projector.}  In a sequel to this paper (joint with M.~Abel \cite{AbHog15}), we construct an idempotent complex $P_{1^n}\in \K^-(\SBim_n)$ corresponding to the one-column partition $1+\cdots +1 =n$, and we study the corresponding colored homology.  Both projectors are homotopy colimits of directed systems $R\rightarrow \FT_n\rightarrow \FT_n^{\otimes 2}\rightarrow \cdots$ (degree shifts omitted), but using different maps.  One may think of $P_n$ as the ``head'' and $P_{1^n}$ as the ``tail;'' we won't attempt to make these terms precise, but instead direct the reader to \cite{AbHog15} for an explanation of the relationship between the two projectors.  There are many important differences between the two stories.  Firstly, the projector $P_{1^n}$ admits an explicit combinatorial description, whereas the construction of $P_n$ is quite mysterious.  Secondly, there are polynomial subalgebras $\ring[U_2,\ldots,U_n]\subset \End(P_{1^n})$ and $\ring[u_2,\ldots,u_n]\subset \End(P_n)$ which exhibit $n-1$-fold periodicity in our descriptions of these complexes.  The variables $U_k$ have homological degree 2, and the resulting 2-periodicity of $P_{1^n}$ is a manifestation of a standard phenomenon in commutative algebra \cite{Eis80}.  On the other hand, the variables $u_k$ have homological degree $2-2k$, and the resulting periodicity is a bit more surprising.  Finally, $P_{1^n}$ is not invariant under the operation of tensoring with a Rouquier complex $F(\b)$, though $F(\b)\otimes P_{1^n}$ depends only on the permutation represented by $\b$ up to homotopy. Thus, we obtain a family of \emph{twisted projectors} $w(P_{1^n})$, indexed by permutations $w\in S_n$. %If $w\in S_n$ is an $n$-cycle, the homology $\HHH$ of the twisted projector $w\cdot P_{1^n}$ is an algebra $\ring[v_1,\ldots,v_n,\phi_1,\ldots,\phi_n]$ where the even variables $v_k$ have tridegree $q^{2k-2}t^{-2}$ and the odd variables $\phi_k$ have tridegree $q^2a\inv$.   This surprising similarity between the $\HHH(w\cdot P_{1^n})$ and $\HHH(P_n)$ is an instance of a more general symmetry which is conjectured by Gorsky-Rasmussen in \cite{GoRa15}.

\textbf{Categorical diagonalization.}
The projectors $P_n$ constructed here play an essential role in \cite{ElHog16a}, where they give a method of computing the triply graded Khovanov-Rozansky homology of several infinite families of links, including the $(n,n)$ torus links.  This work, in turn, plays an important role in an eigenedecomposition of the Soergel category $\SBim_n$, developed in joint work with Ben Elias \cite{ElHog17a-pp}, which gives rise to complexes $P_T\in \KC^-(\SBim_n)$, indexed by standard Young tableaux on $n$ boxes, which categorify the Young symmetrizers.

\textbf{The flag Hilbert scheme.}
Even though the $P_T$ have not yet appeared in the literature, E.~Gorsky, A.~Negu\cb{t}, and J.~Rasmussen \cite{GNR16} have a beautiful series of conjectures regarding the triply graded homologies $\DEnd(P_T)$, coming from flag Hilbert schemes.  Together with explicit computations in \cite{GNR16}, our work here proves a conjecture in \emph{loc.~cit.} in a special case.   The computations in \cite{AbHog15} prove the conjecture in the case of the one-column projector.  One aspect of the conjectures in \cite{GNR16} is that, up to regrading, there should be an isomorphism between the $\DHom(P_T,X\otimes P_T)$ and $\DHom(P_{T^\ast},X\otimes P_{T^\ast})$ up to regrading, where $T$ is a standard tableau on $n$-boxes, $T^\ast$ is the transposed tableau, and $X=\sigma_{n-1}\cdots \sigma_2\sigma_1$ is the positive braid lift of a Coxeter element of $S_n$.  The results of this paper and \cite{AbHog15} demonstrate this duality in the case of the one-column and one-row partitions.

For the interested reader, we state the precise regrading.  Introduce new variables $\bar{q}=q^2$, $\bar{t}=t^2q^{-2}$, and $\bar{a} = aq^{-2}$.  Note for example that the variables in Theorem \ref{thm-introHHH} have degrees $\deg(u_k)=\bar{q}\bar{t}^{1-k}$ and $\deg(\xi_k)=\bar{a}\bar{t}^{1-k}$.  The regrading involved in the transposition symmetry is $\bar{q}\leftrightarrow \bar{t}$ and $\bar{a}\mapsto \bar{a}$.  Indeed, the computations in \cite{AbHog15} involve even variables of degree $\bar{t}\bar{q}^{1-k}$ and odd variables of degree $\bar{a}\bar{q}^{1-k}$, as predicted.

\subsection{Outline of the paper}
Section \S \ref{subsec-soergelCat} begins our story with some relevant background on the Hecke algebra and Soergel category (associated to $S_n$).  In \S \ref{subsec-triangCats} we recall some basics regarding complexes and mapping cones. In \S \ref{subsec-Pn} we give axioms which characterize our categorified Young idempotent $P_n$ and deduce some basic consequences, such as uniqueness and centrality of $P_n$.  In \S \ref{subsec-rouquier} we recall Rouquier complexes.  Section \S \ref{subsec-projAsLimit} constructs $P_n$ as a homotopy colimit of Rouquier complexes, proving Theorem \ref{thm-introHocolim} of the introduction.

Section \S \ref{sec-KRhomology} establishes the relevant theory concerning Hochschild cohomology and triply graded homology.  First, \S \ref{subsec-HHH} introduces derived categories and gives a categorical reformulation of the functor $\HHH$.  In \S \ref{subsec-adjunction} we study the ``partial Hochschild cohomology'' functor, or ``partial trace'' $\Tr$, and proves an important adjunction isomorphism.  Then \S \ref{subsec:markov} establishes the Markov relations for $\Tr$.   In \S \ref{subsec:stabilization} we formulate a precise statement regarding the stabliziation of triply graded homology of torus links.

Section \S \ref{sec:structure} studies the structure of $P_n$, and is the technical heart of this paper.  In \S \ref{subsec-statement} we state our main structural theorems (Theorem \ref{thm-endP} and Theorem \ref{thm-eigenconesAreOrtho}).  Section \S \ref{subsec-constructionOfPn} contains our inductive construction of $P_n$ as outlined in \S \ref{subsec-introIdea} of the introduction.  In \S \ref{subsec:endGeneral} we recall some general results on the endomorphism rings of unital idempotents and apply them to $P_n$.  It is here that we state the isomorphism between $\HHH(P_n)$ and $\DEnd(P_n)$. In \S \ref{subsec-endP} we deduce Theorem \ref{thm-endP} from Theorem \ref{thm-eigenconesAreOrtho}, using the aforementioned construction of $P_n$.     In \S \ref{subsec:hookEigenmap} we complete the proofs of Theorem \ref{thm-endP} and Theorem \ref{thm-eigenconesAreOrtho}.  Finally, in \S \ref{subsec:QnProps} we establish properties of the complex $Q_n\simeq \Cone(P_n\buildrel u_n\over \rightarrow P_n)$, which will be useful in future work.

\subsection{Notation}
\label{subsec-notation}
Let $\mathcal{A}$ denote an additive category.  Throughout this paper we use $\Ch(\mathcal{A})$ to denote the category of chain complexes over $\mathcal{A}$, and $\K(\AC)$ its homotopy category.  Differentials will be homogeneous of homological degree $+1$.  We use superscripts $+,-,b$ to denote the full subcategories consisting of complexes which are bounded from below, respectively  above, respectively above and below, in homological degree.  We denote isomorphism in $\Ch(\mathcal{A})$ by $\cong$ and isomorphism in $\mathcal{K}(\mathcal{A})$ (that is to say, homotopy equivalence) by $\simeq$.

If $A$ is a chain complex, we let $A[1]$ denote the usual suspension: $A[1]_i = A_{i+1}$.  We let $A\ip{1}=A[-1]$ denote the complex obtained by shifting $A$ up in homological degree.  By convention $[1]$ and $\ip{1}$ also negate the differential.  Most of the complexes in this paper will be equipped with an additional grading, which we call $q$-degree.  We denote the \emph{upward} grading shift in $q$-degree by $(1)$, so that $(A(a)\ip{b})_{i,j}=A_{i-a,j-b}$.  A chain map $f:A(i)\ip{j}\rightarrow B$ may be regarded as a map $A\rightarrow B$ of bidegree $(i,j)$.

We will also consider triply graded complexes.  This third grading typically appears in the following way: let $S$ be a graded ring, and let $\CS=D(S\gmod)$ be the derived category of graded left $S$-modules.  The grading on $S$-modules will be called the quantum grading, or $q$-grading, and the homological degree in $\CS$ will be called the derived grading, Hochschild grading, or $a$-grading.  Since $\CS$ is, in particular, an additive category, we may form the homotopy category of complexes over $\CS$, obtaining $\DS:=\KC(\CS)$.  The homological grading in $\DS$ will be called the homological grading, or $t$-grading.  For each triple of integers $i,j,k\in \Z$ and each $X\in \DS$, we have the grading shift $X(i,j)\ip{k}$ which shifts $X$ up by $i,j,k$ in the $q$-grading, $a$-grading, and $t$-grading, respectively.  Note that $S\mod$ embeds in $\CS$, which induces an embedding of $\KC(S\mod)$ into $\DS$.  The shifts in $q$-degree and homological degree in $\KC(S-\mod)$ are related to the shifts in $\DS$ by:
\[
X(i,0)\ip{0} = X(i) \ \ \ \ \ \ X(0,0)\ip{k} = X\ip{k}.
\]
Given $X,Y\in \DS$, we have the enriched hom space
\[
\DHom_\DS(X,Y):=\bigoplus_{i,j,k}\Hom_{\DS}(X(i,j)\ip{k},Y).
\]
The Poincar\'e series of this triply graded abelian group is
\[
\PS_{X,Y}(q,a,t):=\sum_{i,j,k\in\Z} q^ia^jt^k \rank_\ring \Hom_\DS(X(i,j)\ip{k},Y).
\]
The degree of a homogeneous element $f\in \DHom_\DS(X,Y)$ will often be written multiplicatively as $\deg(f) = q^i a^j t^k$.

We also utilize various ``internal homs'' throughout the paper.  As an example, let $S$ be a graded ring.  Suppose $A,B$ are complexes of graded $S$-modules, then we let $\Homg_S(A,B)$ denote the bigraded complex such that $\Homg^{i,j}_S(A,B)$ is the abelian group of bihomogeneous $S$-linear maps $A\rightarrow B$ of $q$-degree $i$ and homological degree $j$.  That is to say $\Homg^{i,j}_S(A,B):=\prod_{k\in \Z}\Hom_{S}(A(i)_k,B_{k+j})$.  There is a differential on $\Homg_S(A,B)$ defined by $f\mapsto d_B\circ f - (-1)^j f\circ d_A$ for each $f\in \Homg^{i,j}_S(A,B)$.

\begin{remark}
The homology of $\Homg_{\SBim_n}(A,B)$ is the $a$-degree zero part of $\DHom_{\DC_n}(A,B)$.
\end{remark}

The usual adjunction isomorphism takes the following form: let $S$ and $R$ be non-negatively graded rings.  Let $M,K,N$ be complexes of graded left $S$-modules (resp.~$(R,S)$-bimodules, resp.~ left $R$-modules).  Then
\[
\Hom_{\Ch(R\mod)}(K\otimes_S M,N) \cong \Hom_{\Ch(S\mod)}(M,\Homg_S(K,N)).
\]

Finally we often abbreviate tensor products $M\otimes N$ simply by writing $MN$.  We also let $M$ denote the identity endomorphism of $M$, so expressions such as $Mf$ and $fM$ will denote $\Id_M\otimes f$ and $f\otimes \Id_M$.

\subsection{Acknowledgements}
The author would like to thank M.~Abel and B.~Elias for their comments on an earlier version of this work, A.~Oblomkov for his encouraging remarks, E.~Gorsky, A.~Negu\cb{t}, and J.~Rasmussen for sharing their wonderful work on flag Hilbert schemes, and all of the aforementioned for the enlightening conversations.   The author would also like to send special thanks to his advisor, Slava Krushkal, for his continued guidance and support, and without whom he might never have discovered the joy of categorification.

\section{A categorified Young symmetrizer}
\label{sec-projector}
In this section we introduce the category of Soergel bimodules, and give a set of axioms which characterizes a unique complex $P_n$ of Soergel bimodules, up to canonical equivalence.  This complex categorifies a certain idempotent in the Hecke algebra.  

\subsection{Soergel's categorification of the Hecke algebra}
\label{subsec-soergelCat}
Fix throughout this section an integer $n\geq 1$.  We use letters such as $s,t,u$ to denote simple transpositions in the symmetric group $S_n$.  Two simple transpositions $s=(i,i+1)$, $t=(j,j+1)$ are \emph{adjacent} if $i=j\pm 1$ and \emph{distant} if $|i-j|\geq 2$.  The Hecke algebra $\Hecke_n$ for $S_n$ is the $\Q(q)$-algebra generated by elements $T_1,\ldots,T_{n-1}$ with defining relations
\begin{enumerate}
\item $(T_i+q^2)(T_i-1) = 0$ for $1\leq i\leq n-1$.
\item $T_iT_{i+1}T_i = T_{i+1}T_iT_{i+1}$ for $1\leq i\leq n-2$.
\item $T_i T_j = T_j T_i$ when $|i-j|\geq 2$.
\end{enumerate}
The Hecke algebra is a $q$-deformation of the group algebra $\Q[S_n]$ in the sense that setting $q=1$ recovers the defining relations in $\Q[S_n]$.  It will be convenient to work with a different set of generators $\{b_i:=q\inv (1 -T_i)\}$, with respect to which the defining relations become
\begin{enumerate}
\item $b_i^2 = (q+q\inv)b_i$ for $1\leq i\leq n-1$.
\item $b_i b_{i+1} b_i - b_i = b_{i+1} b_{i} b_{i+1} - b_{i+1}$ for $1\leq i\leq n-2$.
\item $b_i b_j = b_j b_i$ when $|i-j|\geq 2$.
\end{enumerate}
Our $q$ is denoted by $v\inv$, for instance in \cite{EW14}.  We often abuse notation and write $b_s$ when we mean $b_i$, where $s$ is the simple transposition $(i,i+1)$.

Now we describe a \emph{categorification} of the Hecke algebra; the algebra $\HB_n$ gets replaced by a certain monoidal category $\SBim_n$, the generators $b_i$ get replaced by objects $B_i$, and the defining relations lift to isomorphisms in $\SBim_n$.  More precisely, $\HB_n$ is isomorphic to the split Grothendieck group (actually a ring) of $\SBim_n$.  For this and further background on Soergel bimodules, we refer the reader to \cite{EKh10}.

Put $R_n:=\ring[x_1,\ldots,x_n]$, graded by placing each $x_i$ in degree 2.  We call this the $q$-degree, to distinguish it from the homological degree of complexes, and we write $\deg_q(x_i)=2$.  We let $(1)$ denote the functor which shifts graded objects up by one unit in $q$-degree: $(M(1))_i=M_{i-1}$.  When the index $n$ is understood, we simply write $R_n=R$.

Set $W:=S_n$.  Below, letters such as $s,t,u$ will denote simple transpositions $(i,i+1)\in W$.  There is a left action of $W$ on $R$ given by permuting variables.  Let $R^s\subset R$ denote the subring of polynomials which are fixed-points for the action of $s$.  Let $B_s$ denote the graded $R$ bimodule $B_s:=R\otimes_{R^s} R(-1)$.  We also write $B_i$ when we mean $B_s$ ($s=(i,i+1)$).  
\begin{notation}
We often write $A\otimes_R B$ as $A\otimes B$ or simply $AB$.
\end{notation}
Bimodules of the form $B_s B_t \cdots  B_u(k)$ are called \emph{Bott-Samelson} bimodules.  Let $\SBim_n$ denote the smallest full subcategory of finitely generated, graded $(R,R)$-bimodules which contains the Bott-Samelson bimodules, and which is closed under $\otimes_R$, $\oplus$, $(\pm 1)$, and taking direct summands.  An object of $\SBim_n$ is called a \emph{Soergel bimodule}, and $\SBim_n$ is called the \emph{Soergel category}.  We will not need much background regarding Soergel bimodules.  Everything we need to know is summarized below, and in the sequel we will assume familiarity with the remainder of this section.

Let $s=(i,i+1)\in S_n$ be a simple transposition.  Then there is a canonical map $B_s= R\otimes_{R^s}R(-1)\rightarrow R(-1)$ sending $1\otimes 1\mapsto 1$.  There is also a canonical bimodule map $R(1)\rightarrow B_s$ sending $1\mapsto x_i\otimes 1-1\otimes x_{i+1}$  These maps will be refered to as the \emph{dot maps}.

The following facts ensure that the defining relations in $\HB_n$ (in terms of the $b_i$) lift to isomorphisms in $\SBim_n$.
\begin{enumerate}
\item For any simple transposition $s\in S_n$, we have $B_sB_s \cong B_s(1)\oplus B_s(-1)$.  In fact if $\psi:B_s\rightarrow R(-1)$ is the dot map, then $B_s\psi$ and $\psi B_s$ are projections onto a $B_s(-1)$ summand of $B_sB_s$.
\item If $s,t\in S_n$ are adjacent simple transpositions, then there is a Soergel bimodule $B_{tst}=B_{sts}$ such that $B_s B_t B_s\cong B_{sts}\oplus B_s$.  See equations (4), (5), and (7) in \cite{Kras10}.
\item If $s,u\in S_n$ are distant, then $B_sB_u\cong B_uB_s$.
\end{enumerate}

One more fact we will need is that the bimodules $B_s$ are \emph{self biadjoint}.  That is to say
\[
\Hom(MB_s,N)\cong \Hom(M,NB_s) \ \ \ \ \ \ \ \ \ \text{and} \ \ \ \ \ \ \ \ \ \ \Hom(B_sM,N)\cong \Hom(M,B_sN)
\]
for all $M,N\in \SBim_n$. % This is essentially the reason for our choosing the grading shift in our definition of $B_s$.

We conclude this subsection by collecting some standard definitions and notations:
\begin{definition}\label{def-induction}
Let $R_k=\ring[x_1,\ldots,x_k]$. Note that $R_i\otimes_\ring R_j\cong R_{i+j}$.  Thus we have bilinear functors $(-)\sqcup (-):\SBim_i\times \SBim_j\rightarrow \SBim_{i+j}$ given by tensoring over $\ring$.  The functor $(-)\sqcup \one_{n-k}:\SBim_k\rightarrow \SBim_n$ will be denoted by the letter $I$.  Sometimes we let $\one_n = R_n$ denote the trivial bimodule.  When $n$ is obvious from the context, we will omit it from the notation.
\end{definition}

\subsection{Complexes and mapping cones}
\label{subsec-triangCats}
In this section we recall some basics of mapping cones that will be used throughout.  Let $\CS$ be an additive category.  Associated to any chain map $f:A\rightarrow B$ in $\Ch(\CS)$, one has the \emph{mapping cone} $\Cone(f)\in\Ch(\CS)$.  The chain groups are $\Cone(f)_i = A_{i+1}\oplus B_i$, and the differential is by definition
\[
d_{\Cone(f)} = \matrix{-d_A & \\ f & d_B}
\]
If $f:A\rightarrow B$ is a chain map, then there is a canonical map $B\rightarrow \Cone(f)$ (the inclusion of a subcomplex).  Recall that $A[1]$ is the chain complex $A[1]_i=A_{i+1}$ with differential $-d_A$.  The sign ensures that the projection $\Cone(f)\rightarrow A[1]$ is a chain map.  The sequence of chain maps
\begin{equation}\label{eq-coneTriang}
A\buildrel f\over \longrightarrow B \rightarrow \Cone(f)\rightarrow A[1]
\end{equation}
called a \emph{distinguished triangle} in $\K(\CS)$.  More generally,
\[
A'\buildrel f'\over \longrightarrow B' \buildrel g' \over \longrightarrow C' \buildrel \d'\over \longrightarrow A'[1]
\]
is called a distinguished triangle if there are homotopy equivalences $A\simeq A'$, $B\simeq B'$, $\Cone(f)\simeq C'$, such that the following diagram commutes up to homotopy:
\[
\begin{diagram}
A & \rTo^{f} & B & \rTo^{\iota} & \Cone(f) & \rTo^{\pi} & A[1]\\
\dTo^{\simeq} && \dTo^{\simeq} && \dTo^{\simeq}&&\dTo^{\simeq}\\
A' & \rTo^{f'} & B' & \rTo^{g'} & C' & \rTo^{\d'} & A[1]
\end{diagram}
\]
%If $f:A\rightarrow B$ is a chain map, then the canonical maps $\eta:B\rightarrow \Cone(f)$ and $\e:\Cone(f)\rightarrow A[1]$ satisfy $\Cone(\eta)\simeq A[1]$ and $\Cone(\e)\simeq B[1]$.  In fact, if
%\begin{equation}\label{eq-generalTriang}
%A\buildrel f\over \rightarrow B\buildrel g\over \rightarrow C\buildrel\partial\over \rightarrow A[1]
%\end{equation}
%is a distinguished triangle, then so are
%\begin{equation}\label{eq-rotatedTriang1}
%B\buildrel g\over\rightarrow  C\buildrel \partial\over \rightarrow A[1]\buildrel -f[1]\over \longrightarrow B
%\end{equation}
%and
%\begin{equation}\label{eq-rotatedTriang2}
%C\buildrel \partial \over \rightarrow A[1]\buildrel -f[1]\over \longrightarrow B[1] \buildrel g[1]\over \longrightarrow C[1]
%\end{equation}
The power of mapping cones is that many statements about morphisms can be translated into statements about objects, for example:
\begin{proposition}\label{prop-conesAndEquivs}
A chain map $f:A\rightarrow B$ is a homotopy equivalence if and only if $\Cone(f)\simeq 0$.  Thus, in the distinguished triangle $A\buildrel f\over\rightarrow B \rightarrow C\rightarrow A[1]$, $f$ is a homotopy equivalence if and only if $C\simeq 0$.  
\end{proposition}
\begin{proof}
Standard.  Can be found in Spanier \cite{Spanier}.
\end{proof}

\subsection{Characterization of \texorpdfstring{$P_n$}{Pn}}
\label{subsec-Pn}
Recall the generators $b_i=q\inv(1-T_i)\in \HB_n$ of the Hecke algebra, described in \S \ref{subsec-soergelCat}.  The Young symmetrizer $p_{(n)}$ is uniquely characterized by:
\begin{deth}\label{defthm-decatSymmetrizer}
There is a unique $p_{(n)}\in\Hecke_n$ such that
\begin{itemize}
\item[(p1)] $p_{(n)} b_i =0 = b_ip_{(n)}$ for $1\leq i\leq n-1$.
\item[(p2)] $p_{(n)}-1$ lies in the (non-unital) subalgebra generated by $\{b_1,\ldots,b_{n-1}\}$.
\end{itemize}
\end{deth}
We encourage the reader to compare this with the definition of the Jones-Wenzl projectors, for instance in Chapter 13 of \cite{Lick97}.  The proof that (p1) and (p2) characterize $p_{(n)}$ uniquely is easy, and the proof that $p_{(n)}$ exists follows along the same lines as existence of the Jones-Wenzl projector.  In any case, we will explicitly construct a categorified version of $p_n$.  Note that the standard generators $T_i$ satisfy $T_ip_{(n)}=p_{(n)}$.  Thus, the left $\HB_n$ module $\HB_np_{(n)}$ is a $q$-analogue of the trivial representation.

\begin{remark}
One often considers a different set of generators $H_i = -q\inv T_i$.  In terms of these generators we have $H_ip_{(n)}=-q\inv p_{(n)}$, so that $\HB_np_{(n)}$ can also be considered as a $q$-analogue of the sign representation, depending on one's tastes.
\end{remark}

In the remainder of this subsection we characterize a complex $P_{(n)}\in\K^-(\SBim_n)$ by a similar set of axioms, and we note some consequences of our definitions.  For the remainder of the paper, we write $P_n$ rather than $P_{(n)}$, so that our notation is reminiscent of that for categorified Jones-Wenzl projectors \cite{CK12a,H14a}.

%\begin{definition}\label{def-Kn}
%Let $\K^-(n):=\K^-(\SBim_n)$ denote the homotopy category of Soergel bimodules.
%\end{definition}

\begin{definition}\label{def-N}
Fix an integer $n\geq 1$, and let ${\mathcal{N}_n}\subset \K^-(\SBim_n)$ denote the smallest full subcategory which is closed under homotopy equivalences and contains complexes $N\in\K^-(\SBim_n)$ such that $\one(k)$ does not appear as a direct summand of any chain bimodule, for any $k\in \Z$.  Let ${\mathcal{N}_n}_\perp$ and ${}_\perp{\mathcal{N}_n}$ denote the full subcategories of $\K^-(\SBim_n)$ consisting of complexes $A$ such that $NA\simeq 0$, respectively $AN\simeq 0$, for all $N\in{\mathcal{N}_n}$.  When the index $n$ is understood, we will omit it.
\end{definition}
No nontrivial Bott-Samelson bimodule contains the trivial bimodule $\one$, or any of its shifts, as a direct summand.  Conversely, any nontrivial Soergel bimodule can be resolved by non-trivial Bott-Samelson bimodules.

\begin{example}
Let $s,t\in S_3$ denote the transpositions $s=(1,2)$ and $t=(2,3)$.  Then the longest element of $S_3$ is $sts$.  We have $B_sB_tB_s \cong B_s\oplus B_{sts}$, which implies that $B_{sts}$ is chain homotopy equivalent to the complex $0\rightarrow B_s\rightarrow B_sB_tB_s\rightarrow 0$, where the differential is the inclusion of a direct summand.
\end{example}

Thus we have the following equivalent description of ${\mathcal{N}_n}$:
\begin{proposition}\label{prop-N}
The subcategory ${\mathcal{N}_n}\subset \KC^-(\SBim_n)$ is the smallest full subcategory which is closed under homotopy equivalences and contains those complexes $N$ whose chain bimodules are direct sums of shifts of bimodules of the form $B_{s_1}\cdots B_{s_r}$ with $r\geq 1$.  Thus, the full subcategory ${\mathcal{N}_n}\subset \K^-(\SBim_n)$ is a two-sided tensor ideal.\qed
\end{proposition}

The following is proven by a straightforward limiting argument (see Proposition \ref{prop:simultSimp}).

\begin{proposition}\label{prop:BSkilling}
A complex $D\in\KC^-(\SBim_n)$ lies in ${}_\perp{{\NC_n}}\cap {{\NC_n}}_\perp$ if and only if $B_iD\simeq 0\simeq DB_i$ for all $1\leq i\leq n-1$.
\end{proposition}
\begin{proof}
Clearly an object of ${}_\perp{{\NC_n}}\cap {{\NC_n}}_\perp$ kills $B_i$ on the left and right.  Conversely, suppose $D\in \K^-(\SBim_n)$ kills $B_i$ on the left and right Bott-Samelsons, and let $C\in {{\NC_n}}$ be given.  Without loss of generality, we may assume each chain bimodule $C_i$ is a sum of shifts of Bott-Samelson bimodules $B_{i_1}\cdots B_{i_r}$.  Thus, $D C_i\simeq 0 \simeq  C_iD$ for all $i$.  A straightforward limiting argument (see Proposition \ref{prop:simultSimp}), shows that $DC\simeq 0\simeq CD$.  This completes the proof.
\end{proof}

 A priori it is not obvious whether or not the categories $\NC_\perp$ and ${}_\perp\NC$ are nonzero.  Nonetheless, it will turn out that ${}_\perp{{\NC_n}}={\mathcal{N}_n}_\perp$ is nontrivial and has the structure of a monoidal category, with identity given by $P_n$ defined below.
\begin{deth}\label{defthm-Pn}
There exist a chain complex $P_n \in \K^-(\SBim_n)$ and a chain map $\eta:\one\rightarrow P_n$ such that
\begin{itemize}
\item[(P1)] $P_n\otimes B_i\simeq 0\simeq B_i\otimes P_n$ for all $1\leq i\leq n-1$.
\item[(P2)] $\Cone(\eta)\in \NC_n$.
\end{itemize}
The pair $(P_n,\eta)$ is unique up to canonical equivalence: if $(P_n',\eta')$ is another pair satisfying (P1) and (P2) then there is a unique chain map $\phi:P_n\rightarrow P_n'$ up to homotopy such that $\phi\circ\eta \simeq \eta'$.  This map is a chain homotopy equivalence.  We will call $\eta:\one\rightarrow P_n$ the \emph{unit} of $P_n$.
\end{deth}
The existence of $P_n$ will be proven in \S \ref{subsec-projAsLimit}.  The axioms (P1) and (P2) should be compared with the axioms that define the universal projectors in \cite{CK12a}.  See also \cite{H14a}.

\begin{example}\label{example:P2}
Let $P_2\in \K^-(\SBim_n)$ denote the complex
\[
\begin{diagram}[small]
\cdots &\rTo^{\phi_-} & B_1(5) &\rTo^{\phi_+}& B_1(3)&\rTo^{\phi_-}& B_1(1) &\rTo& \underline{R} &\rTo& 0 
\end{diagram}
\]
where the map $B_1(1)\rightarrow R$ is the dot map, $\phi_+$ is the map sending $1\otimes 1\mapsto x_1\otimes 1 - 1\otimes x_2$, and $\phi_-$ is the map sending $1\otimes 1\mapsto x_2\otimes 1- 1\otimes x_2$.  As usual, we place an underline underneath the term of homological degree zero. The inclusion $R\rightarrow P_2$ defines a chain map $\eta:R\rightarrow P_2$ with $\Cone(\eta)\in\NC_2$.  
\end{example}
\begin{proposition}\label{prop:P2}
The complex $P_2$ and chain map $\eta:R\rightarrow P_2$ satisfy the axioms of Definition/Theorem \ref{defthm-Pn}.
\end{proposition}
\begin{proof}
We prove that $P_2B_1\simeq 0$.  Let $R_s$ denote the bimodule which equals $R$ as a left $R$-module, but whose right $R$-action is twisted by $s$.  That is $f\cdot g\cdot h =fgs(h)$ for $f,h\in R$, $g\in R_s$.  Then the dot maps fit into short exact sequences
\begin{equation}\label{eq-SES2}
\begin{diagram}[small]
0 &\rTo & R_s(1) & \rTo^{1\mapsto x_{i+1}\otimes 1 - 1\otimes x_{i+1}} & B_s & \rTo & R(-1)& \rTo &  0,
\end{diagram}
\end{equation}
\begin{equation}\label{eq-SES1}
\begin{diagram}[small]
0 & \rTo &  R(1) & \rTo &  B_s & \rTo^{1\otimes 1\mapsto 1} &  R_s(-1)  &\rTo & 0.
\end{diagram}
\end{equation}
Note that $P_2$ is the concatenation of infinitely many copies of (\ref{eq-SES1}) and (\ref{eq-SES2}), hence $P_2$ is acyclic.  The bimodule $B_1$ is free as a left (or right) $R$-module, hence the functor $(-)\otimes_R B_1$ is exact.  It follows that $P_2B_1$ has zero homology as well.  On the other hand, $P_2B_1$ is free as a complex of $R\otimes_{R^s} R$-modules.  Thus, $P_2B_1$ is an exact sequence of free modules, hence splits.  The splitting implies that $P_2B_1$ is contractible via a $R\otimes_{R^s} R$-equivariant homotopy.  Since the bimodule action factors through the $R\otimes_{R^s} R$-action, this means that $P_2B_1\simeq 0$, as claimed.  The case of $B_1P_2$ is similar.
\end{proof}
Thus, $P_2$ is a unital idempotent.  The complementary idempotent is the ``tail''
\[
\begin{diagram}
\cdots &\rTo^{\phi_-} & B_1(5) &\rTo^{\phi_+}& B_1(3)&\rTo^{\phi_-}& \underline{B_1}(1) &\rTo& 0
\end{diagram}
\]

\subsection{Some general theory of categorical idempotents}
The existence of $P_n$ will be established in \S \ref{subsec-projAsLimit}. In this section we state some important properties of $P_n$ that follow from general arguments.  All of the relevant facts about categorical idempotents are stated and proven in a separate note \cite{Hog17a}.  Throughout this section, let  $\CS$ be an additive, monoidal category, and set $\AS:=\K^-(\CS)$ (or more generally, $\AS$ can be any triangulated monoidal category).  For us our focus will be on the case $\CS=\SBim_n$.  Recall, isomorphism in $\AS$ is denoted $\simeq$.

\begin{definition}\label{def:imageAndKernel}
A \emph{pair of complementary idempotents} in $\AS$ is a pair of objects $(\counital,\unital)$ such that $\unital\otimes \counital\simeq \counital\otimes \unital\simeq 0$, together with a distinguished triangle
\begin{equation}\label{eq:PQtriang}
\counital\buildrel \e\over \rightarrow \one \buildrel\eta\over \rightarrow \unital \buildrel\d\over \rightarrow \counital[1].
\end{equation}
A \emph{unital idempotent} is an object $\unital$ together with a morphism $\eta:\one\rightarrow \unital$ such that $\unital\otimes \Cone(\eta)\simeq 0\simeq \Cone(\eta)\otimes \unital$.  A \emph{counital idempotent} is an object $\counital$ together with a map $\e:\counital\rightarrow \one$ such that $\Cone(\e)\otimes \counital\simeq 0\simeq Q\otimes \Cone(\e)$.  Often times we refer to $\eta$ and $\e$ as a unital and counital idempotent, respectively.
\end{definition}
It is clear that if $\eta:\one\rightarrow \unital$ is a unital idempotent, then $(\Cone(\eta)\ip{1},\unital)$ is a pair of complementary idempotents.  Similarly, if $\e:\counital\rightarrow \one$ is a counital idempotent, then $(\counital,\Cone(\e))$ is a pair of complementary idempotents.  Conversely, if $(\counital,\unital)$ is a pair of complementary idempotents, then $\counital$ and $\unital$ are counital, respectively unital, idempotents, with structure maps given by \eqref{eq:PQtriang}.  We remark that unital and counital idempotents would be called closed and open idempotents in \cite{BD14}, respectively.

Note that $(P_n,\eta)$ from Definition/Theorem \ref{defthm-Pn} is a unital idempotent in $\K^-(\SBim_n)$.  Many important properties of $P_n$ follow from general arguments.
\begin{definition}
Let $\im P_n$ and $\ker P_n$ denote the full subcategories of $\KC^-(\SBim_n)$ consisting of complexes $C$ such that $P_n\otimes C\simeq C$ and $P_n\otimes C\simeq 0$, respectively. 
\end{definition}

\begin{proposition}\label{prop-uniqueAndCentral}
Recall Definition \ref{def-N}. We have
\begin{enumerate}
\item The uniqueness part of Definition/Theorem \ref{defthm-Pn} follows from the existence part.
\item $P_n$ is central: for any $C\in\K^-(\SBim_n)$ there is a homotopy equivalence $C\otimes P_n\simeq P_n \otimes C$ which is natural in $C$.
\item ${}_\perp{{\NC}} = {{\NC}}_\perp=\im P_n$.
\item $\NC=\ker P_n$.
\end{enumerate}
\end{proposition}
\begin{proof}
This is an immediate consequence of Corollary 4.29 in \cite{Hog17a}.
\end{proof}
Once we prove the existence of $P_n$ in \S \ref{subsec-projAsLimit}, will exclusively use the notation $\im P_n$ and $\ker P_n$, instead of $\NC_\perp = {}_\perp {\NC}$ and $\NC$.

In the remainder of this section we discuss the notion of relative unital idempotents.  The following is Theorem 4.24 in \cite{Hog17a}.
\begin{theorem}\label{thm:fundTheoremOfIdemp}
Let $\IB\in \AS$ and $\PB\in \AS$ be unital idempotents with unit maps $\eta_{\IB}:\one\rightarrow \IB$ and $\eta_{\PB}:\one\rightarrow \PB$.  Then the following are equivalent.
\begin{enumerate}
\item $\IB\otimes \PB\simeq \PB$.
\item $\PB\otimes \IB\simeq \PB$.
\item there exists a map $\nu:\IB\rightarrow \PB$ such that $\eta_{\PB} = \nu\circ \eta_{\IB}$.
\end{enumerate}
If either of these equivalent conditions is satisfied then $\nu$ is unique, and $\Id_{\PB}\otimes \nu$ and $\nu\otimes \Id_{\PB}$ are isomorphisms in $\AS$.
\end{theorem}
There is a similar statement for counital idempotents.
\begin{definition}\label{def:relativeIdempts}
If either of the equivalent conditions of Theorem \ref{thm:fundTheoremOfIdemp} are satisfied, then we say that $\nu:\IB\rightarrow \PB$ gives $\PB$ the structure of a \emph{unital idempotent relative to $\IB$}.
One defines the notion of a \emph{counital idempotent relative to $\idempotent$} in a similar way.
\end{definition}

\begin{proposition}\label{prop-PnRelPk}
Fix integers $a,b,c\geq 0$ let $\iota:\K^-(\SBim_b)\rightarrow \K^-(\SBim_{a+b+c})$ denote the functor $\iota(C)=\one_a\sqcup C\sqcup \one_c$.  Then $P_{a+b+c}$ has the structure of a unital idempotent relative to $\iota(P_b)$.
\end{proposition}
\begin{proof}
Note that $\iota$ is a monoidal functor, $\iota(\eta_b): \iota(\one_b)\cong \one_{a+b+c}\rightarrow \iota(P_b)$ is a unital idempotent.   Note that $\Cone(\iota(\eta_b))$ is homotopy equivalent to a complex built from $\iota(B_1),\cdots,\iota(B_{b-1})$, and tensoring with $P_{a+b+c}$ annihilates each of these.  Thus, $\Cone(\iota(\eta_b))\otimes P_{a+b+c}\simeq 0$, which implies that $\iota(P_b)\otimes P_{a+b+c}\simeq P_{a+b+c}$.
\end{proof}

In \S \ref{subsec:endGeneral} we state some important impliciations of this relationship between $P_n$ and $P_k$ ($1\leq k\leq n$).

\subsection{Rouquier complexes}
\label{subsec-rouquier}
We recall Rouquier's categorification of the group homomorphism $\Br_n\rightarrow \HB_n^\times$.  Rouquier complexes play an important role in our construction of $P_n$.  Certain expressions will look nicer if we use a slightly different normalization for Rouquier complexes than is standard.   

\begin{definition}\label{def-rouquierCx}
Let $\sigma_i$ denote the $i$-th generator in the $n$-strand braid group.  That is, $\sigma_i$ is a positive crossing between the $i$ and $i+1$ strands.  Let $\sigma_i\inv$ denote the inverse crossing.  Define the \emph{Rouquier complexes}
\begin{eqnarray*}
%\Big\llbracket\: \sspic{plusCrossing}{24pt}\Big\rrbracket
F(\sigma_i) &:=& \Big(\begin{diagram}[small]0 & \rTo & {B_i(1)} & \rTo & \underline{R} & \rTo & 0\end{diagram}\Big)\\
F(\sigma_i\inv) &:=& \Big(\begin{diagram}[small]0 & \rTo & \underline{R} & \rTo & {B_i(-1)}  & \rTo & 0\end{diagram}\Big)
\end{eqnarray*}
where we have underlined the degree zero chain bimodules, and the maps are the dot maps.  If $\b=\sigma_{i_1}^{\e_1}\cdots\sigma_{i_r}^{\e_r}$ is expressed as a product in the generators (on $n$ strands), then we will let $F({\b})=F(\sigma_{i_1}^{\e_1})\otimes_R \cdots\otimes_R F(\sigma_{i_r}^{\e_r})$ denote the Rouquier complex associated to ${\b}$.  It is well known that $F({\b})$ depends only on the braid $\b$ (and not on the expression of $\b$ in terms of generators) up to canonical isomorphism in $\K(\SBim_n)$ (see \cite{EKh10} for details).
\end{definition}
\begin{remark}\label{rmk-R}
The bimodule $R$ and its shifts do not appear as a direct summand of any nontrivial tensor product $B_{i_1}\cdots B_{i_r}$.  Thus, our normalization for Rouquier complexes ensures that $F(\b)$ has a unique copy of $R$ appearing with grading shift $(0)\ip{0}$ for each braid $\b$.  If $\b$ is a positive braid, then the inclusion $R\rightarrow F(\b)$ is a chain map.  If, instead, $\b$ is a negative braid, then the projection $F(\b)\rightarrow R$ is a chain map.  In either case, the cone on the map $F(\b)\leftrightarrow R$ lies in $\NC_n=\ker P_n$.
\end{remark}

In the remainder of this section we include some assorted basic results that will be used later.

\begin{proposition}\label{prop-ringActionOnCrossing}
Let $s\in S_n$ be the permutation that swaps $k,k+1$.  Then left multiplication by $x_i\in R$ is chain homotopic to right multiplication by $s(x_i)$, as endomorphisms of $F(\sigma_k)$.
\end{proposition}
\begin{proof}
If $i\not \in \{k,k+1\}$ then the claim is obvious from the fact that $F(\sigma_k)$ is obtained from $F(\sigma_1)$ by extension of scalars from $\ring[x_k,x_{k+1}]$ to $\ring[x_1,\ldots,x_n]$.  In the remaining two cases, we have a homotopy:
\[
\begin{diagram}
0 & \rTo & B_k(3) & \rTo & R(2)& \rTo & 0\\
& \ldTo & & \ldTo^{h_1} & & \ldTo & \\
0 & \rTo &B_k(1) &\rTo & R(0)& \rTo & 0
\end{diagram}.
\]
The top row is $F(\sigma_k)(2)$, the bottom row is $F(\sigma_k)$, and the diagonal arrows indicate the homotopy whose only nonzero component $h_1$ is the dot map.  An easy calculation shows that $dh+hd$ sends $c\mapsto x_1c-c x_2 = cx_1-x_2c$.  This proves the claim.
\end{proof}

\begin{proposition}\label{prop-Bsliding}
$F_uB_s\simeq B_sF_u$ when $u$ and $s$ are distant.  If $s$ and $t$ are adjacent then $F_{st} B_s= B_tF_{st}$.  %In particular $J_n B_s \simeq B_s J_n$ for each simple reflection $s\in S_n\subset S_{n+1}$.
\end{proposition}
\begin{proof}
The first statement is clear.  The second can be checked by direct computation, and we leave this as an exercise.  In any case, this calculation is implicit in the proof of invariance of Rouquier complexes under the braid relation $F_sF_tF_s\simeq F_tF_sF_t$.
\end{proof}
The fact that $F_{st}B_s\simeq B_t F_{st}$ will be referred to by saying that $B_s$ \emph{slides past crossings}.

\begin{proposition}\label{prop-Psliding}
Set $X:=F(\sigma_n\cdots \sigma_2\sigma_1)$ and $Y:=F(\sigma_n\inv\cdots \sigma_2\inv\sigma_1\inv)$.  These are objects in $\K^-(\SBim_{n+1})$.  Then
\[
X(P_n\sqcup \one_1) \simeq (\one_1\sqcup P_n)X \ \ \ \ \ \ \ \ \ \ \ \ \text{and} \ \ \ \ \ \ \ \ \ \ \ \ \ \ \ Y(P_n\sqcup \one_1) \simeq (\one_1\sqcup P_n)Y.
\]
\end{proposition}
This proposition implies that $P_n$ can be used to define a link homology theory (in fact, $\Sym^n$-colored Khovanov-Rozansky homology; see \S \ref{subsec:coloredHomology}).
\begin{proof}
We prove only the first of these; the second is similar.  Since $B_i$ slides through crossings, we have $XB_i\simeq B_{i-1}X$ for $2\leq i\leq n$.  Since $(P_n\sqcup \one_1)$ kills $B_1,\ldots,B_{n-1}$ it follows that $(P_n\sqcup \one_1)XB_i\simeq 0$ for $2\leq i\leq n$, hence
\[
(P_n\sqcup \one_1)X(\one_1\sqcup P_n)\simeq (P_n\sqcup \one_1)X
\]
by Proposition \ref{prop-PnRelPk}.  A similar argument shows that 
\[
(P_n\sqcup \one_1)X(\one_1\sqcup P_n)\simeq X(\one_1\sqcup P_n)
\]
from which the Lemma follows.
\end{proof}

Finally, we have two absorption properties which will be useful:

\begin{proposition}[Projectors absorb Rouquier complexes]\label{prop-PabsorbsRouquier}
The complex $P_n$ absorbs Rouquier complexes: $F(\b)P_n\simeq P_n\simeq P_n F(\b)$ for all braids $\b$.
\end{proposition}
\begin{proof}
Immediately follows from the fact that $P_nB_i\simeq 0\simeq B_iP_n$ for all $1\leq i\leq n-1$ and the definition of the elementary Rouquier complexes $F(\sigma_i^{\pm})$ (Definition \ref{def-rouquierCx}).
\end{proof}

\begin{proposition}[Bott-Samelsons absorb Rouquier complexes]\label{prop-BabsorbsRouquier}
$B_sF_s\simeq B_s(2)\ip{-1}\simeq F_sB_s$.
\end{proposition}
\begin{proof}
Recall that $F_s$ is the complex $(B_s(1)\rightarrow \underline{\one}).$  A basic property of Soergel bimodules states that $B_sB_s\cong B_s(1)\oplus (-1)$.  In fact, the dot map $\phi:B_s(1)\rightarrow R$ is such that $\phi\otimes \Id_{B_s}:B_sB_s(1)\rightarrow B_s$ is the projection onto a direct summand.  Thus, tensoring $F_s$ on the right with $B_s$ gives a complex of the form
\[
B_sF_s\ = \  B_s(2)\oplus B_s(0)\buildrel\matrix{\ast &\Id}\over \longrightarrow \underline{B_s}(0).
\]
The contractible summand $(B_s\buildrel \Id\over \rightarrow B_s)$ can be canceled from the complex by Gaussian elimination.  The result is $F_sB_s\simeq B_s(2)\ip{-1}$, as claimed.  A similar argument takes care of $B_sF_s$.
\end{proof}

\subsection{The projector as an infinite Rouquier complex}
\label{subsec-projAsLimit}
In this section we complete the proof of Definition/Theorem \ref{defthm-Pn} by constructing $P_n$ as the homotopy colimit of Rouquier complexes.  The idea that categorified idempotents can be realized as infinite braids is due to Rozansky \cite{Roz10a}, and was utilized by D.~Rose \cite{Rose14} in constructing the $\sl_3$ categorified idempotents and by S.~Cautis \cite{Cau12} in categorifying all of the $\sl_n$ idempotents.  

In \S \ref{sec-KRhomology} we define triply graded link homology in terms of Rouquier complexes.  The results of this section imply that the triply graded link homology of the $(n,m)$ torus links stablizes as $m\to\infty$, and that the stable limit can be computed from $P_n$.

We first recall the notion of homotopy colimit, or \emph{mapping telescope}.  Recall the shift in homological degree $[1]=\ip{-1}$.

\begin{definition}\label{def-homotopyLimit}
Let $\{A_i\:|\: i\in \Z_{\geq 0}\}$ be chain complexes over an additive category $\CS$, and $f_i:A_i\rightarrow A_{i+1}$ chain maps.  Suppose that the infinite direct sum $\bigoplus_{k\geq 0} A_k$ exists in $\K(\CS)$.  The \emph{homotopy colimit} of the directed system $\{A_k, f_k\}_{i\in \Z_{\geq 0}}$ is the mapping cone
\[
\hocolim A_k := \Cone\Big(\bigoplus_{k\geq 0}A_k \buildrel S-\Id\over\longrightarrow \bigoplus_{k\geq 0} A_k\Big),
\]
where $S$ is the shift map, whose components are given by the $f_i$.
\end{definition}

Below, we will find it convenient to represent certain kinds of total complexes diagrammatically.  For instance, we will denote the homotopy colimit $\hocolim A_k$ in the following way:
\begin{equation}\label{eq-hocolimDiagram}
\left(
\begin{diagram}
A_0[1] && A_1[1] && A_2 [1] && A_3[1] && \cdots\\
\dTo_{-\Id} &\rdTo^{f_0} & \dTo_{-\Id} &\rdTo^{f_1} & \dTo_{-\Id} &\rdTo^{f_2} & \dTo_{-\Id} &\rdTo^{f_3} & \\
A_0 && A_1 && A_2 && A_3 && \cdots
\end{diagram}\right)^\oplus.
\end{equation}
Note that each arrow is homogeneous of homological degree $+1$.  The complex $\hocolim A_k$ is obtained by taking the direct sum of all of the complexes in the above diagram.  The differential on $C$ is by definition given by the differentials internal to each term (with the differential on $A_i[1]$ being $-d_A$ by convention), as well as by the labelled arrows. 

As another example of the diagram notation, let $\CS$ be an additive category and $A_k\in \K(\CS)$ complexes so that the direct sum $\bigoplus_{k\in \Z}A_k$ exists.  We will write
\begin{equation}\label{eq-plusCx}
C \ = \  \ (\cdots \buildrel[1]\over\longrightarrow A_i  \buildrel[1]\over\longrightarrow   A_{i+1}  \buildrel[1]\over\longrightarrow  \cdots)^\oplus 
\end{equation}
whenever $C\in \K(\CS)$ is a chain complex which equals $\bigoplus_{k\in \Z}A_k$ as a graded object, and whose differential is represented by a $\Z\times \Z$-lower triangular matrix with $d_{A_k}$'s along the diagonal.   The shifts $[1]$ are there to remind us that all arrows are actually chain maps $A_i\rightarrow A_j[1]$, for $i\leq j$.

\begin{remark}\label{rmk-hocolim}
The homotopy colimit $\hocolim A_k$ of a directed system $\{A_k, f_k\}$ is isomorphic to a complex of the form
\begin{equation}\label{eq-hocolimEq2}
\hocolim A_k \cong \Big(\cdots \buildrel[1]\over\longrightarrow \Cone(f_2)\buildrel[1]\over\longrightarrow \Cone(f_1)\buildrel[1]\over\longrightarrow  \Cone(f_0) \buildrel[1]\over\longrightarrow  A_0\Big)^{\oplus}.
\end{equation}
This is simply a reassociation of the complex represented by the diagram (\ref{eq-hocolimDiagram}), in which each pair of terms connected by a diagonal arrow is rewritten as a copy of $\Cone(f_k)$.
\end{remark}

\begin{proposition}\label{prop:simultSimp}
Suppose we are given complexes $A_i$, $B_i$, and homotopy equivalences $f_i:A_i\simeq B_i$ ($i\in \Z$).  If $A_i=B_i=0$ for $i\gg 0$, then any complex (\ref{eq-plusCx}) is homotopy equivalent to a complex of the form
\[
C \simeq (\cdots \buildrel[1]\over\longrightarrow B_i \buildrel[1]\over\longrightarrow B_{i+1} \buildrel[1]\over\longrightarrow \cdots)^\oplus.
\]
\end{proposition}
\begin{proof}
Follows from homotopy invariance of mapping cones together with a straightforward colimit argument.%We only give a sketch.  Without loss of generality, we assume that $A_i=B_i=0$ for $i>0$. Now $C$ is a colimit of its truncations.  In other words, there exists a family of complexes $C_0,C_1,\ldots$ such that $C_0=A_0$, and the $C_k$ satisfy a recursion of the form $C_k=\Cone(A_{-k}\ip{1}\rightarrow C_{k-1})$, and so that $C$ is an honest colimit of the $C_k$ as $k\to\infty$.  Note that the $C_k$ fit into distinguished triangles of the form\[C_{k-1} \rightarrow C_{k} \rightarrow A_{-k} \rightarrow C_{k-1}[1]\]
%for all $k\geq 1$.  This is called a Postnikov system in triangulated categories literature (see, for instance \cite{Orl97}).  Assume we are given complexes $B_i$ and homotopy equivalences $f_i:A_i\rightarrow B_i$.  We now inductively construct complexes $D_0,D_1,\ldots,$ and homotopy equivalences $\phi_k:C_k\rightarrow D_k$.  Initially, set $D_0=B_0$, and $\phi_0=f_0$.  Suppose we have constructed $D_k$ and a homotopy equivalence $\phi_k:C_k\rightarrow D_k$.  Then simplification of mapping cones (Proposition \ref{prop:mappingConeSimp}) says that we may modify the terms of $C_{k+1}=\Cone(A_{-1-k}\ip{1}\rightarrow C_{k})$, obtaining $D_{k+1}=(\Cone(B_{-1-k}\ip{1}\rightarrow D_{k})$ and an equivalence $\phi_{k+1}:C_{k+1}\rightarrow D_{k+1}$ such that the square
%
%\[\begin{diagram}C_{k}&\rTo & C_{k+1}\\ \dTo^{\phi_{k-1}} && \dTo^{\phi_k} \\ D_{k}&\rTo & D_{k+1} \end{diagram} \]
%
%commutes, where the horizontal maps are the structure maps associated to mapping cones.  In this way we construct a directed system $\{D_k\}$ and a map of directed systems $\{C_k\}\rightarrow \{D_k\}$.  The map of directed systems induces a map of colimits, which is a homotopy equivalence.
\end{proof}

We have similar notions with direct sum replaced everywhere by direct product.  Correspondingly, if we have two families of complexes $A_i$ and $B_i$ such that $A_i\simeq B_i$ and $A_i=B_i=0$ for $i\ll 0$, then any complex
\[
(\cdots \buildrel[1]\over\longrightarrow A_i \buildrel[1]\over\longrightarrow A_{i+1} \buildrel[1]\over\longrightarrow \cdots)^\Pi
\]
is homotopy equivalent to a complex of the form
\[
(\cdots \buildrel[1]\over\longrightarrow B_i \buildrel[1]\over\longrightarrow B_{i+1} \buildrel[1]\over\longrightarrow \cdots )^\Pi.
\]
We will make such infinite simplifications throughout this subsection. 

\begin{proposition}\label{prop-hocolimBdd}
Let $\CS$ be an additive category, and $\{A_k,f_k\}_{k=0}^\infty$ a directed system of complexes $A_i\in \K^-(\CS)$.  Assume that $\Cone(f_k)\simeq C_k$, where $C_k$ is a complex which is supported in homological degrees $c_k$, for some sequence such that $c_k\to -\infty$ as $k\to \infty$.  Then $\hocolim A_k$ is an object of $\K^-(\CS)$, up to equivalence.
\end{proposition}
The content here is that one need not postulate the existence of infinite direct sums, and one need not leave the world of semi-infinite complexes in order to have a well-defined hocolim, in this situation.
\begin{proof}
We can include $\CS$ (fully faithfully) into a category $\CS^\oplus$ which contains infinite direct sums.  Then $\hocolim A_k$ exists as an object of $\K(\CS^\oplus)$.  Each term $\Cone(f_k)$ in (\ref{eq-hocolimEq2}) can be replaced by $C_k$ up to homotopy equivalence, and we obtain a complex $L\simeq \hocolim A_k\in \K^-(\CS^\oplus)$.  The hypotheses ensure that the infinite direct sum $\bigoplus_{k=0}^\infty C_k$ is finite in each homological degree, hence $L$ is isomorphic to an object of $\K^-(\CS)$.
\end{proof}

The directed system of interest to us is the following:  let $x=\sigma_{n-1}\cdots \sigma_2\sigma_1$ denote the $n$-stranded braid which is the positive lift of the standard $n$-cycle $(n,n-1,\ldots,2,1)\in S_n$ (written in cycle notation).  Note that the braid closure of $x^k$ is the $(n,k)$-torus link.   Let $X=F(x)$ denote the associated Rouquier complex.  We claim that the $X^{\otimes k}$ can be made into a direct system with homotopy colimit $P_n$.   So set $A_k:=X^{\otimes k}$, and $A_0:=\one$.  We have a map $f_0:\one\rightarrow X$ which is the inclusion of the $\one$ summand of $X$ (see Remark \ref{rmk-R}).  Set $f_k:=f_0\otimes X^{\otimes k}:A_k\rightarrow A_{k+1}$.  We wish to prove:

\begin{theorem}\label{thm-PasHocolim}
The homotopy colimit $P_n$ of the directed system $\{A_k, f_k\}_{k\in \Z_{\geq 0}}$ lies in $\K^-(\SBim_n)$.  If $\eta:\one=A_0\rightarrow P_n$ is the structure map, then $(P_n,\eta)$ satisfy the axioms of Definition/Theorem \ref{defthm-Pn}.
\end{theorem}

We first need two Lemmas:
\begin{lemma}\label{lemma-Xlemma}
If $C\in\mathcal{N}$ (Definition \ref{def-N}) is supported in homological degrees $\leq \ell$, then $X^{\otimes 2k}\otimes C$ is homotopy equivalent to a complex which is supported in homological degrees $\leq \ell-2$.
\end{lemma}
\begin{proof}
The complex $X^{\otimes 2n}$ is the Rouquier complex associated to the full twist braid $x^{2n}$.  For each $1\leq i\leq n-1$, the full twist can be factored as $x^{2n} = \sigma_i^2y$ for some positive braid $y$.  Thus, the corresponding Rouquier complex factors as $X^{\otimes 2n}\simeq  F_i^2 F(y)$, where $F(y)$ is supported in non-positive homological degrees.  Thus, the fact that $B_i$ absorbs $F_i$ up to a shift of the form $(2)\ip{-1}$ (Proposition \ref{prop-BabsorbsRouquier}) implies that $X^{\otimes 2n}\otimes B_i$ is homotopy equivalent to a complex supported in homological degrees $\leq -2$ for each $1\leq i\leq n-1$.

Now suppose $C\in {\NC_n}$ is supported in homological degrees $\leq \ell$.  Up to equivalence we may assume that the chain bimodules $C_i$ are direct sums of shifts of bimodules of the form $B_{i_1}\cdots B_{i_r}$, and $C_i=0$ for $i>\ell$ by hypothesis.  The Lemma follows by applying the above simplification to each term $X^{\otimes 2n}C_i\ip{i}$ of $X^{\otimes 2n}C$.
\end{proof}

\begin{lemma}\label{lemma-nullColimits}
Let $\{A_k, f_k\}_{k=0}^\infty$ be a directed system of complexes over an additive category $\CS$ and let $c_k\in \Z$ be so that $A_k$ is supported in homological degrees $\leq c_k$, and $c_k\to -\infty$ as $k\to \infty$.  Then $\hocolim A_k$ exists in $\K(\CS)$, and is contractible by a dual version of Proposition \ref{prop:simultSimp}.
\end{lemma}
\begin{proof}
Let $\{A_k,f_k\}$ and $c_k$ be as in the hypotheses.  Then the infinite direct sum $\bigoplus_{k=0}^\infty A_k$ is finite in each homological degree, hence is equal to the categorical direct product $\prod_{k=0}^\infty A_k$ in $\K(\CS)$ by standard arguments.  This shows that $\hocolim A_k$ exists, and is given by same diagram as in (\ref{eq-hocolimDiagram}), with $\oplus$ replaced by direct product.  Then
\[
\hocolim A_k = (\Cone(\Id_{A_0})\rightarrow \Cone(\Id_{A_1})\rightarrow \Cone(\Id_{A_2})\rightarrow \cdots )^\Pi,
\]
which is contractible.
\end{proof}

\begin{proof}[Proof of Theorem \ref{thm-PasHocolim}]
By definition, $X$ is the tensor product of complexes $F_s:=(B_s\rightarrow \underline{\one})$, hence can be written as a mapping cone
\[
X = \Cone(D_0\rightarrow \one)
\]
where $D_0\in \NC$ is some complex which is supported in non-positive homological degrees.  By definition $f_0:\one\rightarrow X$ is the inclusion of the degree 0 chain bimodule, hence $\Cone(f_0)\simeq D_0[1]$.  Recall that $f_k:X^{\otimes k}\rightarrow X^{\otimes k+1}$ is obtained by tensoring $f_0$ on the left with $X^{\otimes k}$.  Thus,
\[
\Cone(f_k)\cong X^{\otimes k}\otimes \Cone(f_0)\simeq X^{\otimes k}\otimes D_0[1].
\]
By Lemma \ref{lemma-Xlemma}, we can find complexes $D_k\simeq X^{\otimes k} \otimes D_0$ supported in homological degrees $d_k$, such that $d_{k+2n}\leq d_k - 2$ for all $k\geq 0$.  Since $\Cone(f_k)\simeq D_k[1]$, Proposition \ref{prop-hocolimBdd} says that there is a well defined complex $Q\simeq \hocolim_k A_k \in \K^-(\SBim_n)$, of the form
\[
Q = \Big(\cdots \rightarrow D_{-2}[1]\rightarrow D_{-1}[1]\rightarrow D_{0}[1]\rightarrow \one\Big)^{\oplus}.
\]

Let $\nu:\one \rightarrow Q$ denote the inclusion of the $\one$ summand above.  We claim that $(Q,\nu)\simeq (P_n,\eta_n)$.  Clearly $\Cone(\nu)\in \mathcal{N}$ since $D_k\in \NC$ for all $k$.  It remains to check that $Q$ kills $B_s$, for each simple transposition $s\in S_n$.  For this, note that $Q\otimes B_s$ is the homotopy colimit of the directed system $X^{\otimes k}\otimes B_s$.  Iterated application of Lemma \ref{lemma-Xlemma} says that $X^{\otimes k}\otimes B_s$ is equivalent to a complex $C_k$ which is supported in homological degrees $\leq c_k$ such that $c_k\to -\infty$ monotonically as $k\to \infty$. Thus, Lemma \ref{lemma-nullColimits} yields
\[
Q\otimes B_s \simeq \hocolim X^{\otimes k}\otimes B_s \simeq \hocolim C_k \simeq 0.
\]
Therefore, $(Q,\nu)$ satisfy the axioms of Definition/Theorem \ref{defthm-Pn} as claimed.
\end{proof}

\section{Derived categories and triply graded homology}
\label{sec-KRhomology}
Khovanov-Rozansky homology is defined in terms of Rouquier complexes, together with a functor $\HHH:\K^-(\SBim_n)\rightarrow \ring\mod^{\Z\times \Z\times \Z}$, where this latter category is the category of triply graded $\ring$-abelian groups.   The purpose of this section is to introduce, for each pair of complexes $C,D\in \K(\SBim_n)$ a triply graded space of homs $\DHom(C,D)$ which contains the usual group of chain maps mod homotopy as the degree $(0,0,0)$ part.  We then show that $\HHH(P_n)$ is isomorphic to the triply graded ring of endomorphisms $\DEnd(P_n)$.  We establish some techniques which will be used in \S \ref{sec:structure} to compute $\DEnd(P_n)$.  Our main tools are an adjunction isomorphism (Proposition \ref{prop-adjointness1}) and a Markov relation for $\HHH$ (Proposition \ref{prop:markov}).  Although the proofs use facts about derived categories, it does not require much background to understand and use these tools.

\subsection{Derived categories}
\label{subsec-HHH}
Whenever the index $n$ is understood, we will let $\xx$ denote the list of variables $x_1,\ldots,x_n$, and similarly for $\yy$.   We will write $R=\ring[\xx]$, and we will identify $R^e=R\otimes_\ring R$ with $\ring[\xx,\yy]$ via $x_i\otimes 1\mapsto x_i$ and $1\otimes x_i\mapsto y_i$.  We will regard graded $(R,R)$-bimodules as graded $R^e$-modules.  Let $\CS_n$ denote the bounded derived category of graded $R^e$-modules, which is equivalent to the homotopy category of finite complexes of \emph{free} graded $R^e$-modules, the equivalence being given by free resolution.  For each $X\in \CS_n$, let $X(i,j)$ denote the result of shifting $X$ up by $i$ in $q$-degree, and up by $j$ in derived (homological) degree.

Observe that $\CS_n$ is a monoidal category in the usual way, via the derived tensor product $\buildrel L\over\otimes$.  There is a fully faithful inclusion $\SBim_n\rightarrow \CS_n$.  Since Soergel bimodules are free as right (or left) $R$-modules (see \cite{EKh10}), the tensor product of Soergel bimodules $M\otimes N$ coincides with the derived tensor product up to isomorphism in $\CS_n$.  Thus, the inclusion $\SBim_n\rightarrow \CS_n$ is a monoidal functor.  By abuse of notation, we sometimes denote triangulated closure of $\SBim_n$ in $\CS_n$ by $D^b(\SBim_n)$.  

\begin{definition}\label{def:HH}
If $M$ is a graded $(R,R)$-bimodule, let $\HH^{\bullet,\bullet}(M)$ denote the bigraded abelian group $\HH^{\bullet,\bullet}(M)=\bigoplus_{i,j\in \Z}\Hom_{\CS_n}(R(i,j),M)$.  Note that $\HH^{i,j}\cong \Ext^j_{R^e}(R,M(-i))$ is the usual Hochschild cohomology group of $R$ with coefficients in $M(-i)$.  When we wish to emphasize the index $n$, we will write $\HH^{\bullet,\bullet}(R_n;M)$.  In this paper, we exclusively use Hochschild cohomology (never Hochschild homology), so we will drop the superscripts from $\HH^{\bullet,\bullet}$, writing only $\HH$.
\end{definition}

Now, let $\DS_n:=\K^-(\CS_n)$ denote the homotopy category of complexes over $\CS_n$.  The inclusion $\SBim_n\rightarrow \CS_n$ extends to a fully faithful, monoidal functor $\K^-(\SBim_n)\rightarrow \K^-(\CS)$.  Henceforth, we will regard objects of $\K^-(\SBim_n)$ also as objects of $\DS_n$.

Note that an object $C\in\DS_n$ is triply graded: $C=\bigoplus_{i,j,k\in \Z} C_{i,j,k}$.  The first grading $i$ is the bimodule degree, also called $q$-degree.  The grading $j$ will be called derived or Hochschild grading, and is the homological grading in $\CS_n$.  The grading $k$ is called homological grading, as it is the homological grading in $\K^-(\CS_n)$.  We let $(a,b)\ip{c}$ denote the shift in tridegree, so that $(C(a,b)\ip{c})_{i,j,k}=C_{i-a,j-b,k-c}$.  See \S \ref{subsec-notation} for more.

\begin{definition}\label{def-DHom}
For objects $A,B\in \K^-(\SBim_n)$, let $\DHom_n(A,B)$ denote the triply graded abelian group
\[
\DHom_n^{i,j,k}(A,B)  := \Hom_{\DS_n}(A(i,j)\ip{k}, B).
\]
For $f\in \DHom^{i,j,k}(A,B)$, let $|f|=j+k$ denote the total homological degree. Similarly, $\DEnd_n(A)=\DHom_n(A,A)$. As always, when the index $n$ is understood, we will drop it from the notation.
\end{definition}

\begin{definition}\label{def-HHH}
If $C$ is a complex of graded $(R,R)$-bimodules, let $\HH(C)$ denote the complex of bigraded abelian groups obtained by applying the functor $\HH$ to the chain bimodules of $C$.  Set $\HHH(C):=H(\HH(C))$.
\end{definition}

It is clear from the definitions that $\HHH(C)\cong \DHom(R,C)$.

\begin{remark}\label{rmk-derivedSBim}
If $C,D,\in\K^-(\SBim_n)$, then the usual space of chain maps $C\rightarrow D$ modulo homotopy can be recovered as the subspace of $\DHom(C,D)$ consisting of elements of degree $(0,0,0)$.  In general, an element $f\in \DHom^{0, j,0}(C,D)$ should be thought of as a chain map whose components are higher exts: $f_m\in \Ext^j(C_m, D_{m})$.  These are to be regarded modulo homotopies $h$ whose components are exts of the same degree: $h_m\in\Ext^j(C_m,D_{m-1})$.
\end{remark}

Let $\b$ be an $n$-stranded braid, and $L=\hat{\b}\subset \R^3$ the oriented link obtained by connecting the strands of $\b$ (in a planar fashion).  In \cite{Kh07}, Khovanov showed that $\HHH(F(\b))$ depends only on the isotopy class of $L$, up to isomorphism and overall grading shift, and coincides with \emph{Khovanov-Rozansky homology} \cite{KR08} up to regrading.  This fact, together with Theorem \ref{thm-PasHocolim} gives says that $\HHH(P_n)$ is a limit of triply graded homologies of $(n,k)$ torus links.  See Corollary \ref{cor:invariance} for the precise normalization of $\HHH$ that yields a link invariant, and see Corollary \ref{cor:stabilization} for the precise statement regarding the stabilization of triply graded homology of torus links.

In \S \ref{subsec:endGeneral} we show that $\HHH(P_n)\cong \DHom(R,P_n)\cong \DEnd(P_n)$, hence $\HHH(P_n)$ has the structure of a triply graded ring, which is in fact graded commutative.  We will compute this ring in \S \ref{subsec-endP}, using tools developed here and in subsequent sections.

\subsection{An adjunction isomorphism}
\label{subsec-adjunction}
The point of the next few subsections is to take the first steps toward the computation of $\DEnd(P_n)$.  In this section we construct a pair of adjoint functors $I:\DS_{n-1}\rightarrow \DS_n$ and $\Tr:\DS_n\rightarrow \DS_{n-1}$.  The functor $\Tr$ should be thought of as ``partial Hochschild cohomology,'' in that  $\HHH(C)\cong \DHom_{\DS_n}(R,C)\cong \DHom_{\DS_0}(\ring,\Tr^n(C))$ for all $C\in \DS_n$.

\begin{definition}\label{def-Tr}
Define a pair of functors $I:\CS_{n-1}\rightarrow \CS_n$ and $\Tr:\CS_n\rightarrow \CS_{n-1}$ as follows.  For any graded $R_{n-1}^e$-module $M$, let $I(M)\cong C[x_n]$ denote the  graded $R_n^e$ module obtained by extending scalars.  More precisely:
\[
I(M) = \Big(R_n^e/(x_n-y_n)\Big)\otimes_{R_{n-1}^e} M.
\]
Since $R_n^e/(x_n-y_n)$ is free as an $R_{n-1}^e$-module, $I$ gives a well-defined functor $\CS_{n-1}\rightarrow \CS_n$.  For any $N\in \CS_n$, let $\Tr(N)$ denote the total complex:

%We will think of $\Lambda$ as a complex of graded $(R_n^e,R_{n-1}^e)$-modules.  Let $I:\CS_{n-1}\rightarrow \CS_n$ denote the functor $I(C)=\Lambda\otimes_{R_{n-1}^e} C$, and let $\Tr:\CS_n\rightarrow \CS_{n-1}$ denote the functor $\Tr(D)=\Homg_{R_n^e}(\Lambda, D)$ (the notation $\Homg$ is explained in \S \ref{subsec-notation}).  Note that $\Tr(D)$ can also be regarded as the total complex of the bicomplex
\[
\Tr(N) \ \ \cong \ \ \Tot( 0 \rightarrow {N} \buildrel {x_n-y_n}\over\longrightarrow  N(-2)\rightarrow 0),
\]
The grading shifts involved above are characterized by the fact that forgetting the differential on $\Tr(N)$ yields $N\oplus N(-2)[1]$.  We call $\Tr$ the \emph{partial trace} or \emph{partial Hochschild cohomology}.
\end{definition}

\begin{proposition}\label{prop-adjointness1}
The functors $(I,\Tr)$ from Definition \ref{def-Tr} form an adjoint pair.  In other words, there is an isomorphism
\[
\Hom_{\CS_n}(I(M),N)\cong\Hom_{\CS_{n-1}}(M,\Tr(N))
\]
which is natural in $N\in\CS_n$ and $M\in \CS_{n-1}$.
\end{proposition}
%We point out that naturality here means that the adjunction isomorphism $A:\Hom_{\CS_n}(I(M),N) \rightarrow \Hom_{\CS_{n-1}}(M,\Tr(N))$ satisfies  $A(f\circ g\circ I(h)) = \Tr(f)\circ A(g)\circ h$ whenever this makes sense.
\begin{proof}
Let $\Lambda$ denote the complex
\[
0\rightarrow R_n^e(2)  \buildrel x_n - y_n \over \longrightarrow \underline{R_n^e} \rightarrow 0.
\]
Clearly $\Lambda$ is homotopy equivalent (as a complex of graded $R_{n-1}^e$-modules) to $R_{n-1}^e[x_n]$, which implies that $I(M)\cong \Lambda\otimes_{R_{n-1}^e}C$ for all $M\in \CS_{n-1}$. Dually, we have $\Tr(N)\cong\Homg_{R_n^e}(\Lambda,N)$ for all $N\in \CS_n$, where the notation $\Homg$ is explained in \S \ref{subsec-notation}.  The Proposition follows from the usual Hom-tensor adjunction.
\end{proof}

Extending $I$ and $\Tr$ to complexes gives functors $\DS_{n-1}\leftrightarrow \DS_n$, and naturality of the  adjunction isomorphism gives:

\begin{corollary}\label{cor:adjunctionComplexes}
There is an isomorphism $\DHom_n(I(C),D)\cong \DHom_{n-1}(C,\Tr(D))$ which is natural in $C\in \DS_{n-1}$ and $D\in \DS_n$.  In particular, $\HHH(D)\cong \DHom_0(\ring,\Tr^n(D))$ for all $D\in\DS_n$.
\end{corollary}
\begin{proof}
Note that $R_n = I^n(\ring)$.  Repeated application of the adjunction isomorphism gives $\DHom_n(I^n(\Z),D)\cong \DHom_0(\ring,\Tr^n(\ring))$ as claimed.
\end{proof}
 
 The following is immediate:
\begin{lemma}\label{lemma-traceIsLinear}
The partial trace $\Tr:\DS_{n} \rightarrow \DS_{n-1}$ is $\DS_{n-1}$ ``bilinear.''  That is to say, there is a natural isomorphism
\[
\Tr(I(M)\otimes C\otimes I(N)) \cong M \otimes \Tr(C) \otimes N
\]
for all complexes $M,N\in \DS_{n-1}$, $C\in \DS_{n}$.\qed
\end{lemma}

In \S \ref{subsec-endP} we use the setup of this section to study $\DEnd(P_n)$.

\subsection{The Markov move}
\label{subsec:markov}

Now we introduce the tools which will be used in simplifying $\Tr(C)$ for certain complexes $C\in \DS_n$ (particularly $P_n$).  It will be helpful below to visualize an object $M\in \DS_{n}$ as a labelled box with $n$ strands attached to the top and bottom, and $\HHH(M)$ as the diagram obtained by connecting all the loose strands in a planar fashion:
\[
M\ = \  \begin{minipage}{.5in}
\labellist
\pinlabel $M$ at 8 13
\endlabellist
\includegraphics[scale=1]{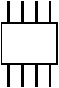}
\end{minipage}
\hskip 1.3in
\HHH(M) \ = \  \begin{minipage}{.5in}
\labellist
\pinlabel $M$ at 8 18
\endlabellist
\includegraphics[scale=1]{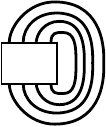}
\end{minipage}.
\]
Then tensor product $\otimes_R$ corresponds to vertical stacking, and external tensor product $\otimes_\ring:\SBim_i\times\SBim_j\rightarrow \SBim_{i+j}$ corresponds to horizontal juxtaposition.  The identity bimodule $R_n$ is denoted by $n$ parallel strands, and the objects $B_i, F_i, F_i\inv\in \DS_n$ will be denoted diagrammatically by
\[
F_i \ = \ \sspic{1in}{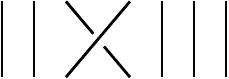}, \ \ \ \ \ \ F_i\inv \ = \ \sspic{1in}{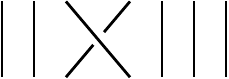}, \ \ \ \ \ \ B_i \ = \ \sspic{1in}{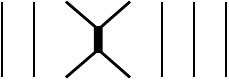} ,
\]
pictured in the case $n=7$, $i=3$.

The functors $I(M)$ and $\Tr(N)$ are indicated diagrammatically by
\[
I(M) \ =
\begin{minipage}{.4in}
{\labellist
\pinlabel $M$ at 8 13
\endlabellist
\includegraphics[scale=1]{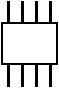}
}\
\includegraphics[scale=1]{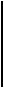}
\end{minipage}
\hskip1in
\Tr(N) \ = \
\begin{minipage}{.4in}
\labellist
\pinlabel $N$ at 8 13
\endlabellist
\includegraphics[scale=1]{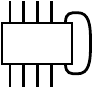}
\end{minipage}.
\]

Then Lemma \ref{lemma-traceIsLinear} can be stated diagrammatically as
\[
\begin{minipage}{.4in}
\labellist
\pinlabel $N$ at 8 10
\pinlabel $C$ at 8 24.5
\pinlabel $M$ at 8 39.5
\endlabellist
\includegraphics[scale=1]{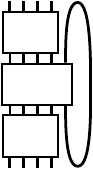}
\end{minipage}
\ \cong \ 
\begin{minipage}{.4in}
\labellist
\pinlabel $N$ at 8 10
\pinlabel $C$ at 8 24.5
\pinlabel $M$ at 8 39.5
\endlabellist
\includegraphics[scale=1]{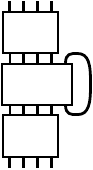}
\end{minipage}
\]

\begin{definition}
For any complex $C\in \DS_n$, set $q^ia^jt^k C:=C(i,j)\ip{k}$.  If $f(q,a,t)$ is a Laurent power series in $q,a,t$ with positive integer coefficients, then let $f(q,a,t)C$ denote the corresponding direct sum of shifted copies of $C$.
\end{definition}
In this section we prove the following diagram relations, which we refer to loosely as Markov moves:
\begin{subequations}
 \begin{equation}\label{eq:trRpic}\spic{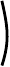}\sspic{.2in}{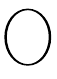} \ \cong \   \frac{1+q^{-2}a}{1-q^2} \ \ \spic{shortStrand}\end{equation}
\begin{equation}\label{eq:trBpic}\sspic{.26in}{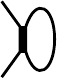} \ \cong \   \frac{q+q^{-3}a}{1-q^2} \ \ \spic{shortStrand}\end{equation}
\begin{equation}\label{eq-trFpic}\sspic{.35in}{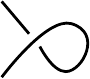} \ \simeq   \ \ \spic{shortStrand}\end{equation}
\begin{equation}\label{eq-trFinvpic}\sspic{.35in}{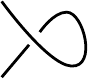}\  \simeq  \ q^{-4}at \ \ \spic{shortStrand}\end{equation}
\end{subequations}
A more precise statement follows:

\begin{proposition}[Markov move]\label{prop:markov}
Let $F_i^{\pm}=F(\sigma_i^\pm)$ denote the Rouquier complexes associated to the elementary braid generators.  For each $n\geq 1$, we have
\begin{subequations}
 \begin{equation}\label{eq:trR}\Tr(R_{n})\cong R_{n-1}[x_{n}]\oplus R_{n-1}[x_{n}](-2,1)\ip{0}\end{equation}
\begin{equation}\label{eq:trB}\Tr(B_{n-1})\cong R_{n-1}[x_{n}](1,0)\ip{0}\oplus R_{n-1}[x_{n}](-3,1)\ip{0}\end{equation}
\begin{equation}\Tr(F_{n-1})   \simeq    R_{n-1} \label{eq-trF}\end{equation}
\begin{equation}\Tr(F_{n-1}\inv )  \simeq   R_{n-1}(-4,1)\ip{1} .\label{eq-trFinv}\end{equation}
\end{subequations}
The first two are isomorphisms in $\CS_{n-1}$, and the second two are isomorphisms in $\DS_{n-1}$.  When used in combination with Lemma \ref{lemma-traceIsLinear} these yield, for instance, $\HHH(R_n;CF_n)\cong \HHH(R_{n-1};C)$ and $\HHH(R_n;CF_n\inv)\cong \HHH(R_{n-1};C)(4,1)\ip{1}$ for all $C\in \K^-(\SBim_{n-1})$.
\end{proposition}
This follows from results in \cite{KR08b} and \cite{Kh07}.  For the reader's convenience, and since the precise grading shifts play an important role in this paper, it is worth including more details.
\begin{definition}\label{def:koszulCx}
Let  $r_1,\ldots,r_k\in R^e$ a collection of elements. The Koszul complex associated to $\mathbf{r}=(r_1,\ldots,r_k)$ is the differential bigraded $R^e$-algebra $\Lambda[\theta_1,\ldots,\theta_k]$ with differential $d(\theta_i)=r_i$.  If $R$ is graded and $\deg(r_i)=q^{c_i}$, then make $K$ bigraded by $\deg(\theta_i)=a\inv q^{c_i}$, where $a\inv$ signifies homological degree $-1$.  We also use the notation $\ring[\xx,\yy,\theta_1,\ldots,\theta_k]$ for Koszul complexes, the understanding being that the odd variables anti-commute and square to zero.
\end{definition}

If the sequence $r_1,\ldots,r_k$ is regular, then the projection $K\rightarrow R^e/I$ is a quasi-isomorphism of complexes of graded $R^e$-modules, where $I\subset R^e$ is the ideal generated by the $r_i$.  In some instances these quasi-isomorphisms can be regarded as homotopy equivalences.
\begin{lemma}\label{lemma:exclusion}
Let $K$ denote the Koszul complex associated to the element $y-x\in \ring[x,y]$.  Then there is a $\ring[x]$-equivariant homotopy equivalence $K\rightarrow \ring[x]$.
\end{lemma}
\begin{proof}
As an abelian group we have $K=\ring[x,y]\oplus \ring[x,y]$ with differential $d(f,g) = (0,(y-x)f)$.  Any element $g\in \ring[x,y]$ can be written uniquely as $g(x,y)=(y-x)g_1(x,y) + g_2(x)$, where $g_2(x)=g(x,x)\in \ring[x]$.  Thus, the subcomplex of $K$ generated by elements of the form $(f,0)$ and $(0,(y-x)g)$, with $f,g\in\ring[x,y]$, is a contractible summand.  Contracting it yields $K\simeq \ring[x]$.
\end{proof}
%\begin{proof}
%By definition we have $K=\ring[x,y,\theta]$ with $d(\theta)=y-x$.  There is a $\ring[x,y]$-linear map $\pi:K\rightarrow \ring[x]$ sending $1\mapsto 1$ and $\theta\mapsto 0$.  There is a $\ring[x]$-linear right inverse $\sigma:\ring[x]\rightarrow K$ sending $1\mapsto 1$.  The composition $\sigma\circ \pi:K\rightarrow K$ is the $\ring[x]$-linear map sending $y^k\mapsto x^k$ and $\theta y^k\mapsto 0$.  Let $h$ denote the $\ring[x]$-linear homotopy defined inductively by $h(1)=0$ and $h(y^k)=y^{k-1}\theta+xh(y^{k-1})$ for $k\geq 1$.  An easy induction shows that $dh+hd=\Id_K-\sigma\circ \pi$.  This proves the lemma.
%\end{proof}

\begin{proof}[Proof of Proposition \ref{prop:markov}]
First, note that Equation (\ref{eq:trR}) is follows easily from the definitions.  

For the remainder of the proof, assume that $n=2$.  The result for arbitrary $n$ reduces immediately to this one. Let $B=B_s=R\otimes_{R^s}R(-1)$ denote the nontrivial Bott-Samelson bimodule.  Let $F=F(\sigma_1)$ denote the elementary Rouquier complex, and $F\inv=F(\sigma_1\inv)$ its inverse.  %The dot map $B(1)\rightarrow R$ sends $1\otimes 1$, while the other dot map $R\rightarrow B(-1)$ sends $1\mapsto x_1\otimes 1-1\otimes x_2$.

Let $K$ denote the Koszul complex $\Lambda[\theta_1,\theta_2]$, where $d(\theta_1)=x_1+x_2-y_1-y_2$ and $d(\theta_2)=x_1x_2-y_1y_2$.  The degrees are necessarily $\deg(\theta_1)=q^{2}a\inv$ and $\deg(\theta_2)=q^4a\inv$.  This complex is quasi-isomorphic to its homology, which is $B(1,0)$. Let $\xi$ be an odd variable of degree $q^2a\inv$ with differential $d(\xi)=x_2-y_2$.  From the definition of $\Tr$ it is clear that $\Tr(K)\cong K\otimes_{R^e}\Lambda[\xi](-2,1)$.  In other words,
\[
\Tr(K)(2,-1) \cong\Lambda[\theta_1,\theta_2,\xi].
\]
Let us change variables, introducing $\theta_1'=\theta_1-\xi$ and $\theta_2'=\theta_2-(x_1-y_2)\xi - y_2\theta_1$.  Note that $d(\theta_1')=x_1-x_2$, and $d(\theta_2)=0$.  An application of Lemma \ref{lemma:exclusion} says we may cancel $\xi$ and set $x_2=y_2$, obtaining
\[
\Tr(K)(2,-1)\cong \ring[x_1,x_2,y_1,\theta_1',\theta_2'].
\]
After factoring out the Koszul complex $\ring[x_1,y_1,\theta_1']\cong R_1$, and recognizing that $\theta_2'$ has degree $q^4a\inv$, we see that $\Tr(K)(2,-1)\cong K(R_1)\oplus K(R_1)(4,-1)$.  Recalling that $K\cong B(1,0)$, we have proven Equation (\ref{eq:trB}).

Now, for Equations (\ref{eq-trF}) and (\ref{eq-trFinv}) we may appeal to the usual invariance of Khovanov-Rozansky homology under the Reidemeister I move.  Specifically, setting the formal variable $a=0$ in Propositions 4 and 5 in \cite{KR08b}, we see that Equations (\ref{eq-trF}) and (\ref{eq-trFinv}) hold up to some degree shift.  To see that we have the correct degree shifts is a simple exercise, using our computations above.  For instance, $F\inv$ is the complex $0\rightarrow \underline{R}(0,0)\rightarrow B(-1,0)\rightarrow 0$.  applying $\Tr$ we obtain
\[
\Tr(F\inv) \ \ \simeq \ \ \Big(0\rightarrow (R_2\oplus R_2(-2,1))\rightarrow (R_2\oplus R_2(-4,1))  \rightarrow 0 \Big)
\]
where $R_2=\ring[x_1,x_2]$ (though we forget the action of $x_2$).  The component of the differential $R_2\rightarrow R_2(-4,1)$ is zero since there are no negative degree exts between $R^e$-modules. In order to have $\Tr(F\inv)\simeq R_1$ up to a shift, it must be that the unshifted $R_2$ terms cancel one another.  Up to isomorphism the remaining component of the differential has to be a unit multiple of $x_1-x_2$, since  multiplication by $x_1-y_2\simeq x_1-x_2$ is a null-homotopic endomorphism of $\Tr(F\inv)$.  Then an application of Lemma \ref{lemma:exclusion} says that the resulting complex is homotopy equivalent to $R_1(-4,1)$.  There is an additional shift of $\ip{1}$, because $B(-1,0)$ is the degree 1 chain bimodule of $F\inv$.  This proves (\ref{eq-trFinv}).  A similar argument proves (\ref{eq-trF}).
\end{proof}

For any $C\in\K^-(\SBim_n)$, let $\PC_c(q,a,t)$ denote the Poincar\'e series of $\HHH(C)$.  If $\b$ is a braid, let $\PC_\b(q,a,t)=\PC_{F(\b)}(q,a,t)$.  We record the following for the sake of posterity:
\begin{corollary}\label{cor:invariance}
If $\b$ is a braid, let $e(\b)$ denote the braid exponent (signed number of crossings).  The normalized Poincar\'e series $(t^{1/2}a^{1/2}q^{-2})^{e(\b)-n}\PC_\b(q,a,t)$ is an invariant of the braid closure $\hat{\b}$.  If we introduce a formal variable $\a = t^{1/2}a^{1/2}q\inv$, then the resulting link invariant assigns the value $\frac{\a\inv+\a/t}{q\inv-q}$ to the unknot.  In general, setting $t=-1$ recovers the usual, unnormalized Homfly invariant in variables $\a,q$.\qed
\end{corollary}
Thus, in order to properly normalize $\HHH(F(\b))$ to obtain an actual link invariant, it is necessary to introduce half-integral homological and Hochschild degrees.

\subsection{Stabilization}
\label{subsec:stabilization}
The results of this section will not be used elsewhere in this paper, but they may be of independent interest.  Recall the notation and results of \S \ref{subsec-projAsLimit}.  We now find a bound on the stable range of $\HHH(X^{\otimes k})$.  The maps of our directed system give rise to $f_{i,j}:X^{\otimes j}\rightarrow X^{\otimes i}$ .  By definition, $f_{i+k,j+k}$ is obtained from $f_{i,j}$ by tensoring on the right with $X^{\otimes k}$.  Thus,
\[
\Cone(f_{i+k,i}) \cong \Cone(f_{k,0})X^{\otimes i}.
\]
Since $f_{k,0}:\one\rightarrow X^{\otimes k}$ is the inclusion of the unique $\one$ summand, $\Cone(f_{k,0})\in\ker P_n$, and this complex is supported in homological degrees $<0$ up to homotopy.  Tensoring with $X^{\otimes i}$ and using Lemma \ref{lemma-Xlemma}, we see that $\Cone(f_{i+k,i})$ is supported in homological degrees $< -2d$, where $d=\floor{i/n}$.  Thus, the ``difference'' between $X^{\otimes i+k}$ and $X^{\otimes i}$ is small when $i$ is large.  That is, the $X^{\otimes i}$ form a categorical analogue of a Cauchy sequence in the sense of Rozansky \cite{Roz10a}.  It follows that, up to equivalence, $X^{\otimes i+k}$ and $X^{\otimes i}$ agree in homological degrees $\geq -2d$.

\begin{proposition}\label{prop:stabilization}
The map $f_{i+1,i}$ induces a map $\HHH(X^{\otimes i})\rightarrow \HHH(X^{\otimes i+1})$ which is an isomorphism through homological degrees $>-2\floor{i/n}$.
\end{proposition}
\begin{proof}
Consider the exact triangle
\[
X^{\otimes i} \rightarrow X^{\otimes i+1} \rightarrow \Cone(f_{i+1,i})\rightarrow X^{\otimes i}\ip{-1}.
\]
The category $\DS_n$ is triangulated, from which it follows that $\HHH(-) = \DHom(R,-)$ is a homological functor.  Thus, there is a long exact sequence
\[
\HHH(X^{\otimes i}) \rightarrow \HHH(X^{\otimes i+1}) \rightarrow \HHH(\Cone(f_{i+1,i}))\rightarrow \HHH(X^{\otimes i}\ip{-1}).
\]
The third term is zero in degrees $\geq -2\floor{i/n}$, from which the proposition follows.
\end{proof}

Considering the normalization of triply graded homology in Corollary \ref{cor:invariance} yields:
\begin{corollary}\label{cor:stabilization}
The triply graded homology of the $(n,k)$ torus link is supported in homological degrees $i$, where $-nk+k-n\leq 2i\leq nk-k-n$.  The stable range is those degrees $i$ with  $nk-k-n-4\floor{n/k} < 2i\leq nk-k-n$.\qed
\end{corollary}

\section{Structure of the projector}
\label{sec:structure}

Our goal in this section is to prove that the triply graded ring of endomorphisms of $P_n$ is the superpolynomial ring $\ring[u_1,\ldots,u_n,\xi_1,\ldots,\xi_n]$, and to interpret the endomorphisms corresponding to the even generators $u_k$.  Essentially, $u_k$ captures the fact that $P_k$ can be constructed as a certain kind of periodic complex built out of $P_{k-1}$.

\subsection{Statement of the Theorems}
\label{subsec-statement}
We wish to prove:
\begin{theorem}\label{thm-endP}
For each integer $n\geq 1$, the triply graded algebra of endomorphisms (Definition \ref{def-DHom})  of $P_n$ satisfies $\DEnd(P_n)\cong \ring[u_1,\ldots,u_n,\xi_1,\ldots,\xi_n]$, where the tridegrees are $\deg(u_k)=q^{2k}t^{2-2k}$ and $\deg(\xi_k)=q^{2k-4}t^{2-2k}a$.  This is an isomorphism of triply graded algebras.
\end{theorem}
Here, the $u_i$ are even variables and the $\xi_i$ are odd variables (so the $\xi_i$ anti-commute and square to zero).

\begin{remark}\label{rmk-groundRingAction}
The image of the ground ring $R=\ring[x_1,\ldots,x_{n}]$ in $\DEnd(P_{n})$ under the isomorphism from Theorem \ref{thm-endP} is precisely $\ring[u_1]\cong R/J$, where $J\subset R$ is the ideal generated by the differences $x_i-x_{i+1}$ ($1\leq i\leq n-1$).  The fact that left multiplication by $x_i$ is homotopic to left multiplication by $x_{i+1}$ follows from Proposition \ref{prop-ringActionOnCrossing} and the fact that $F_iP_{n}\simeq P_{n}$.  The fact that left multiplication by $x_i$ is homotopic to right multiplication by $x_i$ follows from general facts about unital idempotents (see Proposition \ref{prop:rhoProp} with $\PB=P_n$ and $\IB=\one_n$).
\end{remark}

Our computation of this endomorphism ring requires some insight into the structure of $P_n$.  This will be provided by Theorem \ref{thm-eigenconesAreOrtho} (below).  We first set up some notation:

\begin{definition}\label{def-JMandCycles}
Let $X_n=F(\sigma_{n-1}\cdots \sigma_2\sigma_1)$ and $Y_n=F(\sigma_{n-1}\inv\cdots \sigma_2\inv\sigma_1\inv)$ denote the Rouquier complexes associated to the positive and negative braid lifts of the $n$-cycle $(n,n-1,\ldots,2,1)$. These are objects of $\K^-(\SBim_n)$.  Let $J_n:=X_nY_n\inv$ denote the Rouquier complex associated to the \emph{Jucys-Murphy braid}.
\end{definition}
 % Below, it may be useful to picture $I(P_{n-1})$ and $I(P_{n-1})$ diagrammatically as in (\ref{eq-introDiagrammatics}).  The following is crucial:

\begin{notation}
For any $1\leq k\leq n$, let $P_k$ denote the complex $P_k\sqcup \one_{n-k}\in \K^-(\SBim_n)$, by abuse of notation.
\end{notation}
\begin{theorem}\label{thm-eigenconesAreOrtho}
There exist chain maps $\b_{(n-1,1)},\b_{(n)}\in \Homg(P_{n-1}Y_n, P_{n-1}X_n)$ with degrees
\[
\deg(\b_{(n-1,1)}) = q^{2n}t^{2-2n}\hskip1in \deg(\b_{(n)})=q^0t^0,
\]
such that $\Cone(\b_{(n-1,1)})\in \im P_n$ and $\Cone(\b_{(n)})\in \ker P_n$ (Definition \ref{def-N}).
\end{theorem}
\begin{remark}
The reason for this notation will become clear in future work of the author with Ben Elias.  Specifically, it will be shown that for each partition $\l$ of $n$, there is a very special chain map $\a_\l\in \Homg(\one_n,\FT_n)$, where $\FT_n=J_2\cdots J_n$ is the full twist complex.  Tensoring on the left with $P_{n-1}\sqcup \one_1$ gives a family of maps $\a_\l'\in\Homg(P_{n-1},P_{n-1}J_n)$.  Tensoring on the right with $Y_n$ then gives maps $\b_\l\in \Homg(P_{n-1}Y_n,X_n)$.  The maps in the above theorem are then special cases of this construction.
\end{remark}

We do not know of a proof of Theorem \ref{thm-eigenconesAreOrtho} which does not rely on Theorem \ref{thm-endP}.  To avoid a circular argument, we have no choice but to prove Theorem \ref{thm-eigenconesAreOrtho} and Theorem \ref{thm-endP} simultaneously, by induction on $n\geq 1$.  This is accomplished in subsequent sections.

%This section is organized as follows.  In \S \ref{} we show how one constructs $P_n$ from $P_{n-1}$, given Theorem \ref{thm-eigenconesAreOrtho}.  The construction makes it clear that $\ring[u_1,\ldots,u_n]$ acts on $P_n$.  In \S \ref{} we find an interpretation of the maps $\b_{(n-1,1)},\b_{(n)}$ of Theorem \ref{} in terms of ``crossing-change'' maps.  Finally, in \S \ref{} we define a homological obstruction to defining $\b_{(n-1,1)},\b_{(n)}$.  Much later, in \S \ref{}, we will show that this obstruction vanishes.  

\subsection{$P_n$ quasi-isomorphisms, and $\b_{(n)}$}
\label{subsec:easyMap}
In this section we show that $\b_{(n)}:P_{n-1}Y_n\rightarrow P_{n-1}X_n$ is easy to construct.    First, we introduce some terminology.  Let us say that a morphism $f:C\rightarrow D$ in $\KC^-(\SBim_n)$ is a \emph{$P_n$ quasi-isomorphism} if $P_n\otimes f$ is a homotopy equivalence (equivalently, $f\otimes P_n$ is a homotopy equivalence).  It is straightforward to see that the tensor product and composition of $P_n$ quasi-isomorphisms is a $P_n$ quasi-isomorphism.

\begin{lemma}
A chain map $f:C\rightarrow D$ is a $P_n$ quasi-isomorphism if and only if $\Cone(f)\in \ker(P_n)$. 
\end{lemma}
\begin{proof}
Observe: $P_n\otimes \Cone(f)\cong \Cone(P_n\otimes f)$, which is contractible if and only if $P_n\otimes f$ is a homotopy equivalence from properties of triangulated categories. 
\end{proof}

There is a canonical chain map $\one\rightarrow F(\sigma_i)$ which includes $\one$ as the degree zero chain group.  For the negative crossing there is a canonical projection $F(\sigma_i\inv)\rightarrow \one$.  These maps are $P_n$ quasi-isomorphisms.  Tensoring them gives canonical $P_n$ quasi-isomorphisms $\one\rightarrow F(\b)$ and $F(\b\inv)\rightarrow \one$ for any positive braid $\b$.

Similarly for any $a+b+c=n$, let $P':=\one_a\sqcup P_b\sqcup \one_c$, and let $\eta':\one_{n}\rightarrow P'$ be the unit map (induced from the unit map of $P_k$).  Then $\eta'$ is a $P_n$ quasi-isomorphism. 

Now, $Y_n$ and $X_n$ are Rouquier  complexes associated to a negative braid and a positive braid, respectively, so we have a canonical $P_n$ quasi-isomorphism $Y_n\rightarrow \one_n \rightarrow X_n$.

\begin{definition}\label{def:easyEigenmap}
 Let $\b_{(n)}:P_{n-1}Y_n\rightarrow P_{n-1}X_n$ be the map obtained by tensoring the $P_n$ quasi-isomorphism $X\rightarrow Y$ with $P_{n-1}$.
\end{definition}

\begin{remark}\label{rmk:bn}
Let $\e_1:Y\rightarrow \one_n$ and $\e_2:\one_n\rightarrow X$ denote the canonical maps, so that $\b_{(n)}=(P_{n-1}\e_2)\circ(P_{n-1}\e_1)$.  Thus, the composition
\[
\begin{diagram}
P_nP_{n-1} & \rTo^{(P_nP_{n-1}\e_2)\inv} & P_nP_{n-1}Y_n & \rTo^{P_n \b_{(n)}} & P_nP_{n-1}X_n & \rTo^{(P_nP_{n-1}\e_1)\inv} & P_nP_{n-1}
\end{diagram}
\]
is homotopic to the identity, hence $\b_{(n)}$ is a $P_n$ quasi-isomorphism.
\end{remark}

\subsection{Constructing the projector}
\label{subsec-constructionOfPn}
The main idea of our construction was outlined in \S \ref{subsec-introIdea} of the introduction.

Let $u_n$ be a formal indeterminate of bidegree $q^{2n}t^{2-2n}$.  If $C\in \DS_n$ is any chain complex, then we will write
\begin{equation}\label{eq:polyringTensor}
C[u_n] = C\otimes_\Z \Z[u_n] = C\oplus C(2n,0)\ip{2-2n}\oplus C(4n,0)\ip{4-4n}\oplus \cdots.
\end{equation}
This is a convenient way of expressing certain direct sums.  Note that for $n>1$ the above countable direct sum is finite in each homological degree, hence is isomorphic to a direct product in $\DS_n$.  Any chain complex of this form comes equipped with an endomorphism of degree $(2n,0,2-2n)$, given by multiplication $u_n$.  In terms of (\ref{eq:polyringTensor}), multiplication by $u_n$ is simply the obvious inclusion $C[u_n](2n,0)\ip{2-2n}\rightarrow C[u_n]$.

\begin{construction}\label{construction-Pn}
  Suppose that $\b_{(n-1,1)},\b_{(n)}\in \DHom(P_{n-1}Y_n,P_{n-1}X_n)$ as in Theorem \ref{thm-eigenconesAreOrtho} have been constructed. Consider the chain map
\[
\Psi \ :\ \ring[u_n]\otimes_\ring P_{n-1}Y_n(2-2n)\ip{2n}  \ \ \longrightarrow  \ \ \ring[u_n]\otimes_\ring P_{n-1}X_n,
\]
defined by $\Psi = 1\otimes \b_{(n-1,1)} - u_n\otimes \b_{(n)}$, which can also be represented by the following semi-infinite diagram of chain complexes and chain maps:
\vskip 7pt
\begin{equation}\label{eq:bigPn} \ \ \ \ \ \ \ \ \ \ \ \ \ 
\begin{minipage}{3.8in}
\begin{tikzpicture}[baseline=-0.25em]
\matrix (m) [ssmtf,column sep=5em, row sep=5em]
{\sspic{.4in}{PY}(2n)\ip{2-2n} & \sspic{.4in}{PX}(0)\ip{0} \\
\sspic{.4in}{PY}(4n)\ip{4-4n} &  \sspic{.4in}{PX}(2n)\ip{2-2n} \\
\sspic{.4in}{PY}(6n)\ip{6-6n} &  \sspic{.4in}{PX}(4n)\ip{4-4n} \\
\cdots & \cdots\\};
\node (11) at (m-1-1)[xshift=2em,yshift=-1em] {};
\node (21) at (m-2-1)[xshift=2em,yshift=-1em] {};
\node (31) at (m-3-1)[xshift=2em,yshift=-1em] {};
\node (22) at (m-2-2)[xshift=-5em,yshift=1em] {};
\node (32) at (m-3-2)[xshift=-5em,yshift=1em] {};
\node (42) at (m-4-2)[xshift=-5em,yshift=1em] {};
\path[->,font=\scriptsize]
(m-1-1) edge node[auto] {} (m-1-2)
(11) edge node[auto] {} (22)
(m-2-1) edge node[auto] {} (m-2-2)
(21) edge node[auto] {} (32)
(m-3-1) edge node[auto] {} (m-3-2)
(31) edge node[auto] {} (42);
\end{tikzpicture}
\end{minipage}
\end{equation}
in which the horizontal arrows are given by $\b_{(n-1,1)}$ and the diagonal arrows are $-\b_{(n)}$, with the appropriate shift functors applied.  The sign on $\b_{(n)}$ is for later convenience (see Proposition \ref{prop:periodicityMap}).   Set $Q:=\Cone(\Psi)\in\K^-(\SBim_{n})$ and let $\nu:R\rightarrow Q$ denote the inclusion $R\rightarrow P_{n-1}X_n$, followed by the inclusion $P_{n-1}X_n\rightarrow Q$ (as the top right term in Diagram (\ref{eq:bigPn})).
\end{construction}

\begin{proposition}\label{prop-theConstructionWorks}
The pair $(Q,\nu)$ from Construction \ref{construction-Pn} satisfies the axioms for $(P_n,\eta_n)$ in Definition/Theorem \ref{defthm-Pn}. % Furthermore, the periodicity map (multiplication by $u_n$) is homotopic to the $U_n$ from Definition \ref{def:U}. 
\end{proposition}
\begin{proof}
Recall the notation of Construction \ref{construction-Pn}.  Observe that $Q=\Cone(\Psi)$ is built out of shifted copies of $\Cone(\b_{(n-1,1)})$, hence is an object of $\im P_n$ since $\Cone(\b_{(n-1,1)})$ is, by hypothesis.  Thus, $Q$ satisfies axiom (P1) of Definition/Theorem \ref{defthm-Pn}.

Now we must show that $(Q,\nu)$ satisfies axiom (P2).  We have $\Cone(\b_{(n)})\in \ker(P_n)$ by construction, hence the $R$ summands of $P_{n-1}Y_n$ and $P_{n-1}X_n$ cancel one another in $\Cone(\b_{(n)})$.  After canceling all of these $R$-summands from $Q=\Cone(\Psi)$, it follows that there is a unique $R$ summand, in homological degree 0, coming from the $P_{n-1}X_n$ term in the top right of Diagram (\ref{eq:bigPn}).   By definition, $\nu$ is the inclusion of the degree zero chain group $R\rightarrow P_{n-1}X_n$, followed by the inclusion $P_{n-1}X_n\rightarrow Q$.  Thus, $\nu$ is the inclusion of the unique $R$ summand of $Q$.  It follows that $\Cone(\nu)\in \ker(P_n)$, as claimed.
\end{proof}

%Finally, we sketch the argument that $U_n\in\DEnd(P_n)$ coincides with the periodicity map of $P_n$ up to homotopy.   Since $\b_{(n)}:P_{n-1}Y_n\rightarrow P_{n-1}X_n$ is the identity on the $R$-summands, it follows that tensoring $\b_{(n)}$ with $P_n$ and simplifying yields the identity map on $P_n$.  Tensoring $\b_{(n-1,1)}$ on the right $P_n$ and simplifying yields $U_n$ by definition. Thus, tensoring Diagram (\ref{eq:bigPn}) on the right with $P_n$ and simplifying yields an equivalent diagram:
%\[
%\begin{diagram}[width=4em]
%1\otimes P_n\ip{2n,2-2n} & \rTo^{1\otimes U_n} & 1\otimes P_n \\
%& \rdTo^{-u_n\otimes \Id} & \\
%u_n\otimes P_{n}\ip{2n,2-2n}  & \rTo^{1\otimes U_n} & u_n\otimes P_{n}\\
%& \rdTo^{-u_n\otimes\Id} & \\
%u_n^2\otimes P_{n}\ip{2n,2-2n}& \rTo^{1\otimes U_n} & u_n^2\otimes P_{n}\\
%& \rdTo^{-u_n\otimes \Id} & \\
%\cdots & \cdots & \cdots
%\end{diagram}.
%\]
%This corresponds to a simplification of the mapping cone $P_nP_n =\Cone(\Psi)P_n = \Cone(\Psi\Id_{P_n})$ as in Proposition \ref{prop:mappingConeSimp}, and results in an equivalent complex $P_n'$ which is the mapping cone on the map $u_n\otimes \Id - 1\otimes U_n$.  The homotopy equivalence $P_nP_n\simeq P_n'$ can be seen to commute with the evident $u_n$ and $U_n$ actions up to homotopy.  The fact that $u_n\simeq U_n$ is now clear, since these maps are clearly homotopic as endomorphisms of $P_n'$.

\begin{example}\label{example-baseCase}
Recall that $F(\sigma_1) = (B(1)\rightarrow \underline{R})$ and $F(\sigma_1\inv)=(\underline{R}\rightarrow B(-1))$.  We define $\b_{(1,1)}:F(\sigma_1\inv)(4)\ip{-2}\rightarrow F(\sigma_1)$ in such a way that
\[
\Cone(\b_{(1,1)})  \ \ \simeq   \ \ \begin{diagram}  R(4) & \rTo & B_1(3) & \rTo^{x_2\otimes 1 - 1\otimes x_2} & B_1(1) & \rTo & \underline{R} \end{diagram}
\]
where the first and last maps are the canonical maps.  The proof that $\Cone(\b_{(1,1)})B_1\simeq 0$ follows along the same lines as in the case of $P_2$ in Proposition \ref{prop:P2}.

Proposition \ref{prop-theConstructionWorks} states that $P_2$ can be reconstructed from infinitely many copies of this four-term complex.  Doing so explicitly is an easy exercise, and is left to the reader.
\end{example}

In subsequent sections we use this explicit description of $P_n$ to compute $\HHH(P_n)$.  But first recall some general arguments that will aid the computation.

%For the remainder of the section, fix an integer $n\geq 1$.  All of our complexes will regarded as objects of $\K^-(\SBim_{n+1})$.

%The complex $\Cone(\b_{(n-1,1)})$ is the Soergel analogue of the complex $Q_n$ constructed in \cite{H14a}; the reader is invited compare the expression (\ref{eq-coneB1}) for $\Cone(\b_{(n-1,1)})$ with the expression (3.6) for the complex $Q_n$ in \cite{H14a}.  The above construction of $P_{n+1}$ from $\Cone(\b_{(n-1,1)})$ is similar to the construction of the categorified Jones-Wenzl projector from $Q_n$ in \cite{H14a}.

\subsection{Generalities on endomorphism rings of idempotents}
\label{subsec:endGeneral}
Let $(\AS,\otimes,\one)$ be a triangulated monoidal category (for example $\AS=\KC^-(\SBim_n)$ or $\AS=\DS_n$.).  Let $\eta_{\PB}:\one\rightarrow \PB$ be a unital idempotent.

\begin{proposition}\label{prop:resProp}
For any $X\in \AS$ and any $Y\in \AS$ such that $\PB\otimes Y\cong Y$, precomposing with $\eta_\PB\otimes \Id_X$ gives an isomorphism
\[
\Hom_{\AS}(\PB\otimes X, Y)\cong \Hom_{\AS}(\one\otimes X,Y).
\]
In particular, precomposing with $\eta_\PB$ gives an isomorphism $\End(\PB)\rightarrow \Hom(\one,\PB)$. 
\end{proposition}
Compare with Proposition 4.16 in \cite{Hog17a}.    Let's deduce some consequences of this proposition.  We may regard $\eta_{\PB}\otimes \Id_\PB$ and $\Id_\PB\otimes \eta_{\PB}$ as elements of $\Hom(\PB,\PB\otimes \PB)$.  Precomposing with $\eta_{\PB}$ gives an isomorphism
\[
\eta^\ast:\Hom(\PB,\PB\otimes \PB)\rightarrow \Hom(\one,\PB\otimes \PB).
\]
Then $\eta_{\PB}\otimes \Id_\PB = \Id_\PB\otimes \eta_{\PB}$, since the both become $\eta_\PB\otimes \eta_\PB$ upon applying $\eta_\PB^\ast$. 
\begin{proposition}\label{prop:mapHop}
We have $f\otimes \Id_\PB = \Id_\PB\otimes f\in \End(\PB\otimes \PB)$ for all $f\in \End(\PB)$.
\end{proposition}
\begin{proof}
 By definition of unital idempotent, $\Id_\PB \otimes \eta_\PB$ is invertible.  We have an isomorphism $\End(\PB\otimes \PB)\rightarrow \End(\PB)$ sending $g\mapsto (\eta_\PB\otimes \Id_\PB)\inv\circ g\circ (\eta_\PB\otimes \Id_\PB)$.  We show that $\Id\otimes f$ and $f\otimes \Id$ have the same image under this isomorphism.  The image of $\Id_\PB\otimes f$ under this isomorphism is
\[
(\eta\otimes \Id_\PB)\inv\circ (\Id_\PB\otimes f)\circ (\eta_\PB\otimes \Id_\PB) = (\eta\otimes \Id_\PB)\inv \circ (\eta\otimes \Id_\PB) \circ (\Id_{\one}\otimes f),
\]
which equals $f$.  Since $\eta\otimes \Id_\PB=\Id_\PB\otimes \eta$, a similar argument shows that $f\otimes \Id_\PB$ also has image $f$ under the above isomorphism.  This completes the proof.
\end{proof}

\begin{proposition}\label{prop:EndProp}
The algebra $\End(\PB)$ is a commutative algebra.  There is an algebra structure on $\Hom(\one,\PB)$ given by $\a\cdot \b:=\a\otimes \b\in \Hom(\one^{\otimes 2},\PB^{\otimes 2})\cong\Hom(\one,\PB)$.   Precomposing with $\eta$ gives an isomorphism of algebras $\End(\PB)\cong \Hom(\one,\PB)$. 
\end{proposition}
\begin{proof}
Let $f,g\in \End(\PB)$ be arbitrary.  Observe:
\[
(f\otimes \Id_{\PB}) \circ (\Id_{\PB}\otimes g) =  (\Id_{\PB}\otimes g)\circ (f\otimes \Id_{\PB}).
\]
Since $f\otimes \Id_\PB= \Id_\PB\otimes f$, the left-hand side is $\Id_{\PB}\otimes(f\circ g)$ and the right hand side is $\Id_\PB\otimes (g\circ f)$.  Tensoring with $\PB$ gives an isomorphism $\End(\PB)\rightarrow \End(\PB\otimes \PB)$, hence it follows that $f\circ g= g\circ f$ as claimed.

Proposition \ref{prop:resProp} says that precomposing with $\eta$ is an isomorphism $\End(\PB)\rightarrow \Hom(\one,\PB)$.  To see that the induced algebra structure is as claimed, observe that
\begin{eqnarray*}
(f\otimes g)\circ (\eta\otimes \eta) &= & (\Id_\PB\otimes (f\circ g))\circ (\eta\otimes \eta)\\
& =& (\eta\otimes \Id_\PB)\circ (\Id_\one\otimes (f\circ g\circ \eta))\\
\end{eqnarray*}
and also
\begin{eqnarray*}
(f\otimes g)\circ (\eta\otimes \eta) &\simeq & (f\circ \eta)\otimes (g\circ \eta)\\
\end{eqnarray*}
Post-composing with the equivalence $(\eta\otimes \Id_\PB)\inv : \PB\otimes \PB\rightarrow \one\otimes \PB$, we see that
\[
(\eta\otimes \Id_\PB)\inv \circ  \Big((f\circ \eta)\otimes (g\circ \eta)\Big) \ \ \simeq \ \ \Id_\one\otimes (f\circ g\circ \eta)
\]
which implies that the induced algebra structure on $\Hom(\one,\PB)$ is as in the statement.
\end{proof}

\begin{remark}
In the situation of interest to us, we will take $\AS=\DS_n$ and $\PB=P_n$.  We are interested in computing not just $\End_{\DS_n}(P_n)$, but the entire triply graded ring $\DEnd(P_n)$.

The same considerations apply with one important modification.  The shifts $(0,j)\ip{0}$ and $(0,0)\ip{k}$ are both of ``homological type,'' which means that there are signs involved in the isomorphisms
\[
A(0,j)\ip{0}\otimes B\cong A\otimes B(0,j)\ip{0} \ \ \ \ \ \ \ \ \text{and} \ \ \ \ \ \ \ \ \ A(0,0)\ip{k}\otimes B\cong A\otimes B(0,0)\ip{k}.
\]
Thus, the overall sign in commuting $f\otimes \Id$ past $\Id\otimes g$ is $(-1)^{|f||g|}$, where $|\cdot |$ is the total homological degree, as described below.  Consequently, $\DEnd(P_n)$ is \emph{graded commutative}. See \S 3.1 in \emph{loc.~cit.} for more on graded monoidal categories.
\end{remark}

\begin{proposition}\label{prop-HHHcommutative}
$\HHH(P_n)\cong \DHom(\one_n,P_n)\cong \DEnd(P_n)$ is a graded commutative ring with unit $[\eta_n]$, where $\eta_n:\one_n\rightarrow P_n$ is the unit map.  Here, graded commutativity means $[f] [g]=(-1)^{|f||g|}[g] [f]$ where $|f|=j+k$ is the total homological degree of an element $f\in \HHH^{ijk}(P_n)$.
\end{proposition}
\begin{proof}
Since the inclusion $\K^-(\SBim_n)\rightarrow \DS_n$ is a monoidal functor,  $(P_n,\eta_n)$  is a unital idempotent in $\DS_n$.  Thus $\Hom_{\DS_n}(P_n(i,j)\ip{k},P_n)\cong\Hom_{\DS_n}(\one_n(i,j)\ip{k},P_n)$, by Proposition \ref{prop:resProp}.  Thus, $\HHH(P_n)\cong \DHom(\one_n,P_n)\cong \DEnd(P_n)$ as triply graded abelian groups.  Graded commutativity of $\DEnd(P_n)$ also follows by the general arguments above.
\end{proof}

We will want to understand $\HHH(P_n)$ not just as a ring, but as an algebra over $\HHH(P_{n-1})$.  Again, there is much that follows from general arguments.  Let $\eta_{\IB}:\one\rightarrow \IB$ be a unital idempotent, and let $\nu:\IB\rightarrow \PB$ be a morphism such that $\nu\circ \eta_{\IB}=\eta_{\PB}$.  Recall Theorem \ref{thm:fundTheoremOfIdemp}, which implies that $\nu$ becomes an isomorphism after tensoring with $\PB$ on the left or right.  In particular $\IB\otimes \PB\simeq \PB\simeq \PB\otimes \IB$.  Then $\PB$ has the structure of a unital idempotent relative to $\IB$, hence Proposition \ref{prop:resProp} implies that precomposing with $\nu$ gives an isomorphism $\End(\PB)\rightarrow \Hom(\IB,\PB)$.

\begin{example}
We could take $\PB=P_n$ and $\IB=P_{n-1}\sqcup \one_1$.
\end{example}

\begin{definition}\label{def:rho}
Retain notation as above.  Let $\rho:\End_{\AS}(\IB)\rightarrow \End_{\AS}(\PB)$ denote the unique map satisfying
\begin{equation}\label{eq:rhoDef}
\rho(f)\circ \nu = \nu\circ f.
\end{equation}
for all $f\in \End_{\AS}(\IB)$.
\end{definition}

\begin{proposition}\label{prop:rhoProp}
The map $\rho:\End_{\AS}(\IB)\rightarrow \End_{\AS}(\PB)$ is an algebra map.  Furthermore,
\[
f\otimes \Id_{\PB}  = \Id_{\IB}\otimes \rho(f) \in \End_{\AS}(\IB\otimes \PB) \ \ \ \ \ \ \ \ \text{ and }  \ \ \ \ \  \ \ \ \Id_{\PB}\otimes f  = \rho(f)\otimes \Id_{\IB} \in \End_{\AS}(\PB\otimes \IB)
\]
for all $f\in \End_{\AS}(\IB)$.
\end{proposition}
\begin{proof}
Clearly, $\Id_\PB\circ \nu = \nu\circ \Id_\IB$, hence $\rho(\Id_\IB)=\Id_\PB$.  Moreover, $\rho(f)\circ \rho(g) \circ \nu = \nu\circ f\circ g$, hence $\rho(f\circ g)=\rho(f)\circ \rho(g)$ for all $f,g\in \End(\IB)$, by definition of $\rho$.  This proves the first statement.

Now let $f\in \End(\IB)$ be arbitrary.  Proposition \ref{prop:resProp} implies that precomposing with $\Id_\IB\otimes \nu$ is an isomorphism $\End(\IB\otimes \PB)\rightarrow \Hom(\IB\otimes \IB,\IB\otimes \PB)$.  To show that $f\otimes \Id_{\PB}  - \Id_{\IB}\otimes \rho(f) $ is zero, it suffices to show that it becomes zero upon precomposing with $\Id_{\IB}\otimes \nu$.  Compute:
\begin{eqnarray*}
(f\otimes \Id_{\PB}  - \Id_{\IB}\otimes \rho(f)) \circ (\Id_{\IB}\otimes \nu) &=& f\otimes \nu -  \Id_\IB\otimes(\rho(f)\circ \nu)\\
&=&  f\otimes \nu -  \Id_\IB\otimes (\nu\circ f)\\
&=& (\Id_\IB\otimes \nu)\circ (f\otimes \Id_\IB - \Id_\IB\otimes f),
\end{eqnarray*}
which is zero by Proposition \ref{prop:mapHop}. 
\end{proof}

Our goal is to compute $\HHH(P_n)$, since this is the stable homology of torus links.  The above simplifies our task by establishing that $\HHH(P_n)$ is a graded commutative $\HHH(P_{n-1})$-algebra in a canonical way.  In particular, letting $\HHH(P_{n-1})$ act on $1\in \HHH(P_n)$ gives a ring morphism $\HHH(P_{n-1})\rightarrow \HHH(P_n)$.  We will utilize this structure in the computations below.

\subsection{The endomorphism ring of $P_n$}
\label{subsec-endP}
We now compute $\DEnd(P_{n})$ using the Markov move and the description of $P_{n}$ from Construction \ref{construction-Pn}.  In this section $P_{n-1}$ denotes the idempotent in $\KC^-(\SBim_n)$, and we write $I(P_{n-1})=P_{n-1}\sqcup \one_1$ for the corresponding idempotent in $\KC^-(\SBim_n)$.  Proposition \ref{prop-PnRelPk} says that $P_n$ is a unital idempotent relative to $I(P_{n-1})$, and so (a graded version of) Proposition \ref{prop:rhoProp} gives us an algebra map $\rho:\DEnd(I(P_{n-1}))\rightarrow \DEnd(P_n)$.  On the other hand,
\[
\DEnd(I(P_{n-1})\cong \DEnd(P_{n-1})\otimes_\Z \Z[x_n].
\]
We will restrict to $\DEnd(P_{n-1})\otimes 1$, and we write the resulting algebra map $\DEnd(P_{n-1})\rightarrow \DEnd(P_n)$ also by $\rho$, by abuse.

The goal of this section is to prove the following.

\begin{proposition}\label{prop:conditionalEndRings}
Assume that Theorem \ref{thm-eigenconesAreOrtho} holds for $P_1,\ldots,P_{n-1}$.  Then Theorem \ref{thm-endP} holds for $P_{n}$.
\end{proposition}

We will prove this by induction on $n\geq 1$.  In the base case, $P_1=\ring[x_1]$ is the trivial bimodule, and $\DEnd(P_1)=\HHH(P_1)$ is the ordinary Hochschild cohomology of $R=\ring[x_1]$, which is $\ring[x_1]\oplus \ring[x_1](-2)[1]$ by Proposition \ref{prop:markov}.  This is clearly isomorphic to $\ring[u_1,\xi_1]$ as claimed.  Now, let $n\geq 2$ be given and assume by induction that:

\begin{hypothesis}\label{hyp-indHyp}
$\DEnd(P_{k})\cong \ring[u_1,\ldots,u_{k},\xi_1,\ldots,\xi_{k}]$ as triply graded algebras, for $1\leq k\leq n-1$.  Assume also that there exist maps $\b_{(n-1,1)},\b_{(n)}\in\Homg(P_{n-1}Y_n,P_{n-1}X_n)$ as in Theorem \ref{thm-eigenconesAreOrtho}.
\end{hypothesis}

Before continuing, we set up some notation.
\begin{definition}\label{def:rhoAndU}
Let $[U_k],[\Xi_k]\in \DEnd(P_{n-1})$ denote the classes corresponding to $u_k,\xi_k$ under the isomorphism from Hypothesis \ref{hyp-indHyp}.  Denote the image of $[U_k],[\Xi_k]\in \DEnd(P_{n-1})$ under the map $\rho$ also by $[U_k],[\Xi_k]\in \DEnd(P_n)$ ($1\leq k\leq n-1$).  Let $[U_n]\in \DEnd(P_n)$ denote the periodicity endomorphism of $P_n$ (see also Proposition \ref{prop:periodicityMap}).
\end{definition}

Proposition \ref{prop-HHHcommutative}, which says that $\DEnd(P_n)\cong \DHom(\one_n,P_n)$ (the isomorphism is given by precomposing with the unit map $\eta_n:\one_n\rightarrow P_n$).  This is a major simplification indeed, since to compute $\DEnd(P_n)$ from the definitions would require computing the homology of a bi-infinite complex with uncountably generated chain groups, while $\DHom(\one_n,P_n)$ is the homology of a bounded above chain complex, each of whose chain groups is finitely generated as an $R_n$-module. 

The adjunction isomorphism (Corollary \ref{cor:adjunctionComplexes}) then applies, yielding
\begin{equation}\label{eq:DendTrace}
\DEnd(P_n)\cong \DHom(\one_{n-1},\Tr(P_n)).
\end{equation}

Our next task is to compute $\Tr(P_{n})$, assuming Hypothesis \ref{hyp-indHyp}.  Our computation is aided by the fact that $\DEnd(P_{n-1})$ is supported in even homological degrees ($t$-degrees).

By induction, we assume that Theorem \ref{thm-eigenconesAreOrtho} holds for $P_n$.  Thus, we can construct $P_n$ as in Construction \ref{construction-Pn}.  More precisely,
\[
P_n = \Cone\Big(\Z[u_n]\otimes P_{n-1}Y_{n-1}(2n,0)\ip{2-2n} \ \ \longrightarrow \ \ \Z[u_n]\otimes P_{n-1}X_n \Big)
\]
in which the map is $1\otimes \b_{(n-1,1)} - u_n\otimes \b_{(n)}$. 
\begin{definition}\label{def:periodicity}
Let $U_n\in \DEnd(P_n)$ denote the endomorphism induced by multiplication by $u_n$.
\end{definition}

We may apply the linear functor $\Tr$ to the above description of $P_n$, obtaining
\[
\Tr(P_n) = \Cone\Big(\Z[u_n]\otimes\Tr(P_{n-1}Y_{n-1})(2n,0)\ip{2-2n} \ \ \longrightarrow \ \\Z[u_n]\otimes\Tr(P_{n-1}X_n) \Big).
\]
Now, observe that
\begin{equation}\label{eq:TracePieces}
\Tr(P_{n-1}Y) = \sspic{.5in}{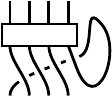} \simeq \sspic{.35in}{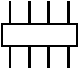}(-4,1)\ip{1}  \hskip.4in \Tr(P_{n-1}X) = \sspic{.5in}{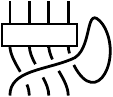} \simeq \sspic{.35in}{P}(0,0)\ip{0}.
\end{equation}
We have used the Markov move (Proposition \ref{prop:markov}) to undo the curls, and then used the fact that $P_{n-1}$ absorbs Rouquier complexes.  Choose homotopy equivalences and denote them by $\psi_Y:\Tr(P_{n-1}Y)\rightarrow P_{n-1}(-4,1)\ip{1}$ and $\psi_X:\Tr(P_{n-1}X)\rightarrow P_{n-1}$.  Tensoring with $\Z[u_n]$ gives
\[
1\otimes \psi_Y  \ :\   \  \ \Z[u_n]\otimes_\Z\Tr(P_{n-1}Y) \buildrel \simeq \over \longrightarrow   \Z[u_n]\otimes_\Z P_{n-1}(-4,1)\ip{1} 
\]
and
\[
1\otimes \psi_X  \ : \ \Z[u_n] \otimes_\Z\Tr(P_{n-1}X) \longrightarrow  \Z[u_n] \otimes_\Z P_{n-1}.
\]
The equivalences and the assoicated homotopies all commute with multiplication by $u_n$.  Thus, our description of $\Tr(P_n)$ becomes
\[
\Tr(P_n) \ \ \simeq \ \ \Cone\Big(\Z[u_n]\otimes P_{n-1}(2n-4,1)\ip{3-2n} \ \ \longrightarrow \ \ \Z[u_n]\otimes P_{n-1}\Big),
\]
in which the differential is
\[
(1\otimes \psi_X)\circ \Big(1\otimes \Tr(\b_{(n-1,1)}) - u_n\otimes \Tr(\b_{(n)})\Big)\circ (1\otimes \psi_Y\inv) \ = \  1\otimes f_1 - u_n\otimes f_0
\]
where $f_0=\psi_X\circ \Tr(\b_{(n-1,1)}\circ \psi_Y\inv$ and $f_1=\psi_X\circ \Tr(\b_{(n)})\circ \psi_Y\inv$ are closed elements of $\DEnd(P_{n-1})$  of degree $\deg(f_0)=(2n-4,1,3-2n)$ and $\deg(f_1)=(-4,1,1)$.  By Hypothesis \ref{hyp-indHyp}, $\DEnd(P_{n-1})$ is supported in even $t$-degrees, so we conclude that $f_0\simeq f_1\simeq 0$.  Thus, $\Tr(P_{n})$ splits as a direct sum of infinitely many copies of $P_{n-1}$. We can keep track of degrees in a convenient way:
\begin{equation}\label{eq-trSimp}
\Tr(P_{n})\simeq \Z[u_n,\xi_n]\otimes P_{n-1}
\end{equation}

We collect our findings below.

\begin{proposition}\label{prop:TrP}
Assuming Hypothesis \ref{hyp-indHyp}, there is an equivalence of the form (\ref{eq-trSimp}) in $\DS_n$, where $u_{n}$ is an even formal indeterminate of degree $q^{2n}t^{2-2n}$ and $\xi_n$ is an odd formal indeterminate of degree $q^{2n-4}t^{2-2n}a$.  This equivalence and the various homotopies can all be chosen to be $\Z[u_n]$-equivariant. \qed
\end{proposition}
%By abuse, we will also regard $u_n$ as the endomorphism of $P_n$ coming from the periodicity of the construction, which includes $P_{n}$ as a proper subcomplex of itself, with appropriate shifts.

\begin{corollary}\label{cor:endCor}
Under Hypothesis \ref{hyp-indHyp}, we have an isomorphism
\begin{equation}\label{eq:THEcomputation}
\Phi:\DHom(\one_{n-1},\Tr(P_n)) \buildrel \cong  \over \rightarrow \Z[u_n,\xi_n]\otimes \DHom(\one_{n-1},P_{n-1})
\end{equation}
as triply graded $\Z[u_n]$-modules.  Consequently, $\DEnd(P_n)\cong \Z[u_1,\ldots,u_n,\xi_1,\ldots,\xi_n]$ as triply graded $\Z[u_n]$-modules.
\end{corollary}
\begin{proof}
This is obtained by applying $\DHom(\one_{n-1},-)$ to (\ref{eq-trSimp}). Since $\one_{n-1}$ is supported in a single homological degree, any homogeneous morphism $\one_{n-1}(a,b)\ip{c}\rightarrow \Z[u_n,\xi_n]\otimes P_{n-1}$ hits only finitely terms $u_n^k\xi_n^\ell\otimes P_{n-1}$, because $u_n$ has negative homological degree for $n>1$.  Thus,
\[
\Hom_{\DS_{n-1}}(\one_{n-1}(a,b)\ip{c}, \Z[u_n,\xi_n]\otimes P_{n-1}) \cong \Hom_{\DS_{n-1}}(\one_{n-1}(a,b)\ip{c}, \Z[u_n,\xi_n]\otimes P_{n-1}).
\]

The last statement follows from $\DHom(\one_{n-1},P_{n-1})\cong \DEnd(P_{n-1})$ together with the induction Hypothesis \ref{hyp-indHyp}. 
\end{proof}

We must promote the statement of Corollary \ref{cor:endCor} to an isomorphism of rings.  First, we bring the $\DEnd(P_{n-1})$ action into the picture. 

\begin{lemma}\label{lemma:equivariance0}
For all $g\in \DEnd(P_{n-1})$ we have
\[
g\otimes \Id_{\Tr(P_n)} \simeq \Id_{P_{n-1}}\otimes \Tr(\rho(g)) \ \in \ \DEnd(P_{n-1}\otimes \Tr(P_n)).
\] 
\end{lemma}
\begin{proof}
Let $g\in \DEnd(P_{n-1})$ be arbitrary.  Consider the following  diagram.
\[
\begin{diagram}
P_{n-1}\Tr(P_n) & \rTo^{\cong} & \Tr(I(P_{n-1})P_n) & \rTo^{\cong} & \Tr(I(P_{n-1})P_n) & \rTo^{\cong} & P_{n-1}\Tr(P_n)\\
\dTo^{g  \Tr(P_n)}&& \dTo^{\Tr(I(g)P_n)} && \dTo^{\Tr(I(P_{n-1})\rho(g))}&& \dTo^{P_{n-1}\Tr(\rho(g))}\\
P_{n-1}\Tr(P_n) & \rTo^{\cong}  & \Tr(I(P_{n-1})P_n) & \rTo^{\cong} &  \Tr(I(P_{n-1})P_n) & \rTo^{\cong} & P_{n-1}\Tr(P_n).
\end{diagram}
\]
The left and right squares commute by naturality of the isomorphism in Lemma \ref{lemma-traceIsLinear}.  The middle square commutes up to homotopy by Proposition \ref{prop:rhoProp}.  Commutativity of the diagram gives the relation in the statement.
\end{proof}

\begin{lemma}\label{lemma:equivariance1}
The isomorphism (\ref{eq:THEcomputation}) commutes with the $\DEnd(P_{n-1})$ action.
\end{lemma}
\begin{proof}
Note that $\DEnd(P_{n-1})$ acts on $P_n$, via the map $\rho:\DEnd(P_{n-1})\rightarrow \DEnd(P_n)$.  Thus, $\DEnd(P_{n-1})$ acts on each of $\DEnd(P_n)$, $\DHom(\one_n,P_n)$, and $\DHom(\one_{n-1},\Tr(P_n))$, via its action on the second arguments (applying the functor $\Tr$ in the case of the latter hom space).  Pre-composition with $\eta_n:\one_n\rightarrow P_n$ commutes with post-composition with $\rho(g)$, hence the isomorphism
\[
\eta_n^\ast:\DEnd(P_n)\cong \DHom(\one_n,P_n)
\]
commutes with the $\DEnd(P_{n-1})$ actions.  Naturality implies that the adjunction isomorphism
\[
\DHom(\one_n,P_n)\cong \DHom(\one_{n-1},\Tr(P_n))
\]
commutes with the $\DEnd(P_{n-1})$ actions. 

Now, consider the sequence of homotopy equivalences (and isomorphisms)
\[
\Tr(P_n)  \ \ \buildrel(1)\over \simeq \ \ P_{n-1}\Tr(P_n) \ \ \buildrel(2)\over \simeq \ \ P_{n-1}\Big(\Z[u_n,\xi_n]\otimes_\Z P_{n-1}\Big) \ \ \buildrel(3)\over\cong \ \ \Z[u_n,\xi_n]\otimes_\Z (P_{n-1}P_{n-1}) %\ \ \buildrel(4)\over\simeq \ \ \Z[u_n,\xi_n]\otimes_\Z P_{n-1}
\]
Lemma \ref{lemma:equivariance0} implies that the the equivalence (1) commutes with the $\DEnd(P_{n-1})$ actions (up to homotopy).  The equivalence (2) commutes with the $\DEnd(P_{n-1})$ action on the first factor $P_{n-1}$.  The isomorphism (3) is merely an application of the distributivity relation between $\otimes$ and $\oplus$; it commutes with the $\DEnd(P_{n-1})$ action.  Finally, applying the equivalence $P_{n-1}P_{n-1}\simeq P_{n-1}$ (which commutes with the $\DEnd(P_{n-1})$-action up to homotopy by the proof of Proposition \ref{prop:mapHop}), we see that the equivalence (\ref{eq-trSimp}) commutes with the $\DEnd(P_{n-1})$ actions up to homotopy.

Applying the functor $\DHom(\one_{n-1},-)$, we find that the equivalence $\DHom(\one_{n-1},\Tr(P_n))\cong \Z[u_n,\xi_n]\otimes_\Z \DHom(\one_{n-1},P_{n-1})$ commutes with the $\DEnd(P_{n-1})$ action.  Putting all of this together proves the lemma.
\end{proof}

We are ready to prove Proposition \ref{prop:conditionalEndRings}.

\begin{proof}[Proof of Proposition \ref{prop:conditionalEndRings}]
To complete the proof we must show that the isomorphism $\Phi$ from (\ref{eq:THEcomputation}) is an algebra map.  Let us recap how $\Phi$ is defined.   Let $f\in \DEnd(P_n)$ be given.  Let $\eta_n:\one_n\rightarrow P_n$ denote the unit map.  First, form $f\circ \eta_n\in \DHom(\one_n,P_n)$.  Now apply the adjunction isomorphism $A:\DHom(\one_n,P_n)\rightarrow \DHom(\one_{n-1},\Tr(P_n))$, obtaining
\[
A(f\circ \eta) = \Tr(f)\circ A(\eta) \ \in \DHom(\one_{n-1},\Tr(P_n)),
\]
by naturality of $A$.  Then post-compose with the simplification $\phi:\Tr(P_n)\simeq \Z[u_n]\otimes P_{n-1}$.  The result is
\begin{equation}\label{eq:PhiDef}
\Phi(f):=\phi\circ \Tr(f)\circ A(\eta_n).
\end{equation}
This is a morphism from $\one_{n-1}$ to $P_{n-1}\otimes \Z[u_n,\xi_n]$ in $\DS_{n-1}$.  Such morphisms can be thought of as elements of $\Z[u_n,\xi_n]\otimes \DHom(\one_{n-1},P_{n-1})$.  This defines $\Phi:\DEnd(P_n)\rightarrow \Z[u_n,\xi_n]\otimes \DHom(\one_{n-1},P_{n-1})$.

\textbf{Claim 1.}  $\Phi(\Id_{P_n}) = 1\otimes [\eta_{n-1}]$.  This is a straightforward verification.  In any case this is forced for degree reasons. The degree $(0,0,0)$ component of $\DEnd(P_n)$ is isomorphic to $\Z$ as a ring, with generator $[\Id_{P_n}]$.  This is an easy consequence of Corollary \ref{cor:endCor}, together with the fact that there are only two unital ring structures on $\Z$ with its standard additive structure (they are related by the automorphism $x\mapsto -x$).  Thus, $\Phi([\Id_{P_n}])=\pm 1$ since $\Phi$ is a group isomorphism, hence sends generators to generators.  Up to replacing $\Phi$ by $-\Phi$, we may arrange that  $\Phi([\Id_{P_n}])=+ 1$.

\textbf{Claim 2.}  $\Phi$ is $\Z[u_n]$ equivariant.  Compute:
\begin{eqnarray*}
\Phi(U_n\circ f) & =& \phi\circ \Tr(U_n)\circ \Tr(f) \circ A(\eta_n) \\
& = & u_n\circ \phi \circ A(f\circ \eta_n).
\end{eqnarray*}
The second equality holds by naturality of the adjunction isomorphism $A$, and the third holds since $\phi$ is equivariant with respect to the periodicity endomorphism (Proposition \ref{prop:TrP}).

Thus, $\Phi:\DEnd(P_n)\rightarrow \Z[u_n,\xi_n]\otimes \DHom(\one_{n-1},P_{n-1})$ sends $1\mapsto 1$ and is equivariant with respect to the action of $\Z[u_n]\otimes \DEnd(P_{n-1})$ (equivariance with respect to the $\DEnd(P_{n-1})$ action was proven in Lemma \ref{lemma:equivariance1}).  On the other hand, the domain and codomain of $\Phi$ are free of rank 2 over $\Z[u_n]\otimes \DEnd(P_{n-1})$, so we are very nearly done.  To finish the proof, we will construct an algebra map
\[
\Psi:\Z[u_n,\xi_n]\otimes \DHom(\one_{n-1},P_{n-1})\rightarrow \DEnd(P_n)
\]
and show that $\Phi\circ \Psi = \Id$.  First, let $\Xi_n\in \DEnd(P_{n})$ denote $\Phi\inv(\xi_n)$.  Proposition \ref{prop-HHHcommutative} says that $g\mapsto g\circ \eta_{n-1}$ is an algebra isomorphism
\[
\DEnd(P_{n-1}) \rightarrow  \DHom(\one_{n-1},P_{n-1}).
\]
Thus, we have an algebra map $\DHom(\one_{n-1},P_{n-1})\rightarrow \DEnd(P_n)$ uniquely characterized by $g\circ \eta_{n-1}\mapsto \rho(g)$.  We also have an algebra map $\Z[u_n,\xi_n]\mapsto \DEnd(P_n)$ characterized by $u_n\mapsto [U_n]$ and $\xi_n\mapsto [\Xi_n]$.  Tensoring these together defines $\Psi$. 

Note that $\Phi\circ \Psi$ acts as the identity on $1\otimes [\eta_{n-1}]$ and $\xi_n\otimes [\eta_{n-1}]$, by construction.  Moreover, $\Phi\circ \Psi$ is equivariant with respect to the $\Z[u_n]\otimes \DEnd(P_{n-1})$ action, and the endomorphism $\Phi\circ \Psi - \Id$ annihilates a set of generators.  It follows that $\Phi\circ \Psi=\Id$, which proves that $\Phi=\Psi\inv$ is an algebra isomorphism.

Thus, $\DEnd(P_n)\cong \Z[u_n,\xi_n]\otimes \DEnd(P_{n-1})\cong \Z[u_1,\ldots,u_n,\xi_1,\ldots,\xi_n]$ as algebras.
\end{proof}

Before moving on, we collect some results regarding the action of $\Z[u_1,\ldots,u_n]$ on $P_n$.

\begin{proposition}\label{prop:periodicityMap}
Assume Hypothesis \ref{hyp-indHyp}, and let $[U_i^{(k)}]\in \DEnd(P_k)$ denote the class (and a chosen representative) corresponding to $u_i$ under the result of Theorem \ref{thm-endP} for $1\leq i\leq k$.  Let $\rho_{nk}:\DEnd(P_k)\rightarrow \DEnd(P_n)$ be the canonical ring map (Proposition \ref{prop-HHHcommutative}).  Then $\rho_{nk}(U_i^{(k)})\simeq U_i^{(n)}$ for all $1\leq i\leq k\leq n$.  The map $U_n^{(n)}$ is homotopic to the map obtained by tensoring $\b_{(n,n-1)}\in \Homg(P_{n-1}Y_n,P_{n-1}X_n)$ on the left with $P_n$ and simplifying. 
\end{proposition}
\begin{proof}
The relations $\rho_{nk}(U_i^{(k)})\simeq U_i^{(n)}$ hold by definition (Definition \ref{def:rhoAndU}).  They are consistent (up to homotopy) since $\rho_{nk}\circ \rho_{ki}\simeq \rho_{ni}$.  

Now, to prove the relation between $U_n^{(n)}$ and $\b_{(n-1,1)}$ we first introduce some abbreviations.  Let $U_n$ denote $U_n^{(n)}$ which by definition is the periodicity map of $P_n$ (Definition \ref{def:periodicity}).  Let $\sigma:\Homg(P_{n-1}Y_n,P_{n-1}X_n)\rightarrow \Endg(P_n)$ denote the map induced by tensoring on the left with $P_n$ and simplifying.  Set $U_n':=\sigma(\b_{(n-1,1)})$.  Note that $\sigma(\b_{(n)})=\Id_{P_n}$ (Remark \ref{rmk:bn}).  We wish to prove that $U_n'\simeq U_n$.

Proposition \ref{prop:mapHop} says that any closed endomorphism $g\in \Endg(P_n)$ satisfies
\[
g P_n\simeq P_n g \in \Endg(P_nP_n).
\]
We will prove that $U_n'P_n \simeq P_nU_n$, from which the proposition will follow.  From Construction \ref{construction-Pn} we have
\[
P_n = \Cone\left(\begin{diagram} P_{n-1}Y_n[u_n] & \rTo^{\b_{(n-1,1)}\otimes 1 - \b_{(n)}\otimes u_n} & P_{n-1}X_n[u_n]\end{diagram}\right),
\]
omitting the grading shifts.  Tensoring with $P_n$ yields
\[
P_nP_n = \Cone\left(\begin{diagram} P_nP_{n-1}Y_n[u_n] & \rTo^{P_n\b_{(n-1,1)}\otimes 1 - P_n\b_{(n)}\otimes u_n} & P_nP_{n-1}X_n[u_n]\end{diagram}\right).
\]
Applying the homotopy equivalence $P_nP_{n-1}Y_n\simeq P_n$ and $P_nP_{n-1}X_n\simeq P_n$ yields
\[
P_nP_n\simeq \Cone\left(\begin{diagram} P_n[u_n] & \rTo^{\sigma(\b_{(n-1,1)}) \otimes 1 - \sigma(\b_{(n)})\otimes u_n} & P_n[u_n]\end{diagram}\right).
\]
Since $\sigma(\b_{(n-1,1)})\simeq U_n'$ and $\sigma(\b_{(n)})\simeq \Id_{P_n}$, we obtain
\[
P_nP_n\simeq \Cone\left(\begin{diagram} P_n[u_n] & \rTo^{U_n' \otimes 1 - \Id_{P_n}\otimes u_n} & P_n[u_n]\end{diagram}\right).
\]
All of the above equivalences (and the homotopies) commute with multiplication by $u_n$ (i.e.~$U_n$). The proposition now follows from general arguments, as we now explain.  Let $C\in \KC(\CC)$ be any complex over an additive category, let $\l$ be any grading shift functor, and let $f:\l C\rightarrow C$ be a chain map.  Then in terms of diagrams, we have
\[
\Cone\left(\begin{diagram} \l C[u] & \rTo^{f \otimes 1 - \Id_{C}\otimes u} & P_n[u]\end{diagram}\right) \ \ \ = \ \ \ 
\left(
\begin{minipage}{1.7in}
\begin{tikzpicture}[baseline=-0.25em]
\matrix (m) [ssmtf,column sep=5em, row sep=5em]
{\l C\ip{1} & C \\
\l^2 C\ip{1}  &  \l C \\
\cdots & \cdots\\};
\path[->,>=stealth',shorten >=1pt,auto,node distance=1.8cm,font=\small]
(m-1-1) edge node[auto] {$f$} (m-1-2)
(m-1-1) edge node[auto] {$-\Id_C$} (m-2-2)
(m-2-1) edge node[auto] {$\l f$} (m-2-2)
(m-2-1) edge node[auto] {$-\Id_{\l C}$} (m-3-2);
\end{tikzpicture}
\end{minipage}
\right).
\]
After Gaussian eliminations, this complex is homotopy equivalent to $C$.  It is a straightforward exercise to check that multiplication by $u$ (which acts by shifting the above diagram down) corresponds to $f$ under this equivalence.  Applying these considerations to the above description of $P_nP_n$ shows that $U_n'\simeq U_n$, as claimed.
\end{proof}
This justifies the sign on $\b_{(n)}$ in our construction of $P_n$.

\subsection{The map $\b_{(n-1,1)}$}
\label{subsec:hookEigenmap}
Throughout this section, we will abuse notation in the following way.  If $C\in \KC(\SBim_k)$, then we will denote $C\sqcup \one_{n-k}\in \KC(\SBim_n)$ also by $C$.

To complete our inductive proof of Theorems \ref{thm-eigenconesAreOrtho} and \ref{thm-endP}, we construct $\b_{(n-1,1)}$, assuming the following.
\begin{hypothesis}\label{hyp:inductEigenmap}
Assume Theorem \ref{thm-endP} and Theorem \ref{thm-eigenconesAreOrtho} hold for $P_1,\ldots,P_{n-1}$.
\end{hypothesis}

In particular $\DEnd(P_{n-1})$ is supported in \emph{even} homological degrees. We will use this parity observation to construct $\b_{(n-1,1)}$.

In this section it will be useful to draw our chain complexes diagrammatically.  For instance, if $A,B,C,D$ are chain complexes, then
\begin{equation}\label{eq:ABCD}
\left(\begin{minipage}{2.59in}
\begin{tikzpicture}
\node (a) at (0,0){$A$ };
\node (b) at (2,0){ $B$};
\node (c) at (4,0){ $C$ };
\node (d) at (6,0){ $D$ };
\path[->,>=stealth',shorten >=1pt,auto,node distance=1.8cm,font=\small]
(a) edge node[above,xshift=4pt] {$d_{BA}$} (b)
(b) edge node {$d_{CB}$} (c)
(c) edge node {$d_{DC}$} (d)
(a) edge[bend left=30] node {$d_{DA}$} (d);
\draw[frontline,->,>=stealth',shorten >=1pt,auto,node distance=1.8cm]
(a) to [bend right] node[below] {$d_{CA}$} (c);
\draw[frontline,->,>=stealth',shorten >=1pt,auto,node distance=1.8cm]
(b) to [bend right] node[below] {$d_{DB}$} (d);
\end{tikzpicture}
\end{minipage}\right)
\end{equation}
will denote the chain complex
\[
A\oplus B\oplus C\oplus D \text{ with differential } \begin{bmatrix}d_A &0&0&0\\ d_{BA} & d_B &0&0\\ d_{CA} & d_{CA} & d_C &0\\ d_{DA} & d_{DB} & d_{DC} & d_D  \end{bmatrix}.
\]
Here,  $d_{YX}$ are degree 1 elements of $\Homg(X,Y)$, for $X,Y\in \{A,B,C,D\}$.  The equation $d^2=0$ holds if and only if
\begin{enumerate}\setlength{\itemsep}{2pt}
\item $d_A^2=d_B^2=d_C^2=d_D^2=0$.
\item $d_{CB}\circ d_{BA} = -(d_C \circ d_{CA} + d_{CA}\circ d_A)$, so that $d_{CB}\circ d_{BA}$ is null-homotopic with homotopy $-d_{CA}$.
\item Similarly $d_{DC}\circ d_{CB} = -(d_D\circ d_{DB} + d_{DB}\circ d_B)$.
\item $d_{DB}\circ d_{BA} + d_{DC}\circ d_{CA} = - (d_D\circ d_{DA} + d_{DA}\circ d_A)$.  The left hand side is a degree 2 cycle in $\Homg(A,D)$, and this equation is interpreted as saying that this cycle is null-homotopic with homotopy $-d_{DA}$.
\end{enumerate} 

Note also that the chain complex (\ref{eq:ABCD}) can be reassociated into a chain complex
\[
\Big((A\rightarrow B)\rightarrow (C\rightarrow D) \bigg).
\]
which is the mapping cone on a map $(A\rightarrow B)\ip{1}\rightarrow (C\rightarrow D)$.  So if we want to construct chain maps between mapping cones, we may as well construct chain complexes of the form (\ref{eq:ABCD}).  This will be our strategy.  We begin with some observations.  First, note that the Rouquier complex associated to a crossing can be written in this language as
\[
F(\sigma_{n-1}) = \Big(B_{n-1}(1)\ip{-1} \rightarrow R\Big)
\]
where the map is the ``dot'' map.  The shift $\ip{-1}$ on $B_{n-1}$ indicates that this term is in homological degree -1.  Tensoring on the right with $F(\sigma_1\cdots \sigma_{n-2})$ gives
\[
X_n \simeq \Big(B_{n-1}X_{n-1}(1)\ip{-1} \rightarrow X_{n-1}\Big).
\]
Now, tensoring with $P_{n-1}$ gives
\[
P_{n-1}X_n \simeq \Big(P_{n-1}B_{n-1}X_{n-1}(1)\ip{-1} \buildrel d_{DC}\over \longrightarrow P_{n-1}\Big)
\]
since $P_{n-1}$ absorbs Rouquier complexes.  We are denoting the horizontal component of the differential by $d_{DC}$ for reasons that will become clear momentarily.  Note that $d_{DC}$ is simply the ``dot'' map $P_{n-1}B_{n-1}X_{n-1}\rightarrow P_{n-1}X_{n-1}$ followed by a homotopy equivalence.  Similar considerations, starting from $F(\sigma_{n-1}\inv)=(R\rightarrow B_{n-1}(-1)\ip{1})$, yield
\[
P_{n-1}Y_n \simeq \Big(P_{n-1} \buildrel d_{BA}\over \longrightarrow P_{n-1}B_{n-1}Y_{n-1}(-1)\ip{1} \Big).
\]
We wish to construct a chain map $P_{n-1}Y_n(2n)\ip{2-2n} \rightarrow P_{n-1}X_n$.  By expanding each term according to the observations above, and taking the mapping cone, we see that we must construct a chain complex $Q_n$ of the form
\begin{equation}\label{eq:Qdiagram1}
Q_n=\left(\begin{minipage}{5.2in}
\begin{tikzpicture}
\node (a) at (0,0){$P(2n)\ip{1-2n} $ };
\node (b) at (4.5,0){ $PBY(2n-1)\ip{2-2n} $};
\node (c) at (8.5,0){ $PBX(1)\ip{-1}$ };
\node (d) at (11.7,0){ $P$ };
\path[->,>=stealth',shorten >=1pt,auto,node distance=1.8cm,font=\small]
(a) edge node[auto] {$d_{BA}$} (b)
(b) edge node {$d_{CB}$} (c)
(c) edge node {$d_{DC}$} (d)
(a) edge[bend left=23] node {$d_{DA}$} (d);
\draw[frontline,->,>=stealth',shorten >=1pt,auto,node distance=1.8cm]
(a) to [bend right] node[below] {$d_{CA}$} (c);
\draw[frontline,->,>=stealth',shorten >=1pt,auto,node distance=1.8cm]
(b) to [bend right] node[below] {$d_{DB}$} (d);
\end{tikzpicture}
\end{minipage}\right)
\end{equation}
where we have abbreviated $P=P_{n-1}$, $X=X_{n-1}$, $Y=Y_{n-1}$, and $B=B_{n-1}$.  We must then show that $Q_n\in \ker(P_n)$.  We construct $Q_n$ piece by piece.  The components $d_{DC}$ and $d_{BA}$ were constructed already.  The most interesting component is $d_{CB}$, and we construct this next.  The ideas are fairly simple; the biggest challenge will be parsing notation.

Observe that
\[
P_{n-1}B_{n-1}P_{n-2}Y_{n-1} \simeq P_{n-1}P_{n-2}B_{n-1}Y_{n-1} \simeq P_{n-1}B_{n-1}Y_{n-1}
\]
and 
\[
P_{n-1}B_{n-1}P_{n-2}X_{n-1} \simeq P_{n-1}P_{n-2}B_{n-1}X_{n-1} \simeq P_{n-1}B_{n-1}X_{n-1},
\]
since $P_{n-2}$ commutes past $B_{n-1}$ and gets absorbed by $P_{n-1}$.

Let $\b_{(n-1)},\b_{(n-2,1)}\in \Homg(P_{n-2}Y_{n-1},P_{n-2}X_{n-1})$ be the morphisms from Theorem \ref{thm-eigenconesAreOrtho}, so that $\Cone(\b_{(n-1)})\in \ker(P)$ and $\Cone(\b_{(n-2,1)})\in \im(P)$.  These exist by the induction hypotheses.
\begin{remark}\label{rmk:factorThroughR}
The map $\b_{(n-1)}$ is obtained from the composition $Y_{n-1}\rightarrow \one \rightarrow X_{n-1}$ by tensoring on the left with $P_{n-2}$.  Thus, $\b_{(n-1)}$ factors through $P_{n-2}$.  This fact will be used in the proof of Proposition \ref{prop:hookEigenmap}.
\end{remark}
Tensoring
\[
\b_{(n-2,1)}: P_{n-2}Y_{n-1}(2n-2)\ip{4-2n}\rightarrow P_{n-2}X_{n-1}
\]
on the left with $P_{n-1}$ and applying $P_{n-1}P_{n-2}Y_{n-1}\simeq P_{n-1}$ and $P_{n-1}P_{n-2}X_{n-1}\simeq P_{n-1}$ gives a morphism $P_{n-1}(2n-2)\ip{4-2n}\rightarrow P_{n-1}$ which we shall denote by $u_{n-1}$ (compare with Proposition \ref{prop:periodicityMap}). Tensoring
\[
\b_{(n-1)}: P_{n-2}Y_{n-1}\rightarrow P_{n-2}X_{n-1}
\]
on the left with $P_{n-1}$ yields the identity morphism $P_{n-1}\rightarrow P_{n-1}$ up to homotopy, since $\b_{(n-1)}$ restricts to the identity map on the unique $\one$ summands of $P_{n-2}Y_{n-1}$ and $P_{n-2}X_{n-1}$ (see Remark \ref{rmk:bn}).

\begin{lemma}\label{lemma:alphaHop}
We have $u_{n-1}\b_{(n-1)}\simeq P_{n-1}\b_{(n-2,1)}$.
\end{lemma}
\begin{proof}
Recall that for any closed morphism $f\in \Endg(P_{n-1})$ we have $fP_{n-1}\simeq P_{n-1}f$.  Thus, it suffices to show that
\[
u_{n-1}P_{n-1}\b_{(n-1)}- P_{n-1} P_{n-1}\b_{(n-2,1)} \ \in \ \ \Homg(P_{n-1}P_{n-1}P_{n-2}Y_{n-1}, P_{n-1}P_{n-1}P_{n-2}X_{n-1})
\]
is null-homotopic.  Let $\phi$ denote the chosen homotopy equivalence $P_{n-1}\rightarrow P_{n-1}P_{n-2}Y_{n-1}$, and let $\psi$ denote the chosen homotopy equivalence $P_{n-1}P_{n-2}X_{n-1}\rightarrow P_{n-1}$.  For any $f\in \Homg(P_{n-2}Y_{n-1},P_{n-2}X_{n-1})$, set $\sigma(f):= \psi\circ (P_{n-1}f) \circ \phi$.  By construction, $\sigma(\b_{(n-1)})\simeq \Id_{P_{n-1}}$ and $\sigma(\b_{(n-2,1)})\simeq u_{n-1}$. Thus,
\begin{eqnarray*}
(P\psi)\circ (uP\b_{(n-1)}- P P\b_{(n-2,1)}) \circ (P\phi) & \simeq & u\sigma(\b_{(n-1)}) - P\sigma(\b_{(n-2,1)})\\
 & \simeq& u P - P u\\
  & \simeq &0\\
\end{eqnarray*}
where $P=P_{n-1}$ and $u=u_{n-1}$.   This proves that $u_{n-1}P_{n-1}\b_{(n-1)}- P_{n-1} P_{n-1}\b_{(n-2,1)}\simeq 0$, as claimed.
\end{proof}

\begin{definition}\label{def:deltaMaps}
Let $\d_1$ denote the composition
\[
\d_1 \ : \ \begin{diagram} PBY(2n-2)\ip{4-2n} & \rTo^{\simeq} & PBP'Y(2n-2)\ip{4-2n} & \rTo^{u B\b_{(n-1)}} & PBP'X  & \rTo^{\simeq} & PBX\end{diagram}
\]
and let $\d_2$ denote the composition
\[
\d_2 \ : \ \begin{diagram} PBY(2n-2)\ip{4-2n} & \rTo^{\simeq} & PBP'Y(2n-2)\ip{4-2n} & \rTo^{PB\b_{(n-2,1)}} & PBP'X  & \rTo^{\simeq} & PBX\end{diagram},
\]
where $P=P_{n-1}$, $B=B_{n-1}$, $P'=P_{n-2}$, and $u=u_{n-1}$.  Set $d_{CB}:=\d_2-\d_1$.  
\end{definition}

\begin{lemma}\label{lemma:Qexists}
With the definitions of $d_{BA}$, $d_{CB}$, $d_{DC}$ above, there exists a complex $Q_n$ as in diagram (\ref{eq:Qdiagram1}).
\end{lemma}
\begin{proof}
We will show that $d_{CB}\circ d_{BA}$ and $d_{DC}\circ d_{CB}$ are null homotopic, which implies the existence of homotopies $d_{CA}$ and $d_{DB}$.   Assuming this is done, it is trivial to check that the expression
\begin{equation}\label{eq:magicallyNull}
d_{DB}\circ d_{BA} + d_{DC}\circ d_{CA}\in \Hom^2\Big((P_{n-1}\sqcup \one_1)(2n)\ip{1-2n}, P_{n-1}\sqcup \one_1\Big)
\end{equation}
is a cycle.  This cycle corresponds to a degree $(2n)\ip{3-2n}$ cycle of $\Endg(P_{n-1}\sqcup \one_1)$.  But notice that
\[
\Endg(P_{n-1}\sqcup \one_1) \cong \Endg(P_{n-1})\otimes_\ring \ring[x_n],
\]
whose homology is supported in even homological degrees by Proposition \ref{prop:conditionalEndRings}.  In particular there are nonzero elements of homological degree $3-2n$.  This shows that the cycle (\ref{eq:magicallyNull}) is null-homotopic, which then proves the existence of $d_{DA}$.

It remains to prove that $d_{CB}\circ d_{BA}\simeq 0$ and $d_{DC}\circ d_{CB}\simeq 0$.  We prove the second of these; the first is similar.  Recall that $\d_2$ is obtained from
\[
\begin{diagram}
P_{n-1}B_{n-1}P_{n-2}Y_{n-1} & \rTo^{ \Id_{P_{n-1}}\Id_{B_{n-1}}\b_{(n-2,1)}-u_{n-1}\Id_{B_{n-1}}\b_{(n-1)}} & P_{n-1}B_{n-1}P_{n-2}X_{n-1}
\end{diagram}
\]
by pre- and post-composing with some homotopy equivalences.  We are omitting shifts for brevity.  Now, if the $B_{n-1}$ factor were not there, then the resulting map would be null-homotopic by Lemma \ref{lemma:alphaHop}.  Thus, $d_{CB}$ becomes null-homotopic upon post-composing (respectively pre-composing) with the ``dot'' map $B_{n-1}\rightarrow \one$ (respectively $\one\rightarrow B_{n-1}$).  This implies that $d_{CB}\circ d_{BA}\simeq 0$ and $d_{DC}\circ d_{CB}\simeq 0$, as claimed.
\end{proof}

\begin{definition}
Let $\b_{(n-1,1)}: P_{n-1}Y_n(2n)\ip{2-2n}\rightarrow P_{n-1}X_n$ be the chain map whose mapping cone is the complex $Q_n$ constructed in Lemma \ref{lemma:Qexists}.
\end{definition}

\begin{proposition}\label{prop:hookEigenmap}
We have $Q_n=\Cone(\b_{(n-1,1)})\in \ker(P_n)$. 
\end{proposition}
\begin{proof}
We must show that $Q_n B_k\simeq 0$ for $1\leq k\leq n-1$.  Recall that by construction $Q_n$ is the mapping cone on a map $P_{n-1}Y_n\rightarrow P_{n-1}X_n$ (shifts omitted).  Tensoring on the right with $B_k$ ($2\leq k\leq n-1$) yields
\[
P_{n-1}Y_nB_k \cong P_{n-1}B_{k-1}Y_n\simeq 0
\]
since $B_k$ slides past $Y_n$ and is annihilated by $P_{n-1}$.  Similarly $P_{n-1}X_nB_k \simeq 0$ for $2\leq k\leq n-1$.  This shows that $Q_nB_k\simeq 0$ for $2\leq k\leq n-1$.

The argument that $Q_nB_1\simeq 0$ depends on $n$.  The case $n=1$ is trivial, and the case $n=2$ is taken care of explicitly by Example \ref{example-baseCase}.  If $n\geq 3$, then $P_{n-1}B_1\simeq 0$.  Thus, the first and fourth terms of $Q_n B_1$ (with respect to diagram (\ref{eq:Qdiagram1})) are null-homotopic, and we see that
\[
Q_nB_1 \simeq \Cone(\d_2-\d_1)B_1,
\]
up to shift.  Before going through the remaining details, we remark that $\Cone(\d_2-\d_1)$ should be thought of as a ``perturbed'' version of $\Cone(\d_2)$.  On the other hand, from the definition of $\d_2$ (Definition \ref{def:deltaMaps}) we see that
\[
\Cone(\d_2) \simeq P_{n-1}B_{n-1}\Cone(\b_{(n-2,1)}) = P_{n-1}B_{n-1}Q_{n-1}.
\]
This complex is annihilated by $(-)\otimes B_1$ up to homotopy, since $Q_{n-1}\in \im P_{n-1}$.  Alas, this does not prove that $\Cone(\d_2-\d_1)B_1\simeq 0$.  But it does give a clue as to how to proceed.

\textbf{Case $n\geq 4$.}    Since $\b_{(n-1)}$ factors through $P_{n-2}$ (Remark \ref{rmk:factorThroughR}), it follows from the definition that $\d_1\in\Homg(P_{n-1}B_{n-1}Y_{n-1},P_{n-1}B_{n-1}X_{n-1})$ factors through $P_{n-1}B_{n-1}$, up to homotopy. Thus, $\d_1B_1$ is null-homotopic, since it factors through $P_{n-1}B_{n-1}B_1 \cong P_{n-1}B_1B_{n-1}\simeq 0$ (we are using $n\geq 4$ in sliding $B_1$ past $B_{n-1}$).  Thus $(\d_2-\d_1)B_1\simeq \d_2 B_1$, and
\[
\Cone(\d_2-\d_1)B_1\cong \Cone((\d_1-\d_2)B_1)\cong \Cone(\d_2B_1)\cong \Cone(\d_2)B_1\simeq 0
\]
by the above arguments.

\textbf{Case $n=3$}.
This case reduces to an explicit computation.  We have $Q_3 B_1\simeq \Cone(\d_2-\d_1)B_1$ up to shifts, where $\Cone(\d_2-\d_1)$ can be written (omitting the homological shifts) as in the following diagram:
\[
\left(\begin{minipage}{5in}
\begin{tikzpicture}
\node (a) at (0,0){$P_2B_2(4)$ };
\node (b) at (3.5,0){ $P_2B_2B_1(3)$};
\node (c) at (7,0){ $P_2B_2B_1(1)$ };
\node (d) at (10.5,0){ $P_2B_2$ };
\path[->,>=stealth',shorten >=1pt,auto,node distance=1.8cm,font=\small]
(a) edge node[auto] {} (b)
(b) edge node {} (c)
(c) edge node {} (d)
(a) edge[bend left=30] node {$-u_2\otimes \Id_{B_2}$} (d);
\end{tikzpicture}
\end{minipage}\right)
\] 
Note that, disregarding the long arrow, the resulting complex is just $\Cone(\d_2)=P_2B_2Q_2$, which becomes contractible on applying $(-)\otimes B_1$.  The long arrow is the contribution from $\d_1=u_2\b_{(2)}$.  Now, $\Cone(\d_2-\d_1)B_1$ is the complex
\[
 \left(\begin{minipage}{5.8in}
\begin{tikzpicture}[baseline=1em]
\node (a) at (.2,0){$P_2B_2B_1(4)$ };
\node (b) at (4,0){ $P_2B_2B_1(4)\oplus P_2B_2B_1(2)$};
\node (c) at (9,0){ $P_2B_2B_1(2)\oplus P_2B_2B_1(0)$ };
\node (d) at (13,0){ $P_2B_2B_1(0)$ };
\path[->,>=stealth',shorten >=1pt,auto,node distance=1.8cm,font=\small]
(a) edge node[auto] {} (b)
(b) edge node {} (c)
(c) edge node {} (d)
(a) edge[bend left=30] node {$u_2B_2B_1$} (d);
\end{tikzpicture}
\end{minipage}\right).
\]
The length 1 components of the differential are precisely the same as in $Q_2B_1$, tensored with the identity map of $P_2B_2$.  The long arrow does not interfere with the Gaussian eliminations which realize the contractibility of $Q_2B_1$, and so the above chain complex is contractible.  We leave the details to the reader.
\end{proof}

\begin{proof}[Proof of Theorems \ref{thm-eigenconesAreOrtho} and \ref{thm-endP}]
This is an induction on $n$.  First, the maps $\b_{(n)}$ are constructed in Definition \ref{def:easyEigenmap}. Theorem \ref{thm-endP} and theorem \ref{thm-eigenconesAreOrtho} hold for $P_1$ trivially, and Theorem \ref{thm-eigenconesAreOrtho} holds for $P_2$ (the maps $\b_{(2)}$ and $\b_{(1,1)}$ are constructed in Example \ref{example-baseCase}.  Assume by induction that Theorems \ref{thm-endP} and \ref{thm-eigenconesAreOrtho} hold for $P_1,\ldots,P_{n-1}$.  Then $\b_{(n)}$ is constructed above, which proves Theorem \ref{thm-eigenconesAreOrtho} for $P_n$ (see Proposition \ref{prop:hookEigenmap}; this holds assuming the induction hypotheses).  Then Proposition \ref{prop:conditionalEndRings} computes the algebra $\DEnd(P_n)$, which proves Theorem \ref{thm-endP} for $P_n$.  This concludes the inductive step in our proof of Theorem \ref{thm-endP} and Theorem \ref{thm-eigenconesAreOrtho}.
\end{proof}

\subsection{Additional properties of $Q_n$}
\label{subsec:QnProps}

The reader should compare the results here with similar results in \cite{H14a}.

\begin{definition}\label{def:Qn}
Let $Q_n:=\Cone(\b_{(n-1,1)})\in \KC^-(\SBim_n)$.
\end{definition} 

\begin{proposition}
Let $U_k^{(n)}\in \Endg(P_n)$ denote chain maps which represent the classes $u_k$ in homology, with respsect to the isomorphism in Theorem \ref{thm-endP}.  Then
\begin{equation}\label{eq:ConeUk}
\Cone(U_n^{(n)})\simeq Q_n \ \ \ \ \ \ \ \ \ \ \ \  \ \ \Cone(U_k^{(n)})\simeq (Q_k\sqcup \one_{n-k})\otimes P_n
\end{equation}
\end{proposition}
\begin{proof}
From Proposition \ref{prop:periodicityMap} we know that  $U^{(n)}_k\simeq \rho_{nk}(U_k^{(k)})$.  From Proposition \ref{prop:rhoProp} we have
\[
(U_k^{(k)}\sqcup \Id_{\one_{n-k}}) \otimes \Id_{P_n} \simeq (\Id_{P_k}\sqcup \Id_{\one_{n-k}})\otimes U_k^{(n)}.
\]
Taking mapping cones gives $(Q_k\sqcup \one_{n-k})\otimes P_n \simeq (P_k\sqcup \one_{n-k})\otimes Q_n$.  But since $Q_n\in \im(P_n)$ by construction, we have $P_nQ_n\simeq Q_n$, hence
\[
(P_k\sqcup \one_{n-k}) Q_n \simeq (P_k\sqcup \one_{n-1})P_n Q_n \simeq P_nQ_n \simeq Q_n.
\]
\end{proof}

Now, recall the various symmetries of $\SBim_n$.  Let $\tau:\SBim_n\rightarrow \SBim_n$ be the covariant functor such that $\tau(B_i)=B_{n-i}$ and $\tau(MN)=\tau(M)\tau(N)$ for all $M,N\in \SBim_n$.  Let $\omega: \SBim_n\rightarrow \SBim_n$ be the covariant functor such that $\omega(B_i)=B_i$ and $\omega(MN)=\omega(N)\omega(M)$ for all $M,N\in \SBim_n$.  Note that $\tau$ comes from the Dynkin automorphism of $A_{n-1}$, and $\omega$ reflects all diagrams across a vertical axis in terms of the Elias-Khovanov diagrammatics.

Note that
\[
\tau((A\sqcup \one_{n-k}) B ) \cong (\one_{n-k}\sqcup \tau(A))\tau(B),
\]
\[
\omega((A\sqcup \one_{n-k}) B ) \cong \omega(B)(\omega(A)\sqcup \one_{n-k}) ,
\]
and
\[
\tau\circ \omega((A\sqcup \one_{n-k}) B ) \cong \tau\circ\omega(B)(\one_{n-k}\sqcup\tau\circ\omega(A)). ,
\]
Note that $\tau(P_n)\simeq P_n \simeq \omega(P_n)$ by the uniqueness statement in Definition/Theorem \ref{defthm-Pn}.  Thus, $\tau$ and $\omega$ define algebra automorphisms of $\DEnd(P_n)\cong \Z[u_1,\ldots,u_n,\xi_1,\ldots,\xi_n]$. The degrees are:
\[
\deg(u_i)  = (2i,0,2-2i) \ \ \ \ \ \ \  \ \  \ \ \ \deg(\xi_i)=(2i-4,1,2-2i)
\]
We collapse gradings by introducing $\deg_s=\deg_q+\deg_t+4\deg_a$.  It is easy to see that the classes $u_i$ and $\xi_i$ ($1\leq i\leq n$) span the subgroup in homology consisting of classes with $\deg_s=2$, and all other homogeneous classes other than $1$ have $\deg_s>2$.  The classes $u_i$ and $\xi$ are further distinguished from one another by $\deg_a$ and $\deg_t$.  Thus, each $u_i$ and $\xi_i$ spans its corresponding homology group, which is isomorphic to $\Z$.  Since these classes are unique up to unit scalar, any automorphism of $\DEnd(P_n)$ must send $u_i\mapsto \pm u_i$ and $\xi_i\mapsto \pm\xi_i$.  As an easy corollary we obtain the following.

\begin{proposition}\label{prop:Qsymmetries}
The object $Q_n$ has the symmetries of a rectangle: $\tau(Q_n)\simeq Q_n$ and $\omega(Q_n)\simeq Q_n$.   Furthermore,
\begin{equation}\label{eq:Qhop}
(Q_k\sqcup \one_{n-k}) Q_n \simeq  Q_n (Q_k\sqcup \one_{n-k}) \simeq (\one_{n-k}\sqcup Q_k) Q_n\simeq  Q_n(\one_{n-k}\sqcup Q_k) 
\end{equation}
\end{proposition}
\begin{proof}
We have $\tau(Q_n)\simeq \Cone(\tau(u_n))\cong \Cone(\pm u_n)\cong Q_n$, since mapping cones satisfy $\Cone(\pm f)\cong \Cone(f)$.  This proves the first statement.  Now, let us prove that
\begin{equation}\label{eq:QP}
(Q_k\sqcup \one_{n-k}) P_n \simeq   (\one_{n-k}\sqcup Q_k)P_n.
\end{equation}
Tensoring on the right with $Q_n\simeq P_nQ_n$ then gives the equivalence $(Q_k\sqcup \one_{n-k}) Q_n \simeq   (\one_{n-k}\sqcup Q_k)Q_n$.  Given this, the remaining equivalences in (\ref{eq:Qhop}) follow by applications of $\tau$ and $\omega$.  Thus, it suffices to prove (\ref{eq:QP}).

Proposition \ref{prop:periodicityMap} says that $(Q_k\sqcup \one_{n-k})P_n$ is homotopy equivalent to the mapping cone of $U_k\in \Endg(P_n)$.  Applying the functor $\tau$ gives
\[
\Cone(\tau(U_k))\simeq (\one_{n-k}\sqcup Q_k) P_n
\]
since $\tau(Q_k)\simeq Q_k$ and $\tau(P_n)\simeq P_n$.  On the other hand, $U_k$ generates the corresponding homology group of $\DEnd(P_n)$, which is isomorphic to $\Z$.  Thus $\tau(U_k)\simeq \pm U_k$, and $\Cone(\tau(U_k))\cong \Cone(U_k)$.  This proves (\ref{eq:QP}), and completes the proof.
\end{proof}

It is possible to show that $Q_n$ is quasi-idempotent: $Q_n^{\otimes 2}\simeq Q_n\oplus Q_n(2n)\ip{1-2n}$.  We won't prove this, since it will follow from later work of the author's and Ben Elias \cite{ElHog17a-pp}.  We collect some remaing properties together for later convenience.
\begin{proposition}
The chain complex $Q_n$ satisfies
\begin{enumerate}\setlength{\itemsep}{3pt}
\item $F(\b)\otimes Q_n\simeq Q_n\simeq Q_n\otimes F(\b)$ for all braids $\b\in \Br_n$.
\item $P_{n-1} J_n \simeq (Q_n\rightarrow P_{n-1}(2n)\ip{2-2n})$, i.e. there exists a distinguished triangle
\[
P_{n-1}(2n)\ip{2-2n}\rightarrow P_{n-1} J_n \rightarrow Q_n \rightarrow P_{n-1}(2n)\ip{1-2n}.
\]
where $J_n=F(\sigma_{n-1}\cdots \sigma_2\sigma_1^2\sigma_2\cdots\sigma_{n-1})$ is the Rouquier complex associated to the Jucys-Murphy braid (Definition \ref{def-JMandCycles}).
\item $\Tr(Q_n)\simeq P_{n-1}\oplus P_{n-1}(2n-4,1)\ip{2-2n}$.
\end{enumerate}
\end{proposition}
\begin{proof}
Statement (1) follows from the fact that $Q_n\simeq P_n\otimes Q_n\simeq Q_n\otimes P_n$, and $P_n$ absorbs Rouquier complexes.

By definition $Q_n$ fits into a distinguished triangle of the form
\begin{equation}\label{eq:TrQ}
P_{n-1}Y_n(2n)\ip{2-2n} \rightarrow P_{n-1}X_n \rightarrow Q_n\rightarrow P_{n-1}Y_n(2n)\ip{1-2n}.
\end{equation}
Tensoring on the right with $Y_n\inv$ gives the distinguished triangle from (2), given that $Q_nY_n\inv\simeq Q_n$.

For (3), apply $\Tr$ to the distinguished triangle (\ref{eq:TrQ}) and use (\ref{eq:TracePieces}) to obtain a distinguished triangle
\[
P_{n-1}(2n-4,1)\ip{3-2n} \rightarrow P_{n-1} \rightarrow \Tr(Q_n)\rightarrow P_{n-1}(2n-4,1)\ip{2-2n}.
\]
Just as in the proof of Proposition \ref{prop:TrP}, the first map is null-homotopic for degree reasons. which forces $\Tr(Q_n)$ to split as claimed.
\end{proof}

Finally we give an interesting alternate description of $Q_n$.

\begin{proposition}\label{prop:4TermCx}
We have
\begin{equation}\label{eq:Qdiagram2}
Q_n=\left(\begin{minipage}{4.3in}
\begin{tikzpicture}
\node (a) at (0,0){$P$ };
\node (b) at (3,0){ $PBP$};
\node (c) at (7.3,0){ $PBP$ };
\node (d) at (10.3,0){ $P$ };
\path[->,>=stealth',shorten >=1pt,auto,node distance=1.8cm,font=\small]
(a) edge node[auto] {} (b)
(b) edge node {$u_{n-1}BP - PBu_{n-1}$} (c)
(c) edge node {} (d)
(a) edge[bend left=23] node {} (d);
\draw[frontline,->,>=stealth',shorten >=1pt,auto,node distance=1.8cm]
(a) to [bend right] node[below] {} (c);
\draw[frontline,->,>=stealth',shorten >=1pt,auto,node distance=1.8cm]
(b) to [bend right] node[below] {} (d);
\end{tikzpicture}
\end{minipage}\right)
\end{equation}
where $P=P_{n-1}$, $B=B_{n-1}$, and the degree shifts are the same as in (\ref{eq:Qdiagram1}).  The first and third maps are constructed from $P\simeq PP$ and the ``dot'' maps.
\end{proposition}
\begin{proof}
Tensor (\ref{eq:Qdiagram1}) on the right with $P$ and simplify.
\end{proof}
Note that in case $n=2$ we recover the expression for $Q_2$ from Example \ref{example-baseCase}.

\printbibliography

\end{document}